\documentclass{elsart1p}
\usepackage[french]{babel}
\usepackage[T1]{fontenc}
\usepackage{amsmath}
\usepackage{amssymb}
\usepackage[dvips]{graphicx}
\usepackage{float}
\usepackage{bbold}
\usepackage{cancel}
\usepackage{dsfont}
\usepackage{amsfonts}

\newtheorem{dfn}[thm]{\textsf{Définition}}

\newcommand{\N}{\mathds{N}}
\newcommand{\R}{\mathds{R}}

\newcommand{\C}{\mathds{C}}

\begin{document}
\begin{frontmatter}

\title{Solutions globales pour l'équation de Schrödinger cubique en dimension 3}
\author[South]{A.Poiret},
\ead{aurelien.poiret@math.u-psud.fr}
\address[South]{Faculté des sciences d'Orsay, Département de mathématiques, Bâtiment 430 Bureau 106, France}

\begin{abstract} 
Dans cet article, on s'intéresse à l'équation de Schrödinger cubique en dimension 3. Grâce à une randomisation des données initiales, comme l'on fait N.Burq et N.Tzvetkov, on construit des ensembles de solutions globales pour des données qui sont dans $ L^2(\mathds{R}^3) $, alors que le problème est $ H^{1/2}( \mathds{R}^3  ) $ critique. Le point clef de la démonstration est l'existence d'une estimée bilinéaire de type Bourgain pour l'oscillateur harmonique.
\end{abstract}


\end{frontmatter}

\textbf{Remerciement} Je tiens à remercier grandement Nicolas Burq sans qui cet article n'aurait jamais vu le jour. Que ce soit pour toutes les idées qu'ils m'aient donné (en particulier la preuve de l'estimée bilinéaire de la section 3 est sienne) ou pour tous ses conseils, je lui en suis grandement reconnaissant.

\newpage

Dans cet article, on explique comment construire des solutions globales pour certaines équations de Schrodinger sur-critiques. On considère les équations de Schrödinger suivantes :
 \begin{equation} \label{schrodinger} \tag{NLS}
  \left\{
      \begin{aligned}
        i \frac{ \partial \tilde{ u }  }{ \partial t } +  \Delta  \tilde{ u }   & =   K | \tilde{u} |^{p-1} \tilde{ u }
       \\ \tilde{u}(0,x)& =u_0(x),
      \end{aligned}
    \right.
\end{equation}
où $ K \in \lbrace -1,1 \rbrace $ et p désigne un entier impair. Alors, on peut trouver dans \cite{tao} les théorèmes suivants : il existe un espace $ X^s_T $ continûment injecté dans $ C^0( [-T,T], H^s(\R^d)) $ vérifiant les faits suivants :
\\\\ - Si $ p < 1 + \frac{4}{d-2s} $ (cas sous-critique) alors pour tout $ R >0$, il existe  $ T_R> 0 $ tel que si $ ||u_0||_{H^s(\R^d)} \leq R $, alors il existe une unique solution locale $ u \in X^s_T $ à l'équation (\ref{schrodinger}). De plus, l'application $ u_0 \in B_{H^s}(0,R) \rightarrow u \in X^s_{  T  } $ est continue. Cela signifie que le problème est localement bien posé. 
\\\\Si T peut être choisi égal à $ \infty $, on dit que le problème est globalement bien posé. C'est le cas par exemple si $ s=0 $ ou si $ s=1 $ et $ K=1$. De manière générale, si $ K =-1 $ ou si s est très proche de 0, le problème de globalisation est délicat.
\\\\Dans le cas où le problème est globalement bien posé, il est naturel d'étudier le comportement de la solution en $ + \infty $. 
\\Par exemple, pour $ p=3, \ d=3 $ et $ s=1 $, nous pouvons obtenir pour tout $  u_0 \in H^1(\R^d) $, l'existence de $ u_+ \in H^1(\R^d) $ tel que $ \underset{t \rightarrow \infty }{\lim} || u(t) - e^{it \Delta} u_+  ||_{  H^s( \R^d ) } = 0 $. Lorsque cette propriété est vérifiée, on parle de "scattering". Il existe des situations où le problème est globalement bien posé mais où le scattering n'est pas vérifié.
 \\\\- Si $ p = 1 + \frac{4}{d-2s} $ (cas critique), on peut démontrer l'existence locale d'une unique solution pour chaque donnée initiale comme dans le cas sous-critique, mais le temps d'existence de la solution ne dépend pas uniquement de $ ||u_0||_{H^s(\R^d)} $ et le problème de globalisation est un problème complexe.
\\\\- Si $ p > 1 + \frac{4}{d-2s} $ (cas sur-critique), il existe une suite de réels $ t_n \in \R $ et une suite de fonctions $ u_n \in H^s(\R^d) $ telles que $ \underset{n \rightarrow \infty }{\lim} ||u_n||_{H^s(\R^d)} = 0 $ et $ \underset{n \rightarrow \infty }{\lim} ||u_n(t_n,.)||_{H^s(\R^d)} = 0 $ où $ u_n(t,.) $ désigne une solution de l'équation (\ref{schrodinger}) avec donnée initiale correspondante. Par conséquent, le flot de l'équation (\ref{schrodinger}) n'est pas continue en 0, le problème n'est pas bien posé et les méthodes usuelles ne permettent pas l'étude de l'équation dans cette situation.
\\\\Dans la majeure partie de cet article, on s'intéresse au cas $ p = 3 $, $ d = 3 $ qui est donc $ H^{1/2}( \mathds{R}^3)$ critique et nous allons construire un nombre infini non dénombrable de solutions globales pour des données initiales qui sont moralement dans $ L^2( \mathds{R}^3)$. Dans la section 9, on expliquera comment généraliser ce résultat (en dimension supérieure à 2) pour toutes équations de moins d'une demi dérivée sur-critique.
\\\\Pour établir ce résultat, on utilise les idées de N.Burq, L.Thomann et N.Tzvetkov développées dans  \cite{burq4} et \cite{burq6} en rendant la donnée initiale aléatoire dans l'espoir de gagner des dérivées dans un certain espace $ L^p $ pour $ p \neq 2  $.
\\\\Tout d'abord, en section 4, on vérifie que la donnée initiale rendue aléatoire ne permet pas de gagner de dérivées dans $ L^2 $ et que le problème reste sur-critique avec le même indice. On prouve également que notre donnée aléatoire est grande dans $ H^s $ pour $ 0 < s \ll 1 $. Ainsi, comme le théorème 2 de \cite{galla3}, nos solutions construites seront également valables pour des équations sur-critiques avec données initiales grandes.
\\\\ La preuve du théorème repose essentiellement sur deux résultats intermédiaires :
\\-la transformation de lentille (utilisée en dimension 1 dans \cite{burq6} et que nous pouvons généraliser en toutes dimensions) qui nous permet de se ramener à prouver un théorème d'existence locale sur $ ] - \frac{\pi}{4} ; \frac{\pi}{4}  [ $ pour l'équation de Schrödinger avec potentiel harmonique (voir section 2.3).
\\-l'existence d'une estimée bilinéaire pour l'oscillateur harmonique (voir section 3) qui permet de gagner la demi dérivée sur les termes d'ordre 1 en $ u_0 $. Cette estimée est l'analogue de l'estimée bilinéaire pour le Laplacian prouvée par Bourgain (voir \cite{staffilani} pour avoir une preuve).
\\\\En sections 5 et 6, par analogie au théorème 1 de \cite{galla2} ou le théorème 2 de \cite{galla1}, en utilisant les deux résultats intermédiaires précédents, on peut prouver que si la donnée initiale est assez petite pour une certaine norme alors nous pouvons appliquer un théorème de point fixe et obtenir une solution globale pour l'équation (\ref{schrodinger}).
\\\\Pour conclure, il suffit de vérifier que la probabilité que la norme de la donnée initiale soit petite n'est pas nulle. En utilisant les méthodes de N.Burq et N.Tzevtkov développées dans les paragraphes 3 et 4 de \cite{burq4}, en ayant choisi comme fonctions propres pour l'oscillateur harmonique les fonctions tenseurs pour avoir la meilleure décroissance possible, on peut évaluer la régularité de la donnée initiale  (section 7) et prouver le théorème (section 8).
\section{Introduction et notations}
En dimension d'espace d quelconque, on pose $ H = -\Delta + x^2 $ l'oscillateur harmonique. On note $ \lambda_n ^2 $ les valeurs propres et  $  h_n $ les fonctions propres de H que l'on indexe par $ n \in \mathds{N} $. On a donc
\begin{equation*}
H h_n = \lambda_n^2 h_n , \ \forall n \in \N.
\end{equation*}
On note $ W^{s,p}(\R^d) $ et $ H^s(\R^d) $ les espaces de Sobolev usuels. Puis, on définit les espaces de Sobolev harmoniques.
\begin{dfn} L'espace $ \overline{H}^s( \R^d ) $ est défini comme la fermeture de l'espace de Schwartz pour la norme
\begin{equation*}
|| u || _{  \overline{H} ^s(  \R^d )   } = || H^{s/2} u ||_ {L^2(   \R^d ) }.
\end{equation*}
\end{dfn}
\begin{dfn}
De manière similaire, l'espace $ \overline{W}^{s,p}( \R^d ) $ est défini comme la fermeture de l'espace de Schwartz pour la norme
\begin{equation*}
|| u || _{  \overline{W} ^{s,p}(  \R^d )   } = || H^{s/2} u ||_ {L^p(   \R^d ) }.
\end{equation*}
\end{dfn}
Dans \cite{sobolev}, nous pouvons trouver la proposition suivante :
\begin{prop} \label{comparaison}
Pour tout $ 1 < p < \infty $, $ s \geq 0 $, il existe une constante $ C > 0 $ telle que 
\begin{equation*} 
\frac{1}{C} || u || _{  \overline{W}^{s,p}( \R^d ) } \leq  || \nabla^s u ||_{L^p(\R^d)} + || <x>^s u ||_{L^p(\R^d ) }  \leq C || u || _{  \overline{W}^{s,p}( \R ^d ) }.
\end{equation*}
\end{prop}
Pour simplifier les calculs, on suppose $ p=3$ et $ d=3 $ dans l'équation (\ref{schrodinger}). Néanmoins, les résultats peuvent être obtenus de manière similaire pour $ d \geq 2 $ car la preuve est basée sur l'estimée bilinéaire que nous démontrons pour $ d  \geq 2 $. 
\\\\On suppose que les fonctions propres de l'oscillateur harmonique sont données par les fonctions propres tenseurs, c'est à dire
\begin{equation} \label{tanseur}
 h_n (x,y,z) = h_{n_1} (x)h_{n_2} (y) h_{n_3} (z) 
\end{equation}
où $  ( - \partial^2_x +x^2 )  h_{n_1}  = \lambda_{n_1}^2 h_{n_1} $ ,  $ ( - \partial^2_y +y^2 ) h_{n_2} = \lambda_{n_2}^2 h_{n_2} $, $( - \partial^2_z +z^2 ) h_{n_3} = \lambda_{n_3}^2 h_{n_3} $ 
\\et $ \lambda_{n_1}^2+\lambda_{n_2}^2 +\lambda_{n_3}^2= \lambda_n ^2 $.
\\\\Ensuite, supposons donnée $ (\Omega , A , P)  $ un espace de probabilité et $ (g_n ( \omega ) )_{n \in \N} $ une suite de variables aléatoires indépendantes et identiquement distribuées telle que $ E( |g_n|^2) < \infty $.
\\Soit $ u_0 \in  \overline{H}^\sigma(\R^3) $, c'est à dire
\begin{equation*}
u_0(x) = \sum_{n \in \N } c_n h_n(x)     \ \mbox{où} \  \sum_{n \in \N } \lambda_n^{2\sigma} |c_n|^2 < \infty,
\end{equation*}
et considérons l'application $ \omega \longrightarrow u_0( \omega,.)  $ de $ ( \Omega , A  ) $ dans $ \overline{H}^\sigma ( \R^3  )  $  que l'on équipe de sa tribu borélienne, définie par $ u_0(\omega,x) = \displaystyle{ \sum_{n \in \N} } c_n g_n( \omega ) h_n(x) $.
\\ On peut vérifier que l'application $ \omega \longrightarrow u_0(\omega,.) $ est dans $ L^2 ( \Omega , \overline{H}^\sigma( \R^3 ) )  $ et on définit $ \mu $  comme la loi de cette variable aléatoire. Il est alors possible d'appliquer le théorème de transfert suivant :
\begin{equation}
P( \omega \in \Omega / \Psi ( u_0(\omega,.) ) \in A ) = \mu ( f \in \overline{H}^\sigma ( \R^3) / \Psi (u_0) \in A   ),
\end{equation}
pour toute application $ \Psi: (\overline{H}^\sigma(\R^3) , \mathcal{B}(  \overline{H}^\sigma(\R^3) )  )\rightarrow (\R, \mathcal{B}(\R) ) $ mesurable et ensemble $ A \in \mathcal{B}(\R) $.
\\\\Premièrement, supposons que $ (g_n ( \omega ) )_{n \in \N} $ est une suite de variables indépendantes gaussiennes complexes standards (c'est à dire que la densité de $ g_n $ est donnée par $  \frac{1}{\pi} e^{-|z|^2}dL$  où $ dL $ désigne la mesure de Lebesgue sur $ \C $) alors nous avons les deux théorèmes suivants :
\\\begin{thm} \label{thm1}
Soient $ \sigma \in ]0, \frac{1}{2} [ $ avec $ u_0 \in \overline{ H }  ^{\sigma}  ( \mathds{R}^3 ) $ et $ s \in ] \frac{1}{2} ,  \frac{1}{2} +\sigma [ $, alors il existe un ensemble $ \Omega ' \subset \Omega $ vérifiant les conditions suivantes :
\\i) $ P( \Omega ' ) > 0$.
\\ii) Pour tout élément $ \omega \in \Omega'$, il existe une unique solution globale $ \tilde{u} $ à l'équation (\ref{schrodinger}) dans l'espace $ e^{it\Delta} u_0(\omega,.) + X^s $  avec donnée initiale $ u_0(\omega,.) $.
\\ iii) Pour tout élément $ \omega \in \Omega'$, il existe $ L ^+$  et $ L_ - \in \overline{H}^s(\mathds{R}^3) $ telles que 
\begin{align*}
\lim_{t \rightarrow \infty } || \tilde{u}(t) - e^{it\Delta} u_0 ( \omega,.) - e^{it\Delta} L^+ ||_{  H^s(\mathds{R}^3) } = 0,
\\ \lim_{t \rightarrow -\infty } || \tilde{u}(t) - e^{it\Delta} u_0 ( \omega,.) - e^{it\Delta} L_- ||_{  H^s(\mathds{R}^3) } = 0.
\end{align*}
De plus, pour tous $ \sigma' > \sigma $ et $ \alpha \in ]0,1] $, il existe deux constantes $ R $ et $ C > 0 $ telles que 
\\- si $  u_0 \notin \overline{H}^{\sigma'}(\mathds{R}^3) $ alors $  P( \omega \in \Omega / u_0(\omega,.) \in H^{\sigma'}(\mathds{R}^3) ) = 0 $,
\\- si $ u_0 \in \overline{H}^{\sigma'}(\mathds{R}^3) $ avec $ ||u_0||_{ \overline{H}^{\sigma'} (\mathds{R}^3) } \geq R $ alors 
\begin{equation*}
P( \Omega') \geq (1-\alpha)  \times  e^{ -C  ||u_0||^{ \frac{6}{\sigma'-\sigma}   }_{ \overline{H}^{\sigma'} (\mathds{R}^3) } \log ||u_0||_{ \overline{H}^{\sigma'} (\mathds{R}^3) }   }.
\end{equation*}
\end{thm}
\begin{thm} \label{thm2} Si $ \frac{1}{2} >  \sigma > 0 $ et $ u_0 \in \overline{ H }  ^{\sigma}  ( \mathds{R}^3 ) $ alors
\begin{equation} \label{proba1}
\lim _{  \eta \rightarrow 0 } \  \mu \left(  u_0 \in \overline{H}^\sigma {(\mathds{R}^3)} / \mbox{ on ait existence globale et scattering }  | \ ||u_0||_{\overline{H}^\sigma {(\mathds{R}^3)}   }   \leq \eta  \right)  = 1.
\end{equation}
De plus, si $ (u^1_0, u^2_0 )  \in \overline{H}^\sigma {(\mathds{R}^3)} \times \overline{H}^\sigma {(\mathds{R}^3)} $ alors pour tout $ \epsilon > 0 $,
\begin{equation}\label{proba2}
\begin{array}{rll}
\underset{\eta  \rightarrow 0}{\lim} & \  \mu_1 \otimes \mu_2  \bigg(  (u^1_0, u^2_0 )  \in \overline{H}^\sigma {(\mathds{R}^3)}  \times \overline{H}^\sigma {(\mathds{R}^3)} / \tilde{u}_1 \mbox{ et }   \tilde{u}_2 
 \\ & \mbox{ sont globales avec } || \tilde{u}_1-  \tilde{u}_2||_{X^0} \leq \epsilon \ |  \ ||u^1_0||_{\overline{H}^\sigma {(\mathds{R}^3)}   }   \leq \eta , \ ||u^2_0||_{\overline{H}^\sigma {(\mathds{R}^3)}   }   \leq \eta  \bigg) = 1,
\end{array}
\end{equation}
où $ \mu_i $ correspond à la loi de la variable aléatoire $ \omega \rightarrow u_0^i(\omega,.) $.
\end{thm}
Deuxièmement, on suppose que $ g_n  = \epsilon * b_n  $ où $ (b_n(\omega))_{n \in \N} $ est une suite de variables aléatoires Bernoulli standards (l'espérance est nulle) et $ \epsilon > 0 $ alors nous avons le théorème suivant :
\\\begin{thm} \label{thm3}
Soient $ \sigma \in ]0, \frac{1}{2} [ $ avec $ u_0 \in \overline{ H }  ^{\sigma}  ( \mathds{R}^3 ) $ et $ s \in ]  \frac{1}{2} ,  \frac{1}{2} +\sigma [ $ alors pour tout $ \alpha \in ]0,1] $ il existe une constante $ \epsilon > 0 $ et un ensemble $ \Omega_{\epsilon,\alpha} $ vérifiant les conditions suivantes :
\\i) $ P( \Omega_{ \epsilon,\alpha} ) \geq 1 - \alpha $.
\\ii) Pour tout élément $ \omega \in \Omega_{\epsilon,\alpha}$, il existe une unique solution globale $ \tilde{u} $ à l'équation (\ref{schrodinger}) dans l'espace $ e^{it\Delta} u_0(\omega,.) + X^s $ avec donnée initiale $ u_0(\omega,.) $.
\\iii) Pour tout élément $ \omega \in \Omega_{\epsilon,\alpha} $, il existe $ L ^+ $ et $ L_ - \in \overline{H}^s(\mathds{R}^3) $ telles que
\begin{align*}
\lim_{t \rightarrow \infty } || \tilde{u}(t) - e^{it\Delta} u_0 ( \omega,.) - e^{it\Delta} L^+ ||_{  H^s(\mathds{R}^3) } = 0,
\\ \lim_{t \rightarrow -\infty } || \tilde{u}(t) - e^{it\Delta} u_0 ( \omega,.) - e^{it\Delta} L_- ||_{  H^s(\mathds{R}^3) } = 0.
\end{align*}
De plus, on a la relation $ \alpha = C_1 e ^{  -\frac{C_2}{\epsilon ^2 ||u_0||^2_{  \overline{H}^\sigma(\mathds{R}^3)  } }  } $ et si $ u_0 \notin \overline{H}^s(\mathds{R}^3) $ alors
\begin{equation*}
P( \omega \in \Omega / u_0(\omega,.) \in H^s(\mathds{R}^3) ) = 0.
\end{equation*}
\end{thm}
\section{Résultats préliminaires}
Dans cette section, excepté dans la quatrième partie, on suppose la dimension d'espace d quelconque.
\subsection{Les estimées de Strichartz pour l'oscillateur harmonique}
Dans ce premier paragraphe, on établit les estimées de Strichartz pour l'oscillateur harmonique.
\begin{dfn}
Soient $ (q , r) \in [2, \infty ]^2 $ alors on dit que $ ( q , r ) $ est admissible si et seulement si
\begin{align*}
 ( q, r, d ) \neq ( 2, \infty , 2 ) \ \mbox{ et } \ \frac{2}{q} = \frac{d}{2} - \frac{d}{r}.
\end{align*}
\end{dfn}
\begin{dfn} Pour $ s \in \R $ et $ T \geq 0$, on définit
\begin{equation*}
    \overline{X}_T^s = \underset{ (q,r) \ admissible }{ \bigcap } L^q( [-T,T] , \overline{W}^{s,r}(\R^d)).
\end{equation*}
\end{dfn}
\begin{prop} \label{Stricharz} Pour tout temps $ T \geq 0 $, il existe une constante $ C_T > 0 $ telle que pour toute fonction $ u \in \overline{H}^s(\R^d) $,
\begin{equation*}
 || e^{-itH} u  ||_{\overline{X}_T^s} \leq C_T || u ||_{  \overline{H}^s(\R^d) }.
\end{equation*}
\end{prop}
\textit{Preuve.} Quitte à remplacer $ u $ par $ e^{iTH} u $, il suffit de prouver l'estimation pour un certain $ T> 0 $, par exemple pour $  T = \frac{\pi}{4}  $. De même, quitte à remplacer $ u $ par $ H ^{  \frac{s}{2} } u $, on peut limiter la preuve au cas où $ s= 0 $.
\\\\On a 
\begin{equation*}
e^{-itH} u ( t ,x ) = \left( \frac{1}{ \cos 2t } \right) ^{d/2} \times e^{it\Delta} u \left( \frac{1}{2}  \tan 2 t ,  \frac{x}{\cos 2 t }  \right) \times e^{  - \frac{i x^2 \tan 2 t }{2} }.
\end{equation*}
Ainsi, si le couple $ (q,r) $ est admissible, on obtient  
\begin{align*}
& || e^{-itH} u  ||_{L^q( ]  - \frac{\pi}{4} , \frac{\pi}{4}  [  , L^r (\R^d))} 
\\  = & \bigg| \bigg| \left( \frac{1}{ \cos 2t } \right) ^{d/2} \times e^{it\Delta} u \left( \frac{1}{2}  \tan 2 t ,  \frac{x}{\cos 2 t }  \right) \times e^{  - \frac{i x^2 \tan 2 t }{2} }  \bigg| \bigg|_{L^q( ]  - \frac{\pi}{4} , \frac{\pi}{4}  [  , L^r (\R^d))}
\\ = & || e^{it\Delta} u  ||_{L^q( \R , L^r(\R^d))}.
\end{align*}
Puis, nous pouvons utiliser les estimées de Strichartz pour le laplacien pour conclure. \hfill $ \boxtimes $
\begin{prop} \label{Stricharz2} Pour tout temps $ T \geq 0 $, il existe une constante $ C_T > 0 $ telle que pour tout couple $ ( q, r ) $ admissible, réel s et fonction $ F \in L^{q'}( [T,T], \overline{W}^{s,r'} (\R^d)) $,
\begin{equation*}
 \bigg| \bigg| \int _0^t e^{-i(t-s)H} F(s) ds   \bigg| \bigg| _{  \overline{X}^s_T} \leq C || F ||_{  L^{q'}( [-T,T], \overline{W}^{s,r'} (\R^d)) }.
\end{equation*}
\end{prop}
\textit{Preuve.} Il s'agit de la même preuve que celle des estimées de Strichartz pour le laplacien que l'on peut trouver dans \cite{tao}. En utilisant la proposition \ref{Stricharz}, par dualité, on obtient
\begin{equation*}
\bigg| \bigg| \int_{-T} ^ T  e^{isH} F(s) \bigg| \bigg| _{  \overline{H}^s(\R^d) } \leq C || F ||_{  L^{q'}( [-T,T], \overline{W}^{s,r'} (\R^d)) }.
\end{equation*}
Et finalement, on établit que
\begin{equation*}
 \bigg| \bigg| \int _{-T} ^ T  e^{-i(t-s)H} F(s) ds   \bigg| \bigg| _{  \overline{X}^s_T} \leq C || F ||_{  L^{q'}( [-T,T], \overline{W}^{s,r'} (\R^d)) },
\end{equation*}
puis nous pouvons conclure en utilisant le lemme de Christ-Kiselev. \hfill $ \boxtimes $
\subsection{Quelques propriétés des espaces de Bourgain}
Dans ce second paragraphe, on définit les espaces de Bourgain puis on établit leurs différentes propriétés.
\begin{dfn}
On définit l'espace $ \overline{X}^{s,b} = \overline{X}^{s,b}( \R * \R^d ) $ comme le complété de $  C_0^\infty ( \R * \R^d   ) $ pour la norme
\begin{align*}
|| u||^2_{   \overline{X}^{s,b} } & = \sum_{n \in \N } || \ <  t + \lambda_n ^2  >^b  \lambda_n^s  \widehat{P_n u } (t)  ||^2_{ L_t^2( \R , L^2_x(\R^d) ) }
\\ & = \sum_{n \in \N } || H^{s/2} e^{itH} P_n u(t,.) ||^2_{  L_x^2(\R^d,H_t^b( \R )) },
\end{align*}
où $ \widehat{P_n u } (t)   $ désigne la transformée de Fourier de $ P_n u := <  u , h_n  >_{L^2(\R^d) \times  L^2(\R^d)} * \ h_n $ par rapport à la variable temps.
\end{dfn}
\textit{Remarque:} Au vu de cette définition, il est important de noter que 
\begin{equation*}
 || H^{s/2} e^{itH} u(t,.) ||_{  L_x^2(\R^d,H_t^b( \R )) } \leq || u||_{   \overline{X}^{s,b} } .
\end{equation*}
Dans \cite{tao}, corollaire 2.10, nous pouvons trouver la proposition suivante :
\begin{prop} \label{bourgain1}
Pour tout temps $ T \geq 0 $ et entier $ b > \frac{1}{2} $, il existe une constante $ C_T > 0 $ telle que pour tout entier $ s \in \R $ et pour tout couple admissible $ (q,r )$, si $ u \in \overline{X}^{s,b} $ alors $ u \in L^q([-T,T] ,  \overline{W}^{s,r} (\R^d)) $ et
\begin{equation*}
||u||_{L^q([-T,T] ,  \overline{W}^{s,r} (\R^d))} \leq C_T ||u||_{  \overline{X}^{s,b} }.
\end{equation*}
\end{prop}
Dans les propositions \ref{bourgain2} et \ref{bourgain3}, quitte à remplacer $ u $ par $ H^{s/2}u $, nous pouvons limiter la preuve au cas $ s =0$.
\\Grâce au lemme 2.4 de  \cite{burq1}, nous pouvons établir la proposition suivante :
\begin{prop} \label{bourgain2}
Pour tout $ \theta \in [0,1] $, si $ b > \frac{1-\theta}{2} $ alors il existe une constante $ C > 0 $ telle que pour tout entier s et toute fonction $ u \in \overline{X}^{s,b} $,
\begin{equation*}
||u||_{ L^{ \frac{2}{\theta} } (\R,    \overline{H}^s (\R^d))} \leq C ||u||_{\overline{X}^{s,b}}.
\end{equation*}
\end{prop}
\textit{Preuve.} \`A l'aide de la transformée de Fourier inverse, on a 
\begin{equation*}
P_n u(t) = \frac{1}{2 \pi} \int_\R \frac{< \tau + \lambda_n ^2>^b}{<\tau+ \lambda_n^2>^b} \times e^{it \tau } \times \widehat{P_n u } ( \tau ) \ d \tau. 
\end{equation*}
Puis, pour $ b > \frac{1}{2} $, on obtient par l'inégalité de Cauchy-Schwarz que
\begin{equation*}
|P_n u (t)| \leq C \left( \int_\R <\tau+\lambda_n^2> ^{2b} | \widehat{P_n u }(\tau) |^2 \ d \tau \right) ^{1/2}.
\end{equation*}
Ainsi, en élevant au carré, en intégrant sur $ \R^d $ puis en sommant pour $ n \in \N $, on obtient pour tout réel $ b > \frac{1}{2} $ et fonction $ u \in \overline{X}^{0,b} $,
\begin{equation*}
||u||_{L^\infty(\R,L^2(\R^d))} \leq C ||u||_{  \overline{X}^{0,b} }.
\end{equation*}
Mais, par définition,
\begin{equation*}
||u||_{L^2(\R,L^2(\R^d))} =  ||u||_{  \overline{X}^{0,0} }
\end{equation*}
donc le résultat suit par interpolation. \hfill $ \boxtimes $
\begin{prop} \label{bourgain3} Pour toute constante $ 1 > \delta > 0 $, il existe deux constantes $ b' < \frac{1}{2} $ et $ C > 0 $ telles que pour tout entier $ s \in \R $ et toute fonction $ u \in L^{1+\delta}(\R,  \overline{H}^s (\R^d)) $,
\begin{equation*}
||u||_{\overline{X}^{s,-b'}} \leq C ||u||_{L^{1+\delta}(\R,  \overline{H}^s (\R^d))}.
\end{equation*}
\end{prop}
\textit{Proof.} D'après la proposition \ref{bourgain2}, on a par dualité que pour tout $ \theta  \in [0,1] $ et $ b > \frac{1-\theta}{2} $, il existe $ C> 0 $ telle que pour toute fonction $ u \in L^{ \frac{2}{2-\theta} } ( \R , L^2( \R^d)) $
\begin{align*}
||u||_{ \overline{X}^{0,-b} } \leq C ||u||_{L^{ \frac{2}{2-\theta} } ( \R , L^2( \R^d))}.
\end{align*}
Puis on choisit $ \theta = \frac{2 \delta}{1+ \delta} $ et $ b= \frac{1-\theta+\delta}{2} < \frac{1}{2} $ pour obtenir la proposition. \hfill $ \boxtimes  $
\\\\Dans \cite{tao}, lemme 2.11, nous pouvons trouver la proposition suivante :
\begin{prop} \label{bourgain4}
Soit $ \psi \in C^\infty_0 ( \R ) $ alors pour tout $ b \geq 0 $, il existe une constante $ C > 0 $ telle que pour tout $ s \in \R $ et toute fonction $ u \in \overline{X}^{s,b} $,
\begin{equation*}
|| \psi(t) u  ||_{\overline{X}^{s,b}} \leq C ||u||_{\overline{X}^{s,b} }
\end{equation*}
\end{prop}
Enfin, on donne une dernière proposition dont la preuve peut être trouvée au lemme 3.2 de \cite{ginibre} (en prenant $b'=1-b$ et $T=1$).
\begin{prop} \label{bourgain5}
Soit $ \psi \in C^\infty_0 ( \R ) $ alors pour tout $ 1 \geq b > \frac{1}{2} $, il existe une constante $ C > 0 $ telle que pour tout $ s \in \R $ et toute fonction $ F \in \overline{X}^{s,b-1}$,
\begin{equation*}
\bigg| \bigg| \psi(t) \int_0^t e^{-i(t-s) H } F(s) \ ds   \bigg| \bigg|_{\overline{X}^{s,b}} \leq C ||F||_{\overline{X}^{s,b-1} }.
\end{equation*}
\end{prop}
\`A partir de maintenant, on pose $ T = \frac{\pi}{4} $, puis on définit le nouvel espace de Bourgain qui va nous intéresser.
\begin{dfn}
On définit l'espace $ \overline{X}^{s,b}_T =  \overline{X}^{s,b}([ - T ; T ] * \R^d ) $ comme le sous ensemble de $ \overline{X}^{s,b} $ pour lequel la norme suivante
\begin{equation*} 
|| u || _{  \overline{X}^{s,b}_T  }  = \inf_{  w \in \overline{X}^{s,b} }  \left\{  || w || _{\overline{X}^{s,b} }  \  avec \  w |_{[-T,T]} = u \right\}
\end{equation*}
est finie.
\end{dfn}
La proposition \ref{bourgain2} nous permet d'obtenir le résultat suivant :
\begin{prop} \label{bourgain7}
Soient $ b > \frac{1}{2} $ et $ s \in \R $ alors $ \overline{X}^{s,b}_T \hookrightarrow C^0( [-T,T], \overline{H}^s(\R^d) ) $.
\end{prop}
\subsection{La transformation de Lentille}
Dans ce troisième paragraphe, on définit la transformation de lentille puis on établit ses propriétés fondamentales.
\begin{dfn}
Pour $ u(t,x)$ une fonction mesurable de $  ] - \frac{\pi}{4} ; \frac{\pi}{4}  [  \times \R^d $, on définit pour $ t \in \R $ et $ x \in \R ^d $, la fonction $ \tilde{u}(t,x) $ de la façon suivante:
\begin{equation*}
\tilde{u}(t,x) =  \left( \frac{1}{\sqrt{1+4t^2}}  \right) ^{d/2} \times u \left( \frac{1}{2} \arctan(2t) , \frac{x}{\sqrt{1+4t^2} } \right)   \times e^{ \frac{ix^2t}{1+4t^2}  }.
\end{equation*}
\end{dfn}
Dans \cite{tao2}, on peut trouver la proposition suivante :
\begin{prop} \label{lechangementdevariable} Soit $ K \in \R $ alors, 
\begin{center} \begin{tabular}{lll}
& u est solution de & $ i \frac{ \partial u }{ \partial t } -H u  = K \cos(2t)^{  \frac{d}{2}(p-1) -2 } | u|^{p-1} u $
\\ si et seulement si & $ \tilde{u} $ est solution de & $ i \frac{ \partial \tilde{u} }{ \partial t } + \Delta \tilde{u}   = K | \tilde{u}|^{p-1} \tilde{u} $ .
\end{tabular}
\end{center}
\end{prop}
\begin{dfn} Pour $ s \in \R $, on définit
\begin{equation*}
   X^s = \underset{ (q,r) \ admissible }{ \bigcap } L^q( \R , W^{s,r}(\R^d)).
\end{equation*}
\end{dfn}
\begin{prop} \label{toutvabienbis}
Soit $ s \geq 0$ alors il existe une constante $ C > 0 $ telle que pour tous $ p \in [1,+\infty ] $ et $ q \in [1,+\infty ]  $ vérifiant $ \frac{2}{p} + \frac{d}{q} - \frac{d}{2} \leq 0   $, on a pour toute fonction $ u \in L^p( [-T,T] , \overline{W}^{s,q} ( \R^d ) ) $,
\begin{align*}
|| \tilde{u} ||_{ L^p( \R , W^{s,q} ( \R^d )  )  } \leq C || u || _{L^p( [-T , T] , \overline{W}^{s,q} ( \R^d ) )   }.
\end{align*}
\end{prop}
\textit{Preuve.} Par interpolation, il suffit de prouver le résultat pour $ s =n \in \N $. Soit $ \alpha \in \N^d$ avec $ | \alpha | \leq n $, alors grâce à la formule de Leibniz, on obtient
\begin{align*}
& \partial ^\alpha_x \tilde{u}(t,x) 
\\ = \ & \partial ^\alpha_x \left( \left( \frac{1}{\sqrt{1+4t^2}}  \right) ^{d/2} \times u \left( \frac{1}{2} \arctan(2t) , \frac{x}{\sqrt{1+4t^2} } \right)   \times e^{ \frac{ix^2t}{1+4t^2}  } \right)
\\ = \ & \left( \frac{1}{\sqrt{1+4t^2}}  \right) ^{d/2} \times \sum_{ 0 \leq \beta \leq \alpha} \binom{\alpha}{\beta} \partial ^\beta_x \left( u \left( \frac{1}{2} \arctan(2t) , \frac{x}{\sqrt{1+4t^2} } \right) \right) \times \partial ^{\alpha-\beta}_x \left( e^{ \frac{ix^2t}{1+4t^2}  }\right).
\end{align*}
Puis, comme 
\begin{equation*}
 | \ \partial ^{\alpha-\beta}_x ( e^{ \frac{ix^2t}{1+4t^2}  } ) | \leq  C_{\alpha,\beta} \left( 1+ |  \frac{x}{\sqrt{1+4t^2}  } |^{|\alpha-\beta|}\right)  ,
\end{equation*}
on établit
\begin{align*}
 | \ \partial ^\alpha_x \tilde{u}(t,x) | \leq   \sum_{ 0 \leq \beta \leq \alpha} C_{\alpha,\beta} & \left( \frac{1}{\sqrt{1+4t^2}}  \right) ^{d/2 + |\beta|} 
 \\ & \times | \partial ^\beta_x  u | \left( \frac{1}{2} \arctan(2t) , \frac{x}{\sqrt{1+4t^2} } \right)  \times \left( 1+ |  \frac{x}{\sqrt{1+4t^2}  } |^{|\alpha-\beta|}\right)  .
\end{align*}
Par conséquent, 
\begin{align*}
& \hspace*{0.3cm}|| \partial ^\alpha_x  \tilde{u} ||_{ L^p( \R , L^{q} ( \R^d )  )  } 
\\  \leq \ & \sum_{ 0 \leq \beta \leq \alpha} C_{\alpha,\beta} || \left( \frac{1}{\sqrt{1+4t^2}}  \right) ^{d/2+d/q + |\beta|}  \times || u  \left( \frac{1}{2} \arctan(2t) , .  \right) ||_{  \overline{W}^{|\alpha|,q}(\R^d) } ||_{L^p(\R)}
\\ \leq \ & \sum_{ 0 \leq \beta \leq \alpha} C_{\alpha,\beta}  || \left( \frac{1}{\sqrt{1+ \tan^2 2t }}  \right) ^{d/2+d/q-2/p + |\beta|}  \times u (t,x)  ||_{ L^p( [-T,T ],  \overline{W}^{|\alpha|,q}(\R^d) ) } .
\end{align*}
Pour conclure, il suffit de remarquer que $ \frac{d}{2} + \frac{d}{q}- \frac{2}{p} + |\beta| \geq 0 $ pour $ \beta \in \N^d $. \hfill $ \boxtimes $
\\\\Finalement, on a établi la proposition suivante :
\begin{prop} \label{toutvabientris}
Soient $ s \geq 0  $ et $ u \in \overline{X}^{s}_T $ alors $ \tilde{u} \in X^{s} $ et il existe une constante $ C> 0 $ telle que pour tout $ u \in \overline{X}^{s}_T $,
\begin{equation*}
||\tilde{u}||_{X^{s}}  \leq C ||u||_{\overline{X}^{s}_T}.
\end{equation*}
\end{prop}
Et enfin, grâce à la proposition \ref{bourgain1}, on obtient la proposition suivante :
\begin{prop} \label{toutvabien} Soient $ b > \frac{1}{2}$, $ s \geq 0 $ et $ u \in \overline{X}^{s,b}_T $ alors $ \tilde{u} \in X^{s} $ et il existe une constante $ c> 0 $ telle que pour tout $ u \in \overline{X}^{s,b}_T $,
\begin{equation*}
||\tilde{u}||_{X^{s}}  \leq c ||u||_{\overline{X}^{s,b}_T}.
\end{equation*}
\end{prop}
\subsection{Propriétés basiques des fonctions propres de l'oscillateur harmonique}
Dans ce quatrième paragraphe, on donne quelques estimations classiques des fonctions propres tenseurs de l'oscillateur harmonique.
\begin{prop}  \label{dispersive}
Pour tout $ \delta > 0 $, il existe une constante $ C_\delta > 0 $ telle que pour tous $ n,m,k \in \N^3 $,
\begin{align}
 & || h_n ||_{L^4(\R^3)}  \leq C \lambda_n^{-1/4} (\log \lambda_n ) ^3 \label{propre1},
 \\ & || h_n ||_{L^\infty(\R^3)}  \leq C \lambda_n^{-1/6} \label{propre2},
\\ & || h_n h_m||_{L^2(\R^3)}  \leq C_\delta  \max ( \lambda_n, \lambda_m)^{-1/2+\delta} \label{propre3},
\\ & || h_{n} h_{m} h_{k} ||_{L^2(\R^3)}  \leq C_\delta  \max ( \lambda_n, \lambda_m, \lambda_k )^{-1/2+\delta} \label{propre4}.
\end{align}
\end{prop}
\textit{Preuve.} Les estimations (\ref{propre1}) et (\ref{propre2}) sont très connues en dimension 1 (voir \cite{koch} pour une démonstration).
\\Dans ce cadre en dimension 3, comme $ \lambda_n ^2 = \lambda_{n_1}^2 + \lambda_{n_2}^2 + \lambda_{n_3}^2 $ alors il existe $ i \in (1,2,3 ) $ tel que $ \lambda_{n_i } ^2 \geq \frac{\lambda_n^2 }{3} $.
\\Ainsi, on obtient 
\begin{align*}
|| h_n ||_{L^4(\R^3)}  & = || h_{n_1} ||_{L_x^4(\R)} || h_{n_2} ||_{L_y^4(\R)} || h_{n_3} ||_{L_z^4(\R)} 
\\ & \leq C \lambda_{n_1}^{-1/4} \lambda_{n_2}^{-1/4}  \lambda_{n_3}^{-1/4} \log(\lambda_{n_1}) \log(\lambda_{n_2})\log(\lambda_{n_3})
\\ & \leq C  \lambda_n^{-1/4} ( \log \lambda_n )^3,
\end{align*}
car $ \lambda_{n_1} \leq \lambda_n, \ \lambda_{n_2} \leq \lambda_n $ et $ \lambda_{n_3} \leq \lambda_n $.
\\Nous pouvons faire la même preuve pour l'estimation (\ref{propre2}). 
\\\\Ensuite, l'estimation (\ref{propre3}) est démontrée en dimension 1 dans \cite{burq6}.
\\On peut supposer que $ \max( \lambda_n,\lambda_m) = \lambda_n $ et $ \max( \lambda_{n_1},\lambda_{n_2}, \lambda_{n_3}) = \lambda_{n_1} $, alors 
\\$ \lambda^2_{n_1} \geq \frac{\lambda_n^2}{3}  \geq \frac{\lambda_m^2}{3} \geq \frac{\lambda_{m_1}^2}{3}. $ Puis, grâce à (\ref{propre1}), on obtient
\begin{align*}
|| h_n h_m ||_{L^2(\R^3)} & = || h_{n_1} h_{m_1} ||_{L_x^2(\R)} || h_{m_2} h_{m_2} ||_{L_y^2(\R)} || h_{n_3} h_{m_3} ||_{L_z^2(\R)} 
\\ & \leq  || h_{n_1} h_{m_1} ||_{L_x^2(\R)} || h_{n_2}||_{L^4_y(\R)} || h_{m_2} ||_{L_y^4(\R)} || h_{n_3} ||_{L^4_z(\R)} ||  h_{m_3} ||_{L^4_z(\R)}
\\ & \leq C_\delta \lambda_{n_1}^{-1/2+\delta} 
\\ & \leq C_\delta \lambda_n^{-1/2+\delta} .
\end{align*}
Pour l'estimation (\ref{propre4}), supposons que $ \max ( \lambda_n, \lambda_m, \lambda_k ) = \lambda_n $. Alors, en utilisant (\ref{propre2}) et (\ref{propre3}), on trouve
\begin{align*}
|| h_{n} h_{m} h_{k} ||_{L^2(\R^3)}  & \leq || h_{n} h_{m} ||_{L^2(\R^3)} || h_{k} ||_{L^\infty(\R^3)}
\\ &\leq C_\delta  \lambda_n ^{-1/2+\delta}
\\ &\leq C_\delta  \max ( \lambda_n, \lambda_m, \lambda_k )^{-1/2+\delta}.
\end{align*}
Ce qui démontre la proposition. \hfill $ \boxtimes $
\begin{lem} \label{inegalite sobolev} Pour tout $ s \in \R  $, il existe une constante $ C> 0 $ telle que pour toutes fonctions $ f,g,h \in \overline{H}^s(\R)$,
\begin{align*}
|| f(x) g(y) h(z)||_{\overline{H}^s(\R^3)} \leq C \times  ( & || f||_{\overline{H}^s (\R)}|| g||_{L^2(\R)}  || h||_{L^2(\R)} + || f ||_{L^2(\R)} || g||_{\overline{H}^s (\R)}  || h||_{L^2(\R)}
\\ & + || f ||_{L^2(\R)} || g||_{L^2(\R)} || h||_{\overline{H}^s (\R)}  ).
\end{align*}
\end{lem}
\textit{Preuve.} Il suffit d'établir le résultat dans le cas où $ f(x)=h_n(x), \ g(y)=h_m(y), \ h(z)=h_k(z) $.
\\Dans ce cas, on a 
\begin{align*}
& || f(x) g(y) h(z)||_{\overline{H}^s(\R^3)} 
\\ = \ & || \ H^{s/2} [ f(x) g(y) h(z)] \ ||_{L^2(\R^3)}
\\ = \ & || \ (  \lambda_n^2 + \lambda_m^2 + \lambda_k^2 )^{s/2}[ h_n(x) h_m(y) h_k(z) ] \ ||_{L^2(\R^3)}
\\ \leq \ &  \lambda_n^s ||  h_n(x) h_m(y) h_k(z) ||_{L^2(\R^3)} + \lambda_m^s || h_n(x) h_m(y) h_k(z)||_{L^2(\R^3)} \\& \hspace*{7cm} + \lambda_k^s || h_n(x) h_m(y) h_k(z)  ||_{L^2(\R^3)}
\\ \leq \ &  || h_n||_{\overline{H}^s (\R)}|| h_m ||_{L^2(\R)}  || h_k ||_{L^2(\R)} + || h_n ||_{L^2(\R)} || h_m ||_{\overline{H}^s (\R)}  || h_k ||_{L^2(\R^3)}
\\ & \hspace*{7cm} + || h_n ||_{L^2(\R)} || h_m ||_{L^2(\R)} || h_k ||_{\overline{H}^s (\R)}. & \boxtimes 
\end{align*}
\begin{prop} \label{propre6} Pour tout $ \delta > 0 $ et $ s \in [ 0, 1 ] $, il existe une constante $ C > 0 $ telle que pour tous $ n,m,k \in \N ^3 $,
\begin{align*}
 & || \  h_n h_m ||_{ \overline{H}^s(\R^3) } \leq C \times \max ( \lambda_n, \lambda_m )^{s-1/2 + \delta },
\\ & || \  h_n h_m h_k ||_{ \overline{H}^s(\R^3) } \leq C \times \max ( \lambda_n, \lambda_m , \lambda_k )^{s-1/2 + \delta }.
\end{align*}
\end{prop}
\textit{Preuve.} Grâce à (\ref{propre3}) et (\ref{propre4}), par interpolation, il suffit d'établir les inégalités pour $ s=1$.
\\Pour la première inégalité, on peut supposer que $ \max( \lambda_n,\lambda_m) = \lambda_n $ et $ \max( \lambda_{n_1},\lambda_{n_2}, \lambda_{n_3}) = \lambda_{n_1} $.
\\En utilisant le lemme \ref{inegalite sobolev}, (\ref{propre1}) et le lemme A.8 de \cite{burq6} avec $ \theta=1 $, on trouve
\begin{align*}
|| h_n h_m ||_{  \overline{H}^1 (\R^3)}  & \leq C  \times  ( ||  h_{n_1} h_{m_1} ||_{ \overline{H}^1 (\R) } + || h_{n_2} h_{m_2} ||_{\overline{H}^1 (\R)} + || h_{n_3}  h_{m_3} ] ||_{  \overline{H}^1 (\R)}  )
\\ & \leq  C \times ( \max(\lambda_{n_1},\lambda_{m_1})^{1/2+\delta}  +  \max(\lambda_{n_2},\lambda_{m_2})^{1/2+\delta} + \max(\lambda_{n_3},\lambda_{m_3})^{1/2+\delta} )
\\ & \leq C \times \max( \lambda_n, \lambda_m ) ^{1/2+\delta}. 
\end{align*}
Pour la seconde inégalité, on peut supposer que $ \max( \lambda_n,\lambda_m, \lambda_k) = \lambda_n $ alors, grâce à l'inégalité précédente, 
\begin{align*}
|| h_n h_m h_k ||_{  \overline{H}^1 (\R^3)}  & \leq || h_n h_m ||_{  \overline{H}^1 (\R^3)} || h_k ||_{  L^\infty (\R^3)} + || h_n h_m ||_{  L^2 (\R^3)} || h_k ||_{  W^{1,\infty }  (\R^3)}
\\ & \leq C \times \max( \lambda_n,\lambda_m)^{1/2+\delta} \lambda_k^{-1/6} + C \times \max( \lambda_n,\lambda_m)^{-1/2+\delta} \lambda_k^{5/6}
\\ & \leq C \times \max( \lambda_n, \lambda_m , \lambda_k ) ^{1/2+\delta}. & \boxtimes 
\end{align*}
\begin{prop} \label{rapidement} Soient $ \delta > 0 $, $ l \geq 4 $ et $ N \geq 1 $, alors il existe une constante $ C_N > 0 $ telle que si on suppose
\begin{equation*}
\lambda_{n_1} \geq \lambda_{n_2}^{1+\delta} \mbox{ et } \lambda_{n_2} \geq \lambda_{n_3} \geq ... \geq \lambda_{n_l} , \mbox{ cela implique que }  \bigg| \int _{\R^3}  \prod_{i=1}^l  h_{n_i}(x) dx  \bigg| \leq C_N \lambda_{n_1}^{-N}.
\end{equation*}
\end{prop}
\textit{Preuve.} En utilisant (\ref{propre1}) et (\ref{propre2}), on obtient
\begin{align*}
\bigg| \int _{\R^3}  \prod_{i=1}^l  h_{n_i}(x) dx  \bigg| & \leq  \lambda_{n_1}^{-2k} \times|| H^k (  \prod_{i=2}^l h_{n_i}) h_{n_1} ||_{L^1(\R^3)}
\\ & \leq  \lambda_{n_1}^{-2k } \times ||  \prod_{i=2}^l h_{n_i} ||_{H^{2k}(\R^3)}
\\ & \leq C_k \times \lambda_{n_1}^{-2k } \lambda_{n_2}^{2k} \times \prod_{i=2}^l || h_{n_i} ||_{L^{2(l-1)}(\R^3)} 
\\ & \leq C_k \times \left( \frac{\lambda_{n_2}}{\lambda_{n_1}} \right) ^{2k} 
\\ & \leq C_k \times \lambda_{n_1}^{ \frac{-k\delta}{1+\delta}   } , \ \forall k \in \N^*.  & \boxtimes
\end{align*}
Soit $ \eta \in C_0^\infty( \R ) $ telle que $  \eta( 0 ) = \eta (1 ) = 1  $ et $ \eta ( 2) = 0 $, alors on définit pour $ N= 2^k $ la suite d'opérateurs suivante:
\begin{equation*}
\Delta_N (u ) = \left\{
    \begin{array}{ll}
         (  \eta ( \frac{H}{N^2} ) -  \eta  ( \frac{4H}{N^2}  )  ) u & \ \mbox{pour} \ N \geq 1,   \\
         \ \ 0  & \ \mbox{sinon}.
    \end{array}
\right.
\end{equation*}
Remarquons que si $ \lambda_n \notin [ \frac{N}{2}  , 2N ] $ alors $ \Delta_N ( h_n ) = 0 $ et que $ \underset{N}{\sum} \  \Delta_N (u ) = u $.
\begin{lem} \label{cas facile} Il existe $ b' < \frac{1}{2} $ tel que pour tous $ \delta > 0 $ et $ K \geq 1 $, on ait l'existence d'une constante $ C_K>0 $ telle que si on suppose $ N_1 \geq N_2^{1+\delta} $ et $ N_2 \geq N_3 \geq N_4 $ alors pour tous $ u_1,  u_2, u_3, u_4 \ \in \overline{X}^{0,b'}$,
\begin{align*}
\hspace*{1cm} \bigg| \int_{\R * \R^3} \Delta_{N_1}(u_1) \Delta_{N_2}(u_2) \Delta_{N_3}(u_3) \Delta_{N_4}(u_4) \bigg|  \leq C_K  N_1^{-K} \prod_{i=1}^4 || \Delta_{N_i}(u_i)||_{ \overline{X}^{0,b'} }.
\end{align*}
\end{lem}
\textit{Preuve.} On commence par étudier le cas où $ u_i(t,x) = c_i (t) h_{n_i}(x) $. D'après les propositions \ref{rapidement} et \ref{bourgain2}, on trouve
\begin{align*}
& \bigg| \int_{\R * \R^3} \Delta_{N_1}(u_1) \Delta_{N_2}(u_2)\Delta_{N_3}(u_3) \Delta_{N_4}(u_4) \bigg|  
\\ = &  \bigg| \int_{\R * \R^3} \prod_{i=1}^4 \phi( \frac{ \lambda_{n_i}^2 }{N_i^2} ) c_i (t) h_{n_i}(x) \ dt \ dx \bigg|
\\ \leq & \prod_{i=1}^4 \  \phi( \frac{ \lambda_{n_i}^2 }{N_i^2} ) \times  \int_\R | c_1(t) ... c_4 (t) | \ dt \times \bigg| \int_{\R^3}  h_{n_1}(x)  ...  h_{n_4}(x) \ dx \bigg|
\\ \leq & C_K N_1^{-K} \times  \prod_{i=1}^4 \  \phi( \frac{ \lambda_{n_i}^2 }{N_i^2} )  \times  \prod_{i=1}^4 || c_i(.)||_{ L_t^4(\R) }
\\  \leq & C_K  N_1^{-K} \prod_{i=1}^4 || \Delta_{N_i}( u_i) ||_{ L^4( \R , L^2(\R^3)) }
\\ \leq & C_K  N_1^{-K} \prod_{i=1}^4 || \Delta_{N_i}(u_i)||_{ \overline{X}^{0,b'} }.
\end{align*}
Pour le cas général, on pose $ u_i (t,x) = \underset{ k \in \N }{\sum} c_{i,k}(t) h_k(x) $, alors 
\begin{align*}
& \bigg| \int_{\R * \R^3} \Delta_{N_1}(u_1) \Delta_{N_2}(u_2)\Delta_{N_3}(u_3) \Delta_{N_4}(u_4) \bigg|  
\\  \leq \ & \sum_{k_1,k_2,k_3,k_4} \bigg| \int_{\R * \R^3} \Delta_{N_1}(c_{1,k_1} h_{k_1}) \Delta_{N_2}(c_{2,k_2} h_{k_2})\Delta_{N_3}(c_{3,k_3} h_{k_3}) \Delta_{N_4}(c_{4,k_4} h_{k_4}) \bigg|
\\  \leq \ & C_K  N_1^{-K}   \sum_{k_1,k_2,k_3,k_4} \prod_{i=1}^4 || \Delta_{N_i}(c_{i,k_i} h_{k_i})||_{ \overline{X}^{0,b'} }
\\  \leq \ & C_K  N_1^{-K+12}  \sqrt{  \sum_{k_1,k_2,k_3,k_4} \prod_{i=1}^4 || \Delta_{N_i}(c_{i,k_i} h_{k_i})||^2_{ \overline{X}^{0,b'} } } .
\end{align*}
Or 
\begin{align*}
\sum_{k_i} || \Delta_{N_i}(c_{i,k_i} h_{k_i})||^2_{ \overline{X}^{0,b'} } & = \sum_{k_i} \sum_n || < \tau + \lambda_n > ^{b'} \widehat{P_n(c_{i,k_i} h_{k_i}  ) } ||^2_{ L_t^2( \R , L^2_x(\R^d) ) }
\\ & = \sum_n || < \tau + \lambda_n > ^{b'} \widehat{P_n(c_{i,n} h_n  ) } ||^2_{ L_\tau^2( \R , L^2_x(\R^d) ) }
\\ & = || u_i ||_{ \overline{X}^{0,b'} } . & \boxtimes 
\end{align*}
\begin{prop}  \label{propre5}
Pour tout $ s \geq 0 $, il existe deux constantes $ C_1 > 0 $ et $ C_2 > 0 $ telles que pour tout $ n \in \N $,
\begin{equation*}
 C_1 \lambda_n^s \leq  || \nabla^s h_n ||_{L^2(\R^3)} \leq C_2 \lambda_n ^s.
\end{equation*}
\end{prop}
\textit{Preuve.} En utilisant la proposition \ref{comparaison}, on a 
\begin{align*}
 C' \left( || \nabla^s h_n ||_{L^2(\R^3)} + ||<x>^s h_n  || _{L^2(\R^3)} \right) &  \leq \lambda_n^s =  || h_n ||_{ \overline{ H}^s(\R^3) } \\ & \leq C \left( || \nabla^s h_n ||_{L^2(\R^3)} + ||<x>^s h_n  || _{L^2(\R^3)} \right),
 \end{align*}
puis
 \begin{equation*}
 C' || \nabla^s h_n ||_{L^2(\R^3)} \leq \lambda_n^s \leq C \left( || \nabla^s h_n ||_{L^2(\R^3)} + ||  \nabla^s ( \hat{h_n}  ) || _{L^2(\R^3)} + || h_n ||_{L^2(\R^3)} \right).
\end{equation*}
Mais comme les fonctions propres sont les fonctions tenseurs alors $ |h_n(x)| = | \hat{h_n}(x)|  $, car cette égalité est vraie en dimension 1, et le résultat suit. \hfill $ \boxtimes $
\section{L'estimée bilinéaire pour l'oscillateur harmonique}
L'objectif de cette section est d'établir une estimée bilinéaire de type Bourgain pour l'oscillateur harmonique. On suppose la dimension d'espace $ d  \geq 2 $ et on propose de prouver le théorème suivant:
\begin{thm} \label{bilis} Pour tout $ \delta \in ]0 , \frac{1}{2} ] $, il existe une constante $ C > 0 $ telle que pour tous $ N,M, u $ et $ v $,
 \begin{align*}
  & || e^{it H} \Delta_N ( v ) \ e^{it  H} \Delta_M (u) ||_{L^2 ( [-1 ; 1 ] , L^2 ( \R^d ))  } 
  \\  \leq \ &  C \times \min  ( N,M)^{ \frac{d-2}{2} } \times  \left( \frac{ \min ( N,M)  }{ \max ( N,M )}  \right) ^{1/2-\delta } \times   || \Delta_N(v) ||_{L^2(\R^d)} || \Delta_M(u) ||_{L^2(\R^d)}.
\end{align*}
\end{thm}
On remarque que pour prouver le théorème \ref{bilis}, il suffit de prouver le théorème suivant:
\begin{thm} \label{bili}
Pour tout $  \delta \in ]0 , \frac{1}{2} ] $, il existe une constante $ C > 0 $ et un réel $ \epsilon > 0 $ tels que pour tous $ N,M, u $ et $ v $,
 \begin{align*}
 & || e^{it H} \Delta_N ( v ) \ e^{it  H} \Delta_M (u) ||_{L^2 ( [-\epsilon ; \epsilon ] , L^2 ( \R^d ))  } 
 \\ \leq \ & C \times  \min ( N,M)^{ \frac{d-2}{2} } \times  \left( \frac{ \min ( N,M)  }{ \max ( N,M )} \right) ^{1/2-\delta } \times || \Delta_N(v) ||_{L^2(\R^d)} || \Delta_M(u) ||_{L^2(\R^d)}.
\end{align*}
\end{thm}
En effet, on peut remplacer $ u $ par $ e^{i\epsilon H}u $ et $ v $ par $ e^{i\epsilon H} v $ pour obtenir 
\begin{align*}
& || e^{i(t+\epsilon) H} \Delta_N ( v ) \ e^{i(t+\epsilon)  H} \Delta_M (u) ||_{L^2 ( [-\epsilon ; \epsilon ] , L^2 ( \R^d ))  } 
\\ \leq \ & C \times \min ( N,M)^{ \frac{d-2}{2} }  \times \left( \frac{ \min ( N,M)  }{ \max ( N,M )} \right)^{1/2-\delta } \times || \Delta_N(v) ||_{L^2(\R^d)} || \Delta_M(u) ||_{L^2(\R^d)}.
\end{align*}
Puis, on utilise le changement de variable $ t \longleftrightarrow t + \epsilon $ et le théorème (\ref{bili}) pour trouver que 
\begin{align*}
& || e^{it H} \Delta_N ( v ) \  e^{i t  H} \Delta_M (u) ||_{L^2 ( [-\epsilon ; 2\epsilon ] , L^2 ( \R^d ))  } \\ \leq \ & C \times \min ( N,M)^{ \frac{d-2}{2} } \times \left( \frac{ \min ( N,M)  }{ \max ( N,M )} \right) ^{1/2-\delta } \times || \Delta_N(v) ||_{L^2(\R^d)} || \Delta_M(u) ||_{L^2(\R^d)}.
\end{align*}
On peut ainsi itérer le procédé $ 2 E(  \frac{1}{\epsilon} ) $ fois pour établir le théorème \ref{bilis} et on cherche donc à montrer le théorème \ref{bili}.
\\\\Soit $ r \ll 1 $ et $ \phi \in C_0^\infty(\R) $ qui vérifie
\begin{equation*}
\phi(x)  = \left\{
    \begin{array}{ll}
        1  \  \mbox{ pour } x \in [1/2;2], \\
         0 \ \mbox{ pour } x \in [0,1/2-r] \cup [2+r , \infty [,
    \end{array}
\right.
 \end{equation*}
et posons $ \Delta'_N = \phi( \frac{H}{N^2} ).  $ Alors en utilisant que $ \phi (x) *(  \eta(x) - \eta(4x)  ) = \eta(x) - \eta(4x)  $ pour tout $ x \in \R $, on a la proposition suivante:
\begin{prop} \label{ouf} Pour tout N, on a  
\begin{equation*}
\Delta'_N o \Delta_N = \Delta_N.
\end{equation*}
\end{prop}
Par conséquent, pour prouver le théorème \ref{bili}, il suffit de montrer le théorème suivant:
\begin{thm} \label{bilibis}
Pour tout $ \delta \in ]0 , \frac{1}{2} ] $, il existe une constante $ C > 0 $ et un réel $ \epsilon > 0 $ tels que pour tous $ N,M, u $ et $ v $,
\begin{align*} 
& || e^{it H} \Delta'_N(v) \ e^{it  H} \Delta'_M(u) ||_{L^2 ( [-\epsilon ; \epsilon ] , L^2 ( \R^d ))  } 
\\ \leq \ &  C  \min ( N,M)^{  \frac{d-2}{2}  } \times \left( \frac{  \min ( N,M)  }{ \max ( N,M )} \right) ^{1/2-\delta } \times  || v ||_{L^2(\R^d)} ||u ||_{L^2(\R^d)}.
\end{align*}
\end{thm}
En effet, si le théorème \ref{bilibis} est vérifié, nous pouvons appliquer cette inégalité à $ v $ remplacé par $ \Delta_N (v) $ et $ u $ remplacé par $ \Delta_M (u) $ puis nous pouvons utiliser la proposition \ref{ouf} pour obtenir le théorème \ref{bili}.
\\\\ \textbf{Cas $ M \sim N $ avec $ M \geq N $ :}
\\\\Pour $ d=2$, nous pouvons utiliser les inégalités de Strichartz (soit le théorème \ref{Stricharz}) pour trouver que
\begin{align*}
& || e^{it H} \Delta'_N(v) \ e^{it  H} \Delta'_M(u) ||_{L^2 ( [-\epsilon ; \epsilon ] , L^2 ( \R^d ))  } 
\\ \leq \ & || e^{it H} \Delta'_N(v)||_{L^4{( [-\epsilon ; \epsilon ] , L^4(\R^d))}} \times || e^{it H} \Delta'_N(v)||_{L^4{( [-\epsilon ; \epsilon ] , L^4(\R^d))} }
\\  \leq \ & || e^{it H} \Delta'_N(v)||_{L^4{( ]-\pi ; \pi [ , L^4(\R^d))} } \times || e^{it H} \Delta'_N(v)||_{L^4{( ]-\pi ; \pi [ , L^4(\R^d))} }
\\ \leq \ & C || \Delta'_N(v) ||_{L^2(\R^d)}\times || \Delta'_N(v)||_{L^2(\R^d )}
\\ \leq \ & C  || v ||_{L^2(\R^d)} \times  ||u ||_{L^2(\R^d)}.
\end{align*}
Pour $ d \geq 3 $, en utilisant encore le théorème \ref{Stricharz} et les injections de Sobolev, on établit que
\begin{align*}
& || e^{it H} \Delta'_N(v) \ e^{it  H} \Delta'_M(u) ||_{L^2 ( [-\epsilon ; \epsilon ] , L^2 ( \R^d ))  } 
\\ \leq \ & || e^{it H} \Delta'_N(v)||_{L^\infty{( [-\epsilon ; \epsilon ] , L^d(\R^d))}} \times || e^{it H} \Delta'_N(v)||_{L^2{( [-\epsilon ; \epsilon ] , L^{  \frac{2d}{d-2} }(\R^d))} }
\\ \leq \ & || e^{it H} \Delta'_N(v)||_{L^\infty{( ]-\pi ; \pi [ , \overline{W}^{ \frac{d-2}{2}   ,2}(\R^d))} } || e^{it H} \Delta'_N(v)||_{L^2{( ]-\pi ; \pi [ , L^\frac{2d}{d-2}  (\R^d))} }
\\  \leq \ & C || \Delta'_N(v) ||_{\overline{H}^{ \frac{d-2}{2}  }(\R^d)}\times || \Delta'_N(v)||_{L^2(\R^d )}
\\ \leq \ & C  N^{ \frac{d-2}{2}  } || v ||_{L^2(\R^d)} \times  ||u ||_{L^2(\R^d)}.
\end{align*}
Finalement, si nous posons $ u _M = \Delta'_M (u) $ et $ v_N = \Delta'_N (v)  $, on se ramène à démontrer que pour tout $ \delta \in ]0 , \frac{1}{2} ] $, il existe une constante $ C > 0 $ et un réel $ \epsilon > 0 $ tels que pour tous $ N,M, u $ et $ v $, si $ M > 10 N $ alors 
\begin{equation} \label{biblitris}
|| e^{it H} v_N  \ e^{it  H} u_M ||_{L^2 ( [-\epsilon , \epsilon ] , L^2 ( \R^d ))  } \leq C N^{ \frac{d-2}{2} } \left( \frac{  N  }{ M }  \right)^{1/2-\delta } || v ||_{L^2(\R^d)} ||u ||_{L^2(\R^d)}.
\end{equation}  
Pour démontrer ce dernier résultat, on commence par donner quelques notions élémentaires sur les opérateurs pseudo-différentiels. 
\subsection{Outils sur les opérateurs pseudo-différentiels et applications aux fonctions propres}
\begin{dfn}
Pour $ m \in \R$, on définit $ T^m $ comme l'espace vectoriel des symboles $ q(x,\xi) \in C^\infty( \R^d \times \R^d ) $ qui vérifient pour tous $ \alpha \in \N^d $ et  $ \beta \in \N^d $, l'existence d'une constante $ C_{\alpha, \beta} $ telle que pour tout $ ( x , \xi ) \in \R^d \times \R^d $, on ait
\begin{equation*}
  | \partial ^\alpha_x  \partial ^\beta_\xi q(x, \xi ) | \leq C_{\alpha, \beta}  (1+|x|+|\xi|)^{m-\beta}.
 \end{equation*} 
\end{dfn}
\begin{dfn}
Pour $ m \in \R$, on définit $ S^m $ comme l'espace vectoriel des symboles $ q(x,\xi) \in C^\infty( \R^d \times \R^d ) $ qui vérifient pour tous $ \alpha \in \N^d $ et  $ \beta \in \N^d $, l'existence d'une constante $ C_{\alpha, \beta} $ telle que pour tout $ ( x , \xi ) \in \R^d \times \R^d $, on ait
\begin{equation*}
  | \partial ^\alpha_x  \partial ^\beta_\xi q(x, \xi ) | \leq C_{\alpha, \beta}  (1+|\xi|)^{m-\beta}.
 \end{equation*} 
\end{dfn}
\begin{dfn}
Pour $ q  \ \in S^m \cup  T^m $ et $ h> 0  $, on pose $ Op_h(q) $ l'opérateur défini par
\begin{align*}
Op_h (q) f(x) &= (2 \pi h )^{-d} \int_{\R^d \times \R^d} e^{i(x-y) \xi / h } q(x , \xi ) f(y) \ dy d \xi
\\ & = (2 \pi )^{-d} \int_{\R^d \times \R^d} e^{ix \xi } q(x , h\xi ) \hat{f}(\xi) \ d \xi.
\end{align*}
\end{dfn}
Dans \cite{Martinez}, on peut alors trouver les deux théorèmes suivants:
\begin{thm}  \label{compose}
Soient $ q_1 \in S^{m1} $ (respectivement $ T^{m1} $) et $ q_2 \in S ^{m2} $ (respectivement $ T^{m2} $) alors il existe un symbole $ q \in S^{m_1+m_2} $ (respectivement $ T^{m_1+m_2} $) tel que
\begin{equation*}
 Op_h(q_1) \circ Op_h(q_2) = Op_h(q) 
 \end{equation*}
 avec
 \begin{equation*}
q = \sum _{| \alpha | \leq N } \frac{ h^{| \alpha | } } { i^{| \alpha| }} \ \partial^\alpha_\xi q_1 \partial^\alpha_x q_2 + h^{N+1} r_N   \ \mbox{où} \  r_N \in S^{m1+m2-(N+1) } \mbox{ (respectivement } T^{m1+m2-(N+1) }).
\end{equation*}
\end{thm} 
\begin{thm} \label{operateur}
Si $ q(x,\xi) \in S^0 $ alors pour tout $ s \in \R $, il existe une constante $ C > 0 $ telle que pour tout $ h \in ]0,1] $ et $ u \in H^s(\R^d) $,
\begin{equation*}
||   Op_h( q(x,\xi) ) u ||_{H^s(\R^d)} \leq C ||u||_{H^s(\R^d)}.
\end{equation*}
\end{thm}
On énonce ensuite la propriété suivante qui va permettre d'inverser l'oscillateur harmonique modulo un terme de reste très régularisant.
\begin{prop}  \label{inverse}
Soit $ \delta > 0 $ et définissons la fonction $ \eta \in C^\infty ( \R^d ) $ telle que  
\begin{equation*}
\eta (x) = \left\{
    \begin{array}{ll}
         0 \ \mbox{si} \ |x| \leq 1 + \delta,  \\
        1 \ \mbox{si} \ |x| \geq 1 + 2 \delta.
    \end{array}
\right.
 \end{equation*}
Posons $ p(x, \xi) = \xi^2 + x^2 -1 $ et définissons $ H_h = Op_h(p) \in Op_h( T^2 ) $ alors pour tout $ N \in \N^*$, il existe deux opérateurs pseudo-différentiels $  E_N \in Op_h( T^{-2} )   $ et $  R_N \in Op_h( T^{-(N+1)}) $ tels que 
\begin{equation*}
E_N \circ H_h = \eta + h^{N+1 } R_N.
\end{equation*}
\end{prop}
\textit{Preuve.} On pose
\begin{equation*}
e_0  =  \frac{\eta}{p} \in T^{-2}.
\end{equation*}
et pour $ n \geq 1 $, on définit $ e_n $ par récurrence de la façon suivante:
\begin{equation*}
e_n  = -   \frac{1}{p} \sum_{|\alpha|+j =n , j \neq n } \frac{1}{i^{|\alpha|}} \partial^\alpha_\xi e_j \partial^\alpha_x p \ \ \in T^{-2-n}.
\end{equation*}
Enfin, on pose
\begin{equation*}
 E_N = Op_h \left( \sum _{0 \leq j \leq N } h^j e_j \right).
\end{equation*}
Alors, par la proposition \ref{compose},
\begin{align*}
E _N \ o \ H_h & = Op_h \left( \sum_{0 \leq j \leq N } h^j e_j \right) \ o \ Op_h( p ) 
\\ & = Op_h \left( \sum_{ |\alpha| + j \leq N } \frac{h^{|\alpha|+j}}{i^{|\alpha|}} \ \partial^\alpha_\xi  e_{j}  \partial^\alpha_x p  + h^{N+1} r_N  \right)
\\ & = Op_h   \left( \sum_{ |\alpha| + j \leq N } \frac{h^{|\alpha|+j}}{i^{|\alpha|}} \ \partial^\alpha_\xi  e_{j}  \partial^\alpha_x  p  \right)  + h^{N+1} R_N  
\end{align*}
avec  $ R_N = Op_h(r_N) \in Op_h \left( T^{-(N+1)}  \right) .  $
\\\\Or 
\begin{align*}
\sum_{ |\alpha| + j \leq N } \frac{h^{|\alpha|+j}}{i^{|\alpha|}} \ \partial^\alpha_\xi  e_{j}  \partial^\alpha_x  p  & = e_0 p + \sum_{  1 \leq l   \leq N } \sum_{|\alpha|+j = l }\frac{h^{l}}{i^{|\alpha|}} \ \partial^\alpha_\xi  e_j  \partial^\alpha_x  p 
\\ & = \eta + \sum_{  1 \leq l   \leq N} h^l \left(  \sum_{|\alpha|+j = l , j \neq  l }\frac{1}{i^{|\alpha|}} \ \partial^\alpha_\xi  e_j  \partial^\alpha_x  p  +  e_l.p \right)
\\ & = \eta.
\end{align*}
\begin{flushright}
$ \boxtimes $
\end{flushright}
On peut donc maintenant établir une propriété fondamentale que vérifie les fonctions propres de l'oscillateur harmonique.
\begin{prop} \label{propre}Pour tous entiers $ K $ et $ N  $, pour $ c > 1 $ et $ 1 \leq p \leq \infty $, il existe une constante $ C > 0  $ telle que pour tout $ n \in \N $,
\begin{equation*}
|| <x >^{K} h_n ||_{L^p (|x| \geq c \lambda_n) } \leq C  \lambda_n ^
{-N}.
\end{equation*}
\end{prop}
\textit{Preuve.} Comme $ ( -\Delta + x^2 -\lambda_n^2) h_n = 0 $, en posant $ h =  \frac{1}{\lambda_n^2} $ et $ \Phi ( x) = h_n( \lambda_n x ) $ alors $ ( -h^2 \Delta + x^2 -1 )  \Phi = 0 $.
\\\\Soient $ \delta \ll 1 $ et $ \chi \in  C^\infty ( \R^d ) $ telle que
\begin{equation*}
\chi(x) = \left\{
    \begin{array}{ll}
        0 \ \mbox{si } \ |x| \leq 1, \\
        1 \ \mbox{si } \ |x| \geq 1 + \delta,
    \end{array}
\right.
\end{equation*}
\\et $ \overline{\chi} \in  C^\infty ( \R^d)  $ telle que
\begin{equation*}
\overline{\chi}(x) = \left\{
    \begin{array}{ll}
        0 \ \mbox{si } \ |x| \leq 1 + 2 \delta, \\
        1 \ \mbox{si} \ |x| \geq 1 + 3\delta.
    \end{array}
\right.
\end{equation*}
Alors
\begin{equation*}
H_h ( \chi \Phi ) = - h^2 \Delta  \chi  \Phi -2h^2 \ \nabla \chi . \nabla \Phi .
\end{equation*}
Puis, grâce à la proposition $ \ref{inverse} $, on trouve 
\begin{equation*}
 \eta \chi \Phi = - E_N ( h^2 \Delta \chi  \Phi + 2 h^2 \nabla \chi . \nabla \Phi )- h^{N+1} R_N ( \chi \Phi ) . \end{equation*} 
Et finalement,
\begin{equation*}
<x>^K \overline{\chi} \eta \chi \Phi = - <x>^K\overline{\chi} E_N ( h^2 \Delta \chi  \Phi+ 2 h^2 \nabla \chi . \nabla \Phi )- h^{N+1} <x>^K \overline{\chi} R_N ( \chi \Phi ) .
\end{equation*}
\\\textbf{ Estimation de $ <x>^K \overline{\chi} E_N (\Delta \chi  \Phi ) $:}
\\On a 
\begin{equation*}
\overline{\chi}(x) (E_N \Delta \chi  \Phi ) (x) = \frac{ \overline{ \chi (x) } }{ (2 \pi h )^d} \int_{\xi,1 \leq  |y| \leq 1+ \delta } e^{i(x-y) \xi / h } E_N(x, \xi ) (\Delta \chi \Phi )(y) \ dy d \xi. 
\end{equation*}
Comme $ |x| > 1 + 2 \delta $ alors $ |x-y | > \delta $. 
\\Puis comme $ \int_{|x-y| >  \delta } \leq \int_{|x_1-y_1| > \delta }+ \int_{|x_2-y_2| >  \delta }+ ... + \int_{|x_d-y_d| >  \delta } $, on peut se ramener à traiter le terme où $ |x_1-y_1 | >  \delta $.
\\\\De $ \frac{h ^ M }{(i(x_1-y_1))^M} \partial^M_{\xi_1}  e^{i(x-y) \xi / h  } = e^{i(x-y) \xi / h  }  $ et d'une intégration par parties, on déduit
\begin{align*}
  & \hspace*{4cm} <x>^K \overline{\chi}(x) (E_N \Delta \chi  \Phi ) (x)
  \\ & = ^+_-\frac{ \overline{ \chi  }(x) }{(2 \pi h )^d} \times  \int_{\xi,1 \leq  |y| \leq 1+ \delta } \frac{h^M}{i^M(x_1-y_1)^M} e^{i(x-y) \xi / h } <x>^K \partial ^M_{\xi_1} E_N(x, \xi ) (\Delta \chi \Phi )(y) \ dy d \xi.
  \end{align*}
Par conséquent, comme $ E_N \in T^{-2} $, on trouve
\begin{equation*}
|  < x > ^K  \overline{\chi}(x)   (E_N \Delta \chi  \Phi ) (x) | \leq C  \times h^{M-d} \times | \overline{\chi}(x) | \times \int_{\xi,1 \leq  |y| \leq 1+ \delta } \frac{ | ( \Delta \chi \Phi|) (y) }{(1+|x|+| \xi| )^{2+M-K}} \ d \xi dy.
\end{equation*}
On peut ensuite supposer que $ 2+M-K > 2d $ (ce qui est possible puisque dans le cas où $ 2+M-K \leq 2d $, on refait le même calcul avec $ M ' > M $ vérifiant $ 2+M'-K > 2d $ puis il suffira de majorer $ h^{M'} $ par $ h^M $) pour obtenir
\begin{equation*}
|  < x > ^K  \overline{\chi}(x)   (E_N \Delta \chi  \Phi ) (x) \times  | \leq C  \times h^{M-d}  \times || \Delta \chi  \Phi ||_{L^2(\R^d)} \times  \frac{| \overline{\chi}(x) | }{(1+|x|)^{2+M-K-d}}.
\end{equation*}
Et finalement, pour tous entiers M et K, il existe une constante $ C > 0  $ telle que pour tout $  h \in ]0,1] $,
\begin{equation*}
||< x > ^K \overline{\chi} E_N ( \Delta \chi  \Phi ) ||_{L^p(\R^d)} \leq C  \times h^{M-d}  ||   \Phi ||_{L^2(\R^d)}.
\end{equation*}
\\\\\textbf{ Estimation de $  <x>^K \overline{ \chi } E_N ( \nabla \chi . \nabla \Phi ) $:}
\\\\Comme $ \Phi $ satisfait $ -h^2 \Delta \Phi + x^2 \Phi = \Phi $ alors $ h || \nabla \Phi ||_{L^2(\R^d) } \leq || \Phi ||_{L^2(\R^d)}  $ puis nous pouvons procéder comme pour le premier terme pour obtenir le même genre d'estimation.
\\\\\textbf{Estimation de $ h^{N+1} <x>^K  \overline{\chi} R_N ( \chi \Phi ) $:}
\\On a 
\begin{align*}
h^{N+1} <x>^K \overline{\chi}(x)  R_N ( \chi \Phi )(x) 
\\ = h^{N+1} \times  \frac{1}{(2 \pi )^d} \times \int_{\R^d } e^{i(x-y) \xi  } <x>^K \overline{\chi}(x) r_N(x,h\xi)   \mathcal{F} (  \chi  \Phi )(\xi) \ d \xi  
\end{align*}
avec
\begin{equation*}
 <x>^K \overline{\chi}(x) r_N (x, \xi ) \in T^{-N-1+K}  \subset T^0 \subset S^0 \mbox{ pour } N \geq K-1.
 \end{equation*}
En utilisant le théorème \ref{operateur} et les injections de Sobolev, on trouve donc 
\begin{align*}
& || h^{N+1} <x>^K \overline{\chi}(x)  R_N ( \chi \Phi )(x) ||_{L^p(\R^d)}
\\ \leq \ & || h^{N+1} <x>^K \overline{\chi}(x)  R_N ( \chi \Phi )(x) ||_{H^{d/2+1}(\R^d)}
\\ \leq \ & C \times h^{N+1} \times  ||  \chi \Phi ||_{H^{d/2+1}(\R^d)}
\\ \leq \ & C  \times h^{N+1} \times ||  \Phi ||_{H^{d/2+1}(\R^d)}.
\end{align*}
Finalement, on obtient pour tous entiers $K$ et $ N  $, pour tout $ p \in [1,\infty[ $ et $ c > 1 $, l'existence d'une constante $ C> 0 $ telle que pour tout $ 0 < h \leq 1$, on ait 
\begin{equation*}
 || <x>^K \Phi ||_{L^p( |x| \geq c ) } \leq C  \times  h^N \times ||  \Phi ||_{H^{d/2+1}(\R^d)} .
\end{equation*}
Puis, en retournant à la variable initiale, on trouve que pour tous entiers $ K $ et $ N  $, pour tout $ p \in [1,\infty[ $ et $ c > 1 $, il existe une constante $ C > 0  $ telle que pour tout $ 0 < h \leq 1$ et $ n \in \N $ avec $ h = \frac{1}{\lambda_n^2} $, on ait
\begin{align*}
 || < \sqrt{h} x >^K h_n ||_{L^p( |x| \geq c \lambda_n ) } & \leq C  \times  h^{N-d/(2p)-d/4-1/2} \times ||  h_n ||_{H^{d/2+1}(\R^d)} 
\\ & \leq C  \times  h^{N-d/(2p)-d/4-1/2} \times ||  h_n ||_{\overline{H}^{d/2+1}(\R^d)} 
\\ & \leq C  \times  h^{N-d/(2p)-d/2-1} \times ||  h_n ||_{L^2(\R^d)}.
\end{align*}
Or
\begin{align*}
|| < x >^K h_n ||_{L^p( |x| \geq c \lambda_n ) } & \leq || h_n ||_{L^p( |x| \geq c \lambda_n ) } + h^{-K/2} || < \sqrt{h} x >^K h_n ||_{L^p( |x| \geq c \lambda_n ) } 
\\ & \leq C  \times  h^{N-d/(2p)-d/2-1-K/2} \times ||  h_n ||_{L^2(\R^d)}
\\ & \leq C  \times  \lambda_n^{-2N+d/p+d+K+2} \times ||  h_n ||_{L^2(\R^d)}.
\end{align*}
Ce qui démontre la proposition. \hfill  $ \boxtimes  $
\\On termine la section avec une propriété de calcul fonctionnel qui explique que certains opérateurs peuvent être approximés par des opérateurs pseudo-différentiels.
\begin{prop}  \label{approx} Soient $ \Phi \in C^\infty_0 ( \R ) $ et $ \chi_2 \in C^\infty_0 ( \R^d ) $ avec $ \chi_2 (x) = 1 $ pour $ x \in B(0,1^+ )$. Alors pour tout $ N \in \N^* $ et $ s \geq 0 $, il existe une constante $ C_{N,s } > 0  $ telle que pour tout $ h \in ]0,1 ] $ et $ u \in L^2(\R^d) $,
\begin{equation*}
 || \Phi ( x^2 + ( h D ) ^2 )u - \sum_{j=0}^{N-1} h^j  Op_h(  \Psi_j(x , \xi ) ) \chi_2  u ||_{H^s(\R^d)} \leq C_{N,s} h^{N-s} ||u||_{L ^2(\R^d)},  
 \end{equation*} 
où $ \Psi_0(x, \xi ) = \Phi ( x^2 + \xi ^2  ) $, $ Supp ( \Psi_j ) \subset ( (x, \xi ) / x^2 + \xi^2 \in Supp ( \Phi ) ) $ et $ \Psi_j \in T^{-j} \subset S^0 $.
\end{prop}
\textit{Preuve.} On se sert de la proposition 2.1 de \cite{burq1}. 
\\Si $ \chi _1 \chi_2 =  \chi_1 $  alors
\begin{equation*}
|| \Phi ( x^2 + ( h D ) ^2 ) \chi_1 u - \sum_{j=0}^{N-1} h^j Op_h( \Psi_j(x , \xi ) ) \chi_2 u ||_{H^s(\R^d)} \leq C_{N,s} h^{N-s} ||u||_{L ^2(\R^d)}, 
\end{equation*} 
avec $ \Psi_0(x, \xi ) = \Phi ( x^2 + \xi^2  ) $, $ Supp ( \Psi_j ) \subset ( (x, \xi ) / x^2 + \xi^2 \in Supp ( \Phi ) ) $ et $ \Psi_j \in T^{-j} $.
\\\\Puis, il suffit de choisir correctement $ \chi_1 $ pour avoir
\begin{equation*}
 || \Phi ( x^2 + ( h D ) ^2 )(1- \chi_1 ) u ||_{H^s(\R^d)} \leq h^\infty ||u||_{L ^2(\R^d)},
 \end{equation*}
c'est à dire
\begin{equation*}
 || \Phi \left( \frac{H}{N^2}  \right)(1- \chi_1 ) \left( \frac{x}{N}  \right)  u ||_{H^s(\R^d)} \leq N^{-\infty} ||u||_{L ^2(\R^d)},
 \end{equation*}
où $ N^2= \frac{1}{h}$. Comme $ \Phi $ est bornée et à support compact, on trouve 
\begin{equation*}
|| \Phi \left( \frac{H}{N^2}  \right) v ||_{H^s(\R^d)} \leq N^s ||v||_{L ^2(\R^d)}, \ \forall v \in L ^2(\R^d)
\end{equation*}
Il suffit de vérifier 
\begin{equation*}
|| (1- \chi_1 ) \left( \frac{x}{N}  \right)  u ||_{L^2(\R^d)} \leq N^{-\infty} ||u||_{L ^2(\R^d)},
\end{equation*}
pour $ u=h_n$, pour tout $ n \in \N $. 
\\Il suffit alors d'utiliser la proposition \ref{propre} avec $ \chi _1 =1 $ sur $ B(0,1^+)$. \hfill $ \boxtimes $
\\\\Revenons à la preuve de l'inégalité (\ref{biblitris}). On peut écrire 
\begin{equation*}
u_M = \chi \left( \frac{4x^2}{M^2} \right) u_M  + (1-\chi ) \left( \frac{4x^2}{M^2} \right) u_M
\end{equation*}
avec $ \chi \in C^\infty_0 (  \R^d ) $  vérifiant 
\begin{equation*}
\chi(x) = \left\{
    \begin{array}{lll}
        1 \ & \mbox{si } \ &|x| \leq \frac{15}{16},  \\
        0 \ & \mbox{si} \ &|x| \geq 1.
    \end{array}
\right.
\end{equation*}
Et, par l'inégalité triangulaire, nous devons estimer les deux termes suivants: 
\begin{equation} \label{term1}
 || e^{it H} \chi \left( \frac{4x^2}{M^2}  \right) u_M \ e^{it  H} v_N ||_{L^2( [-\epsilon, \epsilon] ,   L^2(\R^d))  } 
\end{equation}
et
\begin{equation} \label{term2}
|| e^{it H} (1-\chi) \left( \frac{4x^2}{M^2} \right) u_M \  e^{it H } v_N ||_{ L^2( [-\epsilon, \epsilon] ,   L^2(\R^d))  }. 
\end{equation}
\subsection{Estimation du premier terme: (\ref{term1})}
Le théorème 2.4 de \cite{staffilani} nous donne le résultat suivant:
\begin{prop} \label{steffi} Pour tout $ \delta > 0 $, il existe une constante $ C_\delta $ telle que pour tous $ u \in H^{ -1/2 +  \delta }(\R^d) $ et $ v \in H^{ \frac{d-1}{2} - \delta }(\R^d)$, 
\begin{equation*}
 || e^{it \Delta} u  \ e^{it  \Delta} v ||_{L^2 ( \R , L^2 ( \R^d ))  } \leq C_\delta || u ||_{H^{ -1/2 +  \delta }(\R^d)} ||v ||_{H^{ \frac{d-1}{2} - \delta }(\R^d)} .
\end{equation*}
\end{prop}
Ainsi, à l'aide de la transformation de lentille, on en déduit la proposition suivante:
\begin{prop} \label{steffibis}
Pour tout $ \delta > 0 $, il existe une constante $ C_\delta $ telle que pour tous $ u \in H^{ -1/2 +  \delta }(\R^d) $ et $ v \in \overline{H}^{ \frac{d-1}{2} - \delta }(\R^d)$, 
\begin{equation*}
 || e^{it H } u \ e^{it  H } v ||_{L^2 ( [ - \frac{\pi}{8}, \frac{\pi}{8} ] , L^2 ( \R^d ))  } \leq C_{\delta} || u ||_{H^{ -1/2 +  \delta }(\R^d)} ||v ||_{\overline{H}^{ \frac{d-1}{2} - \delta }(\R^d)} .
 \end{equation*}
\end{prop}
\textit{Preuve.} En utilisant les propositions \ref{lechangementdevariable} et \ref{steffi}, on obtient que
\begin{align*}
& || e^{it H} u \ e^{it  H} v ||^2_{L^2 ( [-\frac{\pi}{8} ; \frac{\pi}{8} ] , L^2 ( \R^d ))  } 
\\  = \ & || e^{-it H} u \ e^{-it  H} v ||^2_{L^2 ( [-\frac{\pi}{8} ; \frac{\pi}{8} ] , L^2 ( \R^d ))  } 
 \\ = \ & \int_ { ]-\frac{\pi}{8} ; \frac{\pi}{8} [ }\int _{ \R ^d } \frac{1}{ | \cos(2t) |^{2d}  } \times | e^{it\Delta} u \ e^{it\Delta}v |^2 \left( \frac{ \tan(2t)}{2} , \frac{x}{ \cos(2t) } \right) \ dt dx  
 \\  = \ & \int_ { ]-\frac{\pi}{8} ; \frac{\pi}{8} [ }\int _{ \R ^d } \frac{1}{ | \cos(2t) |^d  } \times  | e^{it\Delta} u \ e^{it\Delta}v |^2 \left( \frac{ \tan(2t)}{2} , x \right) \ dt dx 
 \\  = \ & \int_ {-1/2}^{1/2}  \int _{ \R ^d } (1+(2t)^2)^{d/2-1} \times | e^{it\Delta} u e^{it\Delta}v |^2 ( t , x ) \ dt dx  
\\  \leq \ &  C  \times \int_\R  \int _{ \R ^d }  | e^{it\Delta} u \ e^{it\Delta}v |^2 ( t , x ) \ dt dx  
 \\ \leq \ &  C \times  || u ||_{H^{-1/2 +  \delta }(\R^d)}^2 \times ||v || ^2 _{H^{ \frac{d-1}{2} - \delta }(\R^d)} 
 \\ \leq \ & C \times  || u ||_{H^{-1/2 +  \delta }(\R^d)}^2 \times ||v || ^2 _{\overline{H}^{ \frac{d-1}{2}  - \delta }(\R^d)},
 \end{align*} 
 où dans la dernière inégalité on utilise la proposition \ref{comparaison}. \hfill $ \boxtimes $
\\\\Nous allons montrer que pour tout $ \delta \in ]0 , \frac{1}{2} ] $, il existe une constante $ C_\delta > 0 $ telle que pour tous $ N,M, u $ et $ v $, si $ M \geq N $ alors
\begin{equation*}
|| e^{itH } v_N e^{itH } \chi \left( \frac{4x^2}{M^2} \right) u_M  ||_{L^2(  [  -\frac{\pi}{8} , \frac{\pi}{8} ], L^2(\R^d)) }  \leq C_\delta N^{  \frac{d-2}{2} } \left( \frac{N}{M} \right) ^{1/2-\delta} ||u||_{L^2(\R^d)} ||v||_{L^2(\R^d)}.
\end{equation*}  
Pour cela, en utilisant la proposition \ref{steffibis}, il suffit de prouver que pour tout $  \delta \in ]0 , \frac{1}{2} ] $, il existe une constante $ C_\delta > 0 $ telle que pour tous $ M$ et $ u $
\begin{equation*} 
|| \chi  \left( \frac{4x^2}{M^2} \right) u_M ||_{H^{-1/2+\delta}(\R^d) } \leq C  M^{-1/2+\delta} || u  ||_{L^2(\R^d)}.
\end{equation*}
Trivialement, nous avons que
\begin{equation*}
|| \chi \left( \frac{4x^2}{M^2} \right) u_M ||_{L^2(\R^d)} \leq || u_M ||_{L^2(\R^d)} \leq C ||u||_{L^2(\R^d)}.
\end{equation*}
Ainsi, par interpolation, il suffit de démontrer qu'il existe une constante $ C > 0 $ telle que pour tout $ M$ et $ u $,
\begin{equation} \label{class}
|| \chi \left( \frac{4x^2}{M^2} \right) u_M ||_{H^{-1}(\R^d)} \leq C M^{-1} ||u||_{L^2(\R^d)}.
\end{equation}
On utilise alors le calcul semi classique. Pour une fonction $ u $, on définit $  \mathfrak{u} : x \longmapsto u( \frac{x}{     \sqrt{ h }  } ) $ où $ h = \frac{1}{M^2} $.
\\\\Remarquons que
\begin{equation} \label{semi-classique1}
\chi \left( \frac{4x^2}{M^2}  \right)  \left[ \frac{H}{M^2} \right] ( u ) (x) =  [   \chi( 4x^2 )  (  -h^2 \Delta + |x|^2  )  ] ( \mathfrak{u} )  (  \sqrt{ h } x  )
\end{equation}
et que
\begin{equation} \label{semi-classique2}
\chi \left( \frac{4x^2}{M^2} \right)  \left[   \phi \left(   \frac{H}{M^2} \right) \right] ( u ) (x) =  [ \chi( 4x^2 ) \phi   (  -h^2 \Delta + |x|^2  )   ] (  \mathfrak{u}  )( \sqrt{  h } x  ).
\end{equation}
Ainsi, pour prouver (\ref{class}), il suffit d'établir l'existence d'une constante $ C > 0 $ telle que pour tout $ h \in ]0,1] $ et $ u \in L^2(\R^d)$,
\begin{equation} \label{classb}
 || \chi \left( 4 x ^2 \right) \phi ( x^2 + ( h D ) ^2 ) u ||_{H^{-1}(\R^d) } \leq C h ||u||_{L^2(\R^d)}.
\end{equation}
En effet,
\begin{align*}
|| \chi \left( \frac{4x^2}{M^2} \right) u_M ||_{H^{-1}(\R^d)} & \leq || \ [ \chi( 4x^2 ) \phi   (  -h^2 \Delta + |x|^2  )   ] (\mathfrak{u}  )( \sqrt{  h } .  )  ||_{H^{-1}(\R^d) }
\\ & \leq h^{- \frac{1}{2} - \frac{d}{4} }  || [ \chi( 4x^2 ) \phi   (  -h^2 \Delta + x^2  )   ] (\mathfrak{u}  )  ||_{ H^{-1}(\R^d) }
\\ & \leq C  h^{\frac{1}{2} - \frac{d}{4} } ||  \mathfrak{u} || _{L^2(\R^d)}
\\ & \leq C h^{1/2}  ||u||_{L^2(\R^d)}
\\ & \leq C M^{-1} ||u||_{L^2(\R^d)}.
\end{align*}
Démontrons ensuite (\ref{classb}). Grâce à la proposition \ref{approx}, on a 
\begin{align*}
 & || \chi( 4 x^2  ) \phi ( x^2 + ( h D ) ^2 ) u ||_{H^{-1}(\R^d)}  
 \\  = \ & ||   \chi( 4x^2 ) [ \phi ( x^2 + ( h D ) ^2 ) u   - Op_h ( \phi ( x^2 +  \xi  ^2 )  ) \chi_2 u + Op_h ( \phi ( x^2 +   \xi  ^2 )  ) \chi_2 u ] ||_{H^{-1}(\R^d)} 
 \\   \leq \ & h || u ||_{L^2(\R^d) }   + || \chi ( 4x^2 ) Op_h ( \phi ( x^ 2 + \xi ^2 ) ) \chi_2 u ||_{H^{-1}(\R^d)}.
 \end{align*}
Ainsi, il suffit d'évaluer $ || \chi ( 4x^2 ) Op_h ( \phi ( x^ 2 + \xi ^2 ) ) \chi_2 u ||_{H^{-1}(\R^d)} $. On a 
\begin{align*}
  & \chi ( 4x^2 ) Op_h ( \phi ( x^ 2 +  \xi  ^2 ) ) \chi_2 u 
  \\  =  & \frac{\chi ( 4x^2 )}{(2 \pi)^d} \times \int_{\R^d }  e^{ix \xi }  \phi ( x^ 2 + ( h \xi ) ^2 )  \mathcal{F} ( \chi_2 u ) ( \xi ) \ d\xi  
  \\  = & \frac{1}{(2 \pi)^d} \times \int_{\R^d } \frac{ih \xi}{ih\xi } e^{ix \xi }   \chi ( 4 x^2 ) \phi ( x^ 2 + ( h \xi ) ^2 )  \mathcal{F} (\chi_2 u) ( \xi ) \ d \xi 
  \\  = & \frac{h}{(2 \pi)^d} \times \int_{\R^d } \nabla_x ( e^{ix \xi } )  \frac{ \chi ( 4 x^2 )  }{ih \xi}  \phi ( x^2 + ( h  \xi )^2 )  \mathcal{F} (\chi_2 u)( \xi ) \ d\xi     
\\   = & \frac{h}{(2 \pi)^d} \times \left( \nabla_x \left( \int_{\R^d }  e^{ix \xi }   \frac{ \chi ( 4 x^2 )  }{ih \xi}  \phi ( x^2 + ( h  \xi )^2 ) \mathcal{F}(\chi_2 u)( \xi ) \ d\xi \right)  \right.
\\ & \hspace*{2cm} -  \left. \int_{\R^d }  e^{ix \xi } \nabla_x \left( \frac{ \chi ( 4 x^2 )  }{ih \xi}  \phi ( x^2 + ( h  \xi )^2 ) \right) \mathcal{F} (\chi_2 u)( \xi ) \ d\xi  \right).  
\end{align*}
Puis, comme $ 4 x^2  \leq 1 $ et $ \frac{1}{2}^- \leq x^2 +  \xi ^2 \leq 2^+   $  implique $ \frac{1}{4}^- \leq \xi^2 $, on en déduit que
\begin{align*}
( x , \xi ) \longrightarrow \frac{ \chi ( 4 x^2 )  }{i\xi}  \phi ( x^2 +   \xi ^2 ) \ \in S^0 \ \mbox{et} \ ( x , \xi ) \longrightarrow \nabla_x \left( \frac{ \chi ( 4 x^2 )  }{i \xi}  \phi ( x^2 +  \xi ^2 ) \right) \ \in S^0.
\end{align*}
Ainsi, par le théorème \ref{operateur}, on a alors
\begin{align*}
& \bigg| \bigg| \frac{h}{(2 \pi)^d} \times \left( \nabla_x \left( \int_{\R^d }  e^{ix \xi }   \frac{ \chi ( 4 x^2 )  }{ih \xi}  \phi ( x^2 + ( h  \xi )^2 ) \mathcal{F}(\chi_2 u)( \xi ) \ d\xi \right)  \right) \bigg| \bigg|_{H^{-1}(\R^d)}
\\  \leq \ & \frac{h}{(2 \pi)^d} \times  \bigg| \bigg|  \int_{\R^d }  e^{ix \xi }   \frac{ \chi ( 4 x^2 )  }{ih \xi}  \phi ( x^2 + ( h  \xi )^2 ) \mathcal{F}(\chi_2 u)( \xi ) \ d\xi \bigg| \bigg|_{L^2(\R^d)} 
\\  \leq \ & C h  || \chi_2  u  ||_{L^2(\R^d)} 
\\  \leq  \ & C h  ||  u  ||_{L^2(\R^d)}
\end{align*}
et
\begin{align*}
& \bigg| \bigg| \frac{h}{(2 \pi)^d} \times \int_{\R^d }  e^{ix \xi } \nabla_x \left( \frac{ \chi ( 4 x^2 )  }{ih \xi}  \phi ( x^2 + ( h  \xi )^2 ) \right) \mathcal{F} (\chi_2 u)( \xi ) \ d\xi \bigg| \bigg|_{  H^{-1}(\R^d)  }
\\  \leq \ & \bigg| \bigg| \frac{h}{(2 \pi)^d} \times \int_{\R^d }  e^{ix \xi } \nabla_x \left( \frac{ \chi ( 4 x^2 )  }{ih \xi}  \phi ( x^2 + ( h  \xi )^2 ) \right) \mathcal{F} (\chi_2 u)( \xi ) \ d\xi \bigg| \bigg|_{  L^2 (\R^d)  }
\\  \leq \ & C h || \chi_2  u  ||_{L^2(\R^d)} 
\\ \leq \ & C h ||  u  ||_{L^2(\R^d)}.
\end{align*}
Ce qui démontre (\ref{classb}).
\subsection{Estimation du second terme: (\ref{term2})}
\begin{prop} \label{clef} Il existe un temps $ T \in ] 0 , \frac{\pi}{4} [$ tel que pour tout $ N \geq 0 $, il existe une constante $ C_N > 0 $ telle que pour tout $ M \geq 1 $ et $u \in L^2(\R^d) $,
\begin{equation*}
 || \chi \left( \frac{8 x^2}{M^2} \right) e^{itH} (1-\chi) \left( \frac{4 x^2}{M^2} \right) u_M   ||_{L^2 ( [-T , T ] , L^2 ( \R^d ) )}  \leq C_N M^{-N} ||u||_{L^2(\R^d)}.  
 \end{equation*}
\end{prop}
Montrons que la proposition \ref{clef} implique l'estimation du terme (\ref{term2}).
\\\\De
\begin{equation*}
||u||_{L^\infty( [a,b] ,  L^\infty (\R^d))} \leq C ||u||_{L^\infty ( [a,b] , \overline{H}^{ \frac{d}{2}+1 }(\R^d)) } 
\end{equation*}
et la proposition \ref{clef}, on déduit pour $ M \geq 10 N $ que
\begin{align*}
 & || \chi \left( \frac{8 x^2}{M^2} \right) e^{itH} (1-\chi) \left( \frac{4 x^2}{M^2} \right) u_M \ e^{itH } v_N  ||_{L^2([-T, T] , L^2(\R^d))} 
 \\ \leq \ & || \chi \left( \frac{8 x^2}{M^2} \right) e^{itH} (1-\chi) \left( \frac{4 x^2}{M^2} \right) u_M||_{L^2([-T, T] , L^2(\R^d)) }   \times || e^{itH} v_N||_{L^\infty([-T, T] , L^\infty(\R^d))}   
 \\ \leq \ &  C \times || \chi \left( \frac{8 x^2}{M^2} \right)  e^{itH} (1-\chi) \left( \frac{4 x^2}{M^2} \right) u_M||_{L^2([-T, T] , L^2(\R^d))  }  \times  ||  v_N||_{ \overline{H} ^{\frac{d}{2} +1 } ( \R^d )  }  
 \\ \leq \ & C_K M^{-K} N^{ \frac{d}{2} +1 } ||u||_{L^2(\R^d)}   ||v||_{L^2(\R^d)}  
 \\ \leq \ & C_K M^{-K + \frac{d}{2} +1 }   ||u||_{L^2(\R^d)} ||v||_{L^2(\R^d)}. 
\end{align*}
Et donc, il suffit d'estimer convenablement pour $ M \geq 10 N$, le terme suivant:
\begin{equation}
|| e^{itH} v_N    \  (1-\chi ) \left( \frac{8 x^2}{M^2} \right) e^{itH} (1-\chi) \left( \frac{4 x^2}{M^2} \right) u_M ||_{L^2([-T, T] , L^2(\R^d))  }. 
\end{equation}  
Grâce au théorème \ref{Stricharz}, on trouve pour tout $ R \geq 1$ ,
\begin{align*}
& || e^{itH} \ v_N   (1-\chi ) \left( \frac{8 x^2}{M^2} \right) e^{itH} (1-\chi) \left( \frac{4 x^2}{M^2} \right) u_M ||_{L^2([-T, T] , L^2(\R^d))  } 
\\ \leq \ &   || <x> ^ R e^{itH} v_N  \  <x>^{-R}(1-\chi ) \left( \frac{8 x^2}{M^2} \right) e^{itH} (1-\chi) \left( \frac{4 x^2}{M^2} \right) u_M ||_{ L^2([-T, T] , L^2(\R^d)) } 
\\  \leq \ & M^{-R} \times  || <x> ^ R e^{itH} v_N  (1-\chi ) \left( \frac{8 x^2}{M^2} \right) \ e^{itH} (1-\chi) \left( \frac{4 x^2}{M^2} \right) u_M ||_{ L^2([-T, T] , L^2(\R^d)) } 
\\ \leq \ & M^{-R} \times  || (1-\chi ) \left( \frac{8 x^2}{M^2} \right)  \ <x> ^ R e^{itH} v_N ||_{ L^4([-T, T] , L^4(\R^d)) }
\\ &  \hspace*{5cm}  \times ||  e^{itH} (1-\chi) \left( \frac{4 x^2}{M^2} \right) u_M ||_{ L^4([-T, T] , L^4(\R^d)) } 
\\   \leq \ & M^{-R} \times || (1-\chi )\left( \frac{8 x^2}{M^2} \right)  <x> ^ R e^{itH} v_N ||_{ L^4([-T, T] , L^4(\R^d)) } 
\\ & \hspace*{5cm} \times ||  e^{itH} (1-\chi) \left( \frac{4 x^2}{M^2} \right) u_M ||_{ L^4([-T, T] , \overline{W} ^{  \frac{d-2}{4} , \frac{2d}{d-1} }(\R^d)) } 
 \\  \leq \ & M^{-R} \times ||  <x> ^ R (1-\chi )\left( \frac{8 x^2}{M^2} \right) e^{itH} v_N   ||_{L^\infty([-T, T] , L^4(\R^d))} \times M^{ \frac{d-2}{4}} ||  u_M ||_{L^2(\R^d) }.
  \end{align*} 
Puis, comme $ Supp \left\{ (1-\chi )\left( \frac{8 x^2}{M^2} \right) \right\} \subset \left\{ |x|^2 \geq  \frac{M^2}{8} \times  \frac{15}{16} \right\} \subset \left\{ |x|^2  \geq 5 N^2 \right\} $ alors par la proposition \ref{propre}, on déduit
\begin{equation*}
|| (1-\chi ) \left( \frac{8 x^2}{M^2} \right) <x> ^ R e^{itH} v_N ||_{L^4([-T, T] , L^4(\R^d))} \leq C  ||v||_{L^2(\R^d)}.
\end{equation*}
\textit{Preuve de la proposition \ref{clef}.} En utilisant l'analyse semi classique comme en (\ref{semi-classique2}), il suffit de prouver l'existence d'un temps $ T \in ] 0 , \frac{\pi}{4}[ $ tel que pour tout $ N \geq 1 $, il existe une constante $ C_N > 0 $ telle que pour tout $ u \in L^2(\R^d) $ et $ h \in ]0,1] $,
\begin{equation}
 || \chi( 8 x^2 ) e^{itH_h/h} (1-\chi) ( 4 x^2 ) \phi ( x^2 + (h D )^2  ) u  ||_{L^2( [-T , T ] , L^2 ( \R^d ) )}  \leq C_N h^{N} ||u||_{L^2(\R^d)}.
 \end{equation}
En utilisant la proposition \ref{approx}, il suffit d'établir le résultat suivant:
\\Pour toute fonction $ g(x, \xi )  \in  C_0^\infty( \R^d * \R^d ) $ vérifiant $ Supp ( g )  \subset  \lbrace \frac{1}{4} \leq \xi ^2 + x^2 \leq 4 \rbrace $ et tout entier $ N \geq 1 $, il existe une constante $ C_N > 0 $ telle que pour tout $ u \in L^2(\R^d) $ et $ h \in ]0,1] $,
\begin{equation}
|| \chi( 8 x^2 ) e^{ihH_h/h} (1-\chi) ( 4 x^2 ) Op_h( g ( x , \xi  ) ) u  ||_{L^2 ( [-T , T] ,L^2 ( \R^d ) )}  \leq C_N h^{N}  ||u||_{L^2(\R^d)}.
\end{equation}
\\En effet, par la proposition \ref{approx},
\begin{align*}
|| \chi( 8 x^2 ) e^{itH_h/h} (1-\chi) ( 4 x^2 ) \phi ( x^2 + (h D )^2  ) u  ||_{L^2( [-T , T ] , L^2 ( \R^d ) )}
\\\leq  || \chi( 8 x^2 ) e^{itH_h/h} (1-\chi) ( 4 x^2 ) [ \phi ( x^2 + (h D )^2  ) u  -  \sum_{j=0}^{N-1} h^j Op_h(\Psi_j(x , \xi ))  u ]  ||_{L^2( [-T , T ] , L^2 ( \R^d ) )}  
\\ + \sum_{j=0}^{N-1} h^j || \chi( 8 x^2 ) e^{itH_h/h} (1-\chi) ( 4 x^2 ) Op_h ( \Psi_j ( x , \xi) ) u ||_{L^2([-T,T], L^2( \R ^ d)) } 
\\  \leq  C_N h^{N} ||u||_{L^2 ( \R^d ) }  + \sum_{j=0}^{N-1} h^j || \chi( 8 x^2 ) e^{itH_h/h} (1-\chi) ( 4 x^2 ) Op_h( \Psi_j ( x, \xi  ) ) u ||_{L^2([-T,T], L^2 ( \R ^ d)) } 
\\ \leq C_N h^{N} ||u||_{L^2(\R^d)}.
\end{align*}
\begin{lem} Il existe un temps $ T \in ] 0 , \frac{\pi}{4}[ $ tel que si $ g $ est une fonction $ C_0^\infty( \R^d * \R^d )$ vérifiant $ Supp ( g )  \subset  \lbrace \frac{1}{4} \leq \xi ^2 + x^2 \leq 4 \rbrace $ alors pour tout entier $ N \geq 1 $, il existe une constante $ C_N > 0 $ telle que pour tout $ u \in L^2(\R^d) $ et $ h \in ]0,1] $,
\begin{equation*}
|| \chi( 8 x^2 ) e^{it H_h/h } (1-\chi) ( 4 x^2 ) Op_h( g ( x, \xi ) ) u  ||_{L^2 ( [-T , T] , L^2 ( \R^d ) )}  \leq C_N h^{N} ||u||_{L^2(\R^d)}.
\end{equation*}
\end{lem}
\textit{Preuve.} On définit
\begin{equation*}
w(s,x) = \int_{\R^d} e^{ \frac{i}{h} \Phi(s,x,\xi )  }  a ( s,x, \xi , h ) \hat{u} ( \frac{\xi}{h} ) \frac{d \xi }{ ( 2 \pi h )^d }
\end{equation*}  où
\begin{equation*}
a(s,x, \xi,h) = \sum_{j=0} ^N h^j a_j ( s,x,\xi).
\end{equation*}
Supposons que
\begin{equation*}
\left\{
    \begin{array}{ll} \label{hamilton}
        \Phi ( 0, x, \xi )= x . \xi, \\
        \partial_s  \Phi - | \nabla \Phi |^2 - x^2 = 0, 
    \end{array}
\right. 
\end{equation*}
\begin{equation*}
\left\{
    \begin{array}{ll}
        \ a_0 ( 0, x, \xi )= (1-\chi)( 4 x^2 ) g ( x, \xi  ), \\
        \partial_s a_0 - 2 \nabla \Phi . \nabla a_0 - \Delta ( \Phi ) a_0 = 0,
    \end{array}
\right. 
\end{equation*}
et
\begin{equation*}
\left\{
    \begin{array}{ll}
        \ a_j ( 0, x, \xi )= 0,  \\
        \partial_s a_j - 2 \nabla \Phi . \nabla a_j - \Delta ( \Phi ) a_j =  -i\Delta ( a_{j-1} ),  
    \end{array}
\right. 
\end{equation*}
\hspace*{7cm} pour $ 1 \leq j \leq N. $ 
\\\\Alors
\begin{align*}
 i h  \partial _s w -h^2 \Delta w + x^2 w & = - h^{N+2} \int_{\R^d} e^{ \frac{i}{h} \Phi(s,x,\xi )  } \Delta ( a_N ( s,x, \xi) ) \hat{u} ( \frac{\xi}{h} ) \frac{d \xi }{ ( 2 \pi h )^d }  
 \\  :&= h^{N+2 } f.
 \end{align*}
Par conséquent 
\begin{align*} 
& \chi( 8 x^2 ) e^{it H_h/h } (1-\chi) ( 4 x^2 ) Op( g ( x,h \xi  ) ) u 
\\   = \ & \chi(8x^2) w ( t, x)  - i h^{N+2} \chi ( 8x^2 ) \int_0^t e^{i(t -s ) H_h/h } f(s) ds.  
\end{align*}
Remarquons que si $ \Psi $ est solution de l'équation $  \partial  _t \Psi + | \nabla \Psi  | ^2 = 0  $ avec donnée initiale $ \Psi ( 0,x,\xi) = x. \xi $ alors $ \Phi(t,x,\xi ) = \Psi (  \frac{ -\tan 2 t }{2 }  , \frac{x}{ \cos 2t } , \xi ) +  \frac{x^2  \tan 2 t }{ 2 } $ est solution de l'équation $ \partial_s  \Phi - | \nabla \Phi |^2 - x^2 = 0 $ avec même donnée initiale.
\\\\Par la méthode des caractéristiques, on trouve
\begin{equation*}
\Psi ( t ,x, \xi ) = - t | \xi |^2 +x.\xi,
\end{equation*}
puis on déduit que
\begin{equation*}
 \Phi ( t, x , \xi ) = \frac { \tan(2t) }{2} ( \xi^2  + x^2 ) + \frac{x.\xi}{\cos( 2t) }.
\end{equation*}
Ainsi
\begin{equation*}
 \nabla \Phi = \frac{ \xi}{ \cos ( 2t) } + x \tan ( 2 t )  \ \mbox{et}  \  \Delta \Phi =  d \tan( 2t ). 
\end{equation*}
En utilisant la méthode des caractéristiques, on trouve 
\begin{equation} \label{a0}
 a_0 \left( t , x - 2 \int_0^t \nabla \Phi , \xi \right) = \frac{ a_0 ( 0 ,x,\xi ) }{ | \cos 2 t |^{  \frac{d}{2}  } },   
 \end{equation}
et
\begin{equation} \label{aj}
a_j \left( t, x - 2 \int_0^t \nabla \Phi , \xi \right) = -i \int_0^t  \bigg| \frac{ \cos ( 2 s ) }{ \cos ( 2 t )  } \bigg|^{  \frac{d}{2} } \  \Delta a_{j-1} \left( s,x - 2 \int_0^s \nabla \Phi , \xi \right) ds.
\end{equation}
Or, pour tout $ \xi \in \R^d $ et $ | x | \leq \frac{1}{2} $,
\begin{equation*}
a_0 ( 0 ,x,\xi ) =0
\end{equation*}
Par conséquent, pour tout $ \xi \in \R^d $, $ t \in ] - \frac{\pi}{4} , \frac{\pi}{4}  [ $, $ | x | \leq \frac{  \sqrt{15} }{8} $ et $ j \in \N$
\begin{equation*}
a_j \left( t, x - 2 \int_0^t \nabla \Phi , \xi \right) = 0 .
\end{equation*}
Or $ \int _0^t \nabla \Phi =  \xi F ( t) - \frac{x \log \cos (2 t ) }{2} $ avec F une fonction continue vérifiant $ F ( 0) = 0 $. Ainsi, pour tout $ \epsilon > 0 $, il existe un temps $ T \in ]0, \frac{\pi}{4}  ] $ tel que si $ |t| \leq T $ alors $ |F(t)| \leq \epsilon $ et $ |  \log \cos 2 t  | \leq \epsilon $.
\\\\Remarquons que
\begin{align*}
 & y   = x - 2 \int_0^t \nabla \Phi = x ( 1+ \log \cos 2t ) - 2 \xi F(t)
\\ \Leftrightarrow \hspace*{1cm} & x = \frac{y+2 \xi F(t)}{1+\log \cos 2t }
\end{align*}
et que si $ | y | \leq \frac{1}{ \sqrt{8} } $, $ \xi^2 \leq  4 $ et $ |t| \leq T $ alors $ |x| \leq \frac{1/ \sqrt{8} + 4 \epsilon}{1-\epsilon}  \leq \frac{  \sqrt{15} }{8} $, si $ \epsilon $ choisi correctement.
\\\\Cela implique que pour tout $ j \in \N$ et $ |t| \leq T $
\begin{align*}
 Supp ( a_j( t) )  \subset  B_x \left( 0, \frac{1}{  \sqrt{8} }  \right) ^c \times B_\xi \left( 0, 2  \right) .
\end{align*}
Par conséquent, si $ |t| \leq T $, comme $ Supp ( \chi(8 x^2) ) \subset B_x( 0, \frac{1}{\sqrt{8} }) $, on déduit que
\begin{equation*}
\chi( 8 x^2 ) e^{it H_h/h} (1-\chi) ( 4 x^2 ) Op_h( g(x,\xi) )  u = - ih^{N+2} \chi ( 8x^2 ) \int_0^t e^{i(t-s) H_h/h} f(s)  \ ds.
\end{equation*}
Puis, par le théorème \ref{Stricharz2}, on trouve
\begin{align*}
& || \chi( 8 x^2 ) e^{itH_h/h } (1-\chi) ( 4 x^2 ) Op_h( g ( x,\xi ) ) u ||_{L^2( [-T,T] ,  L^2(\R^d))}  
\\  \leq \ & h^{N+2} \bigg| \bigg|     \chi ( 8x^2 ) \int_0^t e^{i(t-s)  H_h/h } f(s) \ ds   \bigg| \bigg| _{L^2( [-T,T] ,  L^2(\R^d)) } 
\\ \leq \ & h^{N+2}  ||f||_{L^1( [-T,T  ] ,  L^2(\R^d)) }  
\\ \leq \ & h^{N+2} || \Delta a_N ||_{L^1_t( [-T,T ]  , L^2_x(\R^d , L ^2_\xi(\R^d))) } \times  ||u||_{L^2(\R^d)}.
\end{align*}
Le lemme est donc démontré si $ \Delta a_N \in L^1_t( [-T,T] , L^2_x(\R^d  , L ^2_\xi(\R^d))) $.
\\\\On démontre par récurrence sur $ N \in \N $, que pour tout $ \alpha \in \N^d, \\ \hspace*{2cm}  \partial^\alpha_x  a_N \in L^1_t( [-T,T ] , L^2_x(\R^d  , L ^2_\xi(\R^d))) $.
\\\\Pour $ N=0 $, à l'aide de (\ref{a0}), on voit par changement de variables que 
\\ \hspace*{2cm} $ \partial^\alpha_x  a_0\in L^1_t( [-T,T ] , L^2_x(\R^d  , L ^2_\xi(\R^d))) $ si $ \partial^\alpha_x a_0(0) \in L^2_x(\R^d, L ^2_\xi(\R^d)) $. 
\\Or $ a_0 ( 0 ) \in C ^\infty (\R^d \times \R^d )$ avec $ Supp ( a_0 ( 0 ) ) \subset \lbrace (x ,\xi ) / x^2 \leq 1 , \xi^2 \leq  4  \rbrace $ et le cas $ N= 0 $ est évident. 
\\\\Supposons le résultat établi au rang $ N -1 $ et montrons le au rang $ N $. \`A l'aide de (\ref{aj}), on note que $ \partial^\alpha_x  a_N \in L^1_t( [-T,T ] , L^2_x(\R^d  , L ^2_\xi(\R^d))) $ si  $ \partial^{\alpha+2}_x  a_{N-1} \in L^1_t( [-T,T ] , L^2_x(\R^d  , L ^2_\xi(\R^d))) $. Cette dernière affirmation étant claire par hypothèse de récurrence. \hfill  $  \boxtimes  $
\subsection{Estimées bilinéaires et espaces de Bourgain}
L'objectif de cette section consiste à écrire l'estimée bilinéaire du théorème \ref{bilis} dans les espaces de Bourgain. Plus précisément, on cherche à établir les deux théorèmes suivants:
\begin{thm} \label{bilibourgain1}
Il existe $ \delta_0 \in ]0 , \frac{1}{2} ]$ tel que pour tout $ \delta \in ]0, \delta_0 ] $, il existe $ b' < \frac{1}{2} $ et une constante $ C > 0 $ tels que pour tous $ u,v,M,N $,
\begin{align*}
& ||  \Delta_N ( v )   \Delta_M (u) ||_{L^2 ( \R , L^2 ( \R^d ))  } 
\\ \leq \ & C \times \min ( N,M)^{ \frac{d-2}{2} +\delta } \times  \left( \frac{ \min ( N,M)  }{ \max ( N,M )}   \right) ^{1/2-\delta } \times  || \Delta_N(v) ||_{  \overline{X}^{0,b'}} || \Delta_M(u) ||_{  \overline{X}^{0,b'}}.
\end{align*}
\end{thm}
\begin{thm} \label{bilibourgain2}
Soit $ \psi \in C_0^\infty (\R) $ alors il existe $ \delta_0 \in ]0 , \frac{1}{2} ]$ tel que pour tout $ \delta \in ]0, \delta_0 ] $, il existe $ b' < \frac{1}{2} $ et une constante $ C > 0 $ tels que pour tous $ u,u_0,M,N $,
\begin{align*}
& ||  \Delta_N ( \psi(t)  e^{-itH}  u_0 )  \Delta_M (u) ||_{L^2 ( \R , L^2 ( \R^d ))  } 
\\ \leq \ & C \times \min ( N,M)^{ \frac{d-2}{2} +\delta } \times  \left( \frac{ \min ( N,M)  }{ \max ( N,M )} \right) ^{1/2-\delta } \times || \Delta_N(u_0) ||_{ L^2(\R^d)} || \Delta_M(u) ||_{  \overline{X}^{0,b'}}.
\end{align*}
\end{thm}
Pour démontrer ces théorèmes, on adapte la preuve du lemme 4.4 de \cite{burq3}. Commençons par remarquer qu'il suffit de démontrer les deux propositions suivantes:
\begin{prop} \label{bilibourgain3}
Pour tout $ b \in ] \frac{1}{2} , 1 ] $ et $ \delta \in ]0, \frac{1}{2} ] $, il existe une constante $ C > 0 $ telle que pour tous $ u,v,M,N $,
\begin{align*}
& ||  \Delta_N ( v )   \Delta_M (u) ||_{L^2 ( \R , L^2 ( \R^d ))  } 
 \\ \leq \ & C \times  \min ( N,M)^{ \frac{d-2}{2} } \times  \left( \frac{ \min ( N,M)  }{ \max ( N,M )} \right) ^{1/2-\delta } \times || \Delta_N(v) ||_{  \overline{X}^{0,b}} || \Delta_M(u) ||_{  \overline{X}^{0,b}}.
\end{align*}
\end{prop}
\begin{prop} \label{bilibourgain4}
Soit $ \psi \in C_0^\infty (\R) $ alors pour tout $ b \in ] \frac{1}{2} , 1 ] $ et $ \delta \in ]0, \frac{1}{2} ] $, il existe une constante $ C > 0 $ telle que pour tous $ u_0,v,M,N $,
\begin{align*}
& ||  \Delta_N ( \psi(t) e^{-itH} u_0  )   \Delta_M (u) ||_{L^2 ( \R , L^2 ( \R^d ))  }
\\  \leq \ & C \times \min ( N,M)^{ \frac{d-2}{2} } \times  \left( \frac{ \min ( N,M)  }{ \max ( N,M )} \right) ^{1/2-\delta } \times || \Delta_N(u_0) ||_{ L^2 (\R^d)} || \Delta_M(u) ||_{  \overline{X}^{0,b}}.
\end{align*}
\end{prop}
En effet, pour tout $ \epsilon > 0 $, d'après la proposition \ref{bourgain2} (avec $ \theta = \frac{1}{2} $), on obtient
\begin{align*}
& || \Delta_N ( v )   \Delta_M (u)  ||_{L^2(\R,L^2(\R^d))} 
\\  \leq \ & ||\Delta_N ( v )||_{L^4(\R,L^2(\R^d))} \times ||\Delta_M (u)||_{L^4(\R,L^\infty(\R^d))}
\\  \leq \ & C  || \Delta_N ( v )||_{  \overline{X}^{0,1/4+\epsilon} } \times  || \Delta_M (u) ||_{  \overline{X}^{d/2+\epsilon,1/4+\epsilon} },
 \\ \\  & || \Delta_N ( \psi(t) e^{itH} u_0 )   \Delta_M (u)  ||_{L^2(\R,L^2(\R^d))} 
 \\  \leq \ & ||\Delta_N ( \psi(t) e^{itH} u_0 )||_{L^4(\R,L^\infty(\R^d))} \times ||\Delta_M (u)||_{L^4(\R,L^2(\R^d))}
\\  \leq \ & C  || \Delta_N ( u_0 )||_{  \overline{H}^{d/2-1/2+\epsilon}(\R^d) } \times  || \Delta_M (u) ||_{  \overline{X}^{0,1/4+\epsilon}(\R*\R^d) }
\\ \leq \ & C  || \Delta_N ( u_0 )||_{  \overline{H}^{d/2+\epsilon}(\R^d) } \times  || \Delta_M (u) ||_{  \overline{X}^{0,1/4+\epsilon}(\R*\R^d) },
\end{align*}
et
\begin{align*}
& || \Delta_N ( \psi(t) e^{itH}u_0 )   \Delta_M (u)  ||_{L^2(\R,L^2(\R^d))} 
\\ \leq \ & ||\Delta_N ( \psi(t) e^{itH} u_0  )||_{L^4(\R,L^2(\R^d))} \times ||\Delta_M (u)||_{L^4(\R,L^\infty(\R^d))}
\\ \leq \ & C  ||\Delta_N ( \psi(t) e^{itH} u_0  )||_{L^\infty(\R,L^2(\R^d))} \times  || \Delta_M (u) ||_{  \overline{X}^{d/2+\epsilon,1/4+\epsilon} }
\\ \leq \ & C  ||\Delta_N ( u_0  )||_{L^2(\R^d)} \times  || \Delta_M (u) ||_{  \overline{X}^{d/2+\epsilon,1/4+\epsilon} }.
\end{align*}
Par conséquent, par interpolation, pour tout $ \theta \in [0,1] $, on trouve
\begin{align*}
||  \Delta_N ( v )   \Delta_M (u)  ||_{L^2 ( \R , L^2 ( \R^d ))  }  \leq & C \times \min( M,N) ^{ \frac{d-2}{2} +\theta ( 1+ \epsilon ) } \times  \left( \frac{ \min(M,N)  }{ \max(M,N) } \right) ^{(1/2-\delta)(1-\theta) }
\\ & \times  || \Delta_N(v) ||_{  \overline{X}^{0,b(1-\theta)+\theta(1/4+\epsilon)}} || \Delta_M(u) ||_{  \overline{X}^{0,b(1-\theta)+\theta(1/4+\epsilon)}}
\end{align*}
et
\begin{align*}
||  \Delta_N & ( \psi(t)  e^{itH} u_0  )   \Delta_M (u)  ||_{L^2 ( \R , L^2 ( \R^d ))  } \leq   C \times  \min( M,N) ^{ \frac{d-2}{2} +\theta ( 1+ \epsilon ) }
\\ &  \times  \left( \frac{ \min(M,N)  }{ \max(M,N) } \right) ^{(1/2-\delta)(1-\theta) }  \times || \Delta_N(u_0) ||_{  L^2(  \R^d ) } \times || \Delta_M(u) ||_{  \overline{X}^{0,b(1-\theta)+\theta(1/4+\epsilon)}}.
\end{align*}
Choisissons $ \delta = \frac{\epsilon}{2} $ et $ \theta = \frac{\epsilon}{4} $ alors
\begin{align*}
b(1-\theta)+\theta( \frac{1}{4} +\epsilon) & =  b - \frac{b\epsilon}{4} + \frac{\epsilon}{4} (  \frac{1}{4}+ \epsilon )
\\ & \leq b - \frac{\epsilon}{8} + \frac{\epsilon}{16} + \frac{\epsilon^2}{4} \leq b - \frac{\epsilon}{17}.
\end{align*}
Il suffit alors de prendre $ b = \frac{1}{2} + \frac{\epsilon}{34} $ et de poser $ b ' = b(1-\theta)+\theta( \frac{1}{4} +\epsilon) < \frac{1}{2} $ pour obtenir
\begin{align*}
& ||  \Delta_N ( v )   \Delta_M (u)  ||_{L^2 ( \R , L^2 ( \R^d ))  } 
\\ \leq \ &  C  \times  \min(N,M)^{ \frac{d-2}{2} + \epsilon }  \times  \left( \frac{ \min(M,N)  }{ \max(N,M) } \right) ^{(1/2-\delta)(1-\theta) } \times  || \Delta_N(v) ||_{  \overline{X}^{0,b'}} || \Delta_M(u) ||_{  \overline{X}^{0,b'}}
\end{align*} 
et
\begin{align*}
& ||  \Delta_N ( \psi(t) e^{-itH} u_0  )   \Delta_M (u)  ||_{L^2 ( \R , L^2 ( \R^d ))  }
\\  \leq \ &  C  \times \min(N,M)^{ \frac{d-2}{2} + \epsilon } \times  \left( \frac{ \min(M,N)  }{ \max(M,N) } \right) ^{(1/2-\delta)(1-\theta) } \times  || \Delta_N(u_0) ||_{ L^2(\R^d)} || \Delta_M(u) ||_{  \overline{X}^{0,b'}}.
\end{align*}
Pour terminer, il suffit de remarquer que 
\begin{equation*}
\left( \frac{1}{2} -\delta \right) (1-\theta) = \frac{1}{2} -\frac{5\epsilon}{8} + \frac{\epsilon^2}{8} \geq \frac{1}{2} - \epsilon  
\end{equation*}
et les théorèmes \ref{bilibourgain1} et \ref{bilibourgain2} suivent avec $ \delta_0 = \epsilon $.
\\Puis, comme pour le lemme 4.4 de \cite{burq3}, pour prouver la proposition \ref{bilibourgain3}, il suffit d'établir la proposition suivante :
\begin{prop} \label{bilibourgain5}
Pour tout $ b \in ] \frac{1}{2} , 1 ] $ et $ \delta \in ]0, \frac{1}{2} ]$, il existe une constante $ C > 0 $ telle que pour tous $ u,v,M,N $,
\begin{align*}
& ||  \Delta_N ( v )   \Delta_M (u) ||_{L^2 ( [0,1] , L^2 ( \R^d ))  } 
\\ \leq \ & C \times \min ( N,M)^{ \frac{d-2}{2} } \times  \left( \frac{ \min ( N,M)  }{ \max ( N,M )} \right) ^{1/2-\delta } \times  || \Delta_N(v) ||_{  \overline{X}^{0,b}} || \Delta_M(u) ||_{  \overline{X}^{0,b}}.
\end{align*}
\end{prop}
Enfin, pour obtenir les propositions \ref{bilibourgain4} et \ref{bilibourgain5}, il suffit d'utiliser le lemme 2.1 de \cite{burq1} et le théorème \ref{bilis}.
\section{Données initiales aléatoires et espaces de Sobolev}
Le but de cette partie est de montrer que la donnée initiale rendue aléatoire ne permet pas de gagner de dérivée dans $ L^2( \mathds{R}^3 )  $.
\begin{thm} \label{sobolev} Pour tout $ s \geq 0 $, si
\begin{equation*}
 u_0 \notin \overline{H}^s(\R^3)
\end{equation*}
alors
\begin{equation*}
u_0 (\omega,.) \notin  H^s( \R^3 )  \ \omega \ \mbox{ presque surement.}
\end{equation*}
\end{thm}
On désigne par X la loi commune des variables aléatoires $ (g_n)_n $ et pour prouver le théorème, il suffit de considérer les cas où $ X \thicksim  \mathcal{N}_\C (0,1 ) $ ou $ X \backsim \mathcal{B} ( \frac{1}{2} )$.
\\\\Soit $ \chi \in C^\infty_0 ( \R ^3 )  $ telle que $ \chi(x)=1 $ si $ |x| \leq 1 $, $ \chi(x)=0 $ si $ |x| \geq 2 $ et $ 0 \leq \chi \leq 1 $ et définissons $ \sigma_N^2 = \displaystyle{ \sum_{ n \in \N}  }  \ \chi^2 \left(  \frac{\lambda_n^2}{N^2} \right) |c_n|^2 \lambda_n^{2s} $. Comme $ \sigma _N^2  \geq \displaystyle{ \sum_{ \lambda_n \leq N  } } \ |c_n|^2 \lambda_n^{2s} $ alors $ \lim\limits_{\substack{N \to \infty }} \sigma_N^2  = \infty $.
\begin{lem} \label{zygmounde}
Soit $ X $ une variable aléatoire dans $ L^2( \Omega ) $ alors pour tout $ \lambda \geq 0 $,
\begin{equation*}
P( X \geq \lambda E(X) ) \geq (1-\lambda)^2 \frac{E(X)^2}{E(X^2)}.
\end{equation*}
\end{lem}
\textit{Preuve.} Il suffit d'appliquer l'inégalité de Cauchy Schwarz, on pose $ A = \lbrace X \geq \lambda E(X) \rbrace $ et on obtient
\begin{align*}
E(X)=E( X \mathds{1}_A  + X \mathds{1}_{A^c}  ) \leq \sqrt{E(X^2) P(A)}  +  \lambda E(X).
\end{align*}
Donc
\begin{equation*}
(1-\lambda) E(X) \leq \sqrt{E(X^2) P(A)},
\end{equation*} 
et le résultat suit en élevant au carré. \hfill $  \boxtimes  $
\begin{prop} \label{kolmogorov} Pour tout $ s \geq 0 $, on a 
\begin{equation*}
P  \left( \omega \in \Omega /  \sup_{ \ N \in \N^* } \  || \chi \left( \frac{H}{N^2}  \right) u_0^\omega ||_{H^s(\R^3)}  = \infty  \right) =1 \ \mbox{ou} \  = 0.
\end{equation*}
\end{prop}
\textit{Preuve.} Rappelons que $ u_0^\omega = \displaystyle{  \sum _ i } \ X_i(\omega) $ avec $ X_i $ indépendants et $ X_i(\omega) \in H^s( \R^3 ) $ $ \omega$ presque surement. Par conséquent, pour tout $ K \in \N $, on a 
\begin{align*}
 & \sup_{ \ N \in \N^* } \  \bigg| \bigg| \chi \left( \frac{H}{N^2}  \right) u_0^\omega \bigg| \bigg|_{H^s(\R^3)} = \infty  
 \\ \mbox{ si et seulement si}   & \sup_{ \ N \in \N^* } \  \bigg| \bigg| \chi \left( \frac{H}{N^2}  \right) \left( \sum_{i  \geq K }  X_i(\omega) \right) \bigg| \bigg|_{H^s(\R^3)} = \infty,
\end{align*}
donc si nous posons $ F_i = \sigma ( X_i , X_{i+1} , ...  ) $ on a que
\begin{equation*}
\bigg\{   \omega \in \Omega / \sup_{ \ N \in \N^* } \  || \chi \left( \frac{H}{N^2}  \right) u^\omega_0 ||_{H^s(\R^3)}   = \infty \bigg\} \in \underset{K \in N }{ \bigcap }  F_K.
\end{equation*}
Par conséquent $ \left\{   \omega \in \Omega / \displaystyle{ \sup_{ \ N \in \N^* } } \  || \chi \left( \frac{H}{N^2}  \right) u^\omega_0 ||_{H^s}  = \infty  \right\} $ est dans la tribu asymptotique et le lemme est prouvé par la loi du 0-1. \hfill $  \boxtimes  $
\begin{prop} \label{divergence} Pour $ s \geq 0 $, si $ \underset{n  \in \N }{ \sum  }    |c_n|^2 \lambda_n^{2s} = +   \infty $ alors
\begin{equation*}
P \left( \omega \in \Omega /  \sup_{N \in \N^* } || \chi \left( \frac{H}{N^2}  \right) u^\omega_0 ||_{H^s(\R^3)}  = \infty   \right) =1.
\end{equation*}
\end{prop}
\textit{Preuve.}
\begin{equation*}
\mbox{On pose} \ M = \sup_{N \in \N^*}  || \chi \left( \frac{H}{N^2} \right) u_0 ||_{H^s(\R^3)} \ \mbox{et} \ S_N = || \chi \left( \frac{H}{N^2} \right) u_0 ||_{H^s(\R^3)}.
\end{equation*}
En utilisant le lemme \ref{zygmounde}, on obtient
\begin{align*}
P \left( M ^2 \geq \frac{1}{2} \ E( X^2 ) \times C_1^2 \ \sigma_{N}^2  \right) & \geq P \left( S_N^2 \geq \frac{1}{2} \ E( X^2 ) \times C_1^2 \ \sigma_{N}^2  \right) 
\\ & \geq P \left( S_N^2 \geq \frac{1}{2} \ E \left(  ||  \chi \left( \frac{H}{N^2}  \right) u_0  ||^2_{H^s(  \R^3 )}  \right)  \right)
\\ & \geq  \frac{1}{4}   \frac{  E \left(  ||  \chi \left( \frac{H}{N^2} \right) u_0  ||^2_{H^s(  \R^3 )}  \right) ^2   }{E \left(  ||  \chi \left( \frac{H}{N^2} \right) u_0  ||^4_{H^s(  \R^3 )}  \right) }.
\end{align*}
En effet, grâce à la proposition \ref{propre5}, nous avons
\begin{align*} 
& E \left(  ||  \chi \left( \frac{H}{N^2} \right) u _0 ||^2_{H^s(  \R^3 )}  \right) 
\\ \geq \ & E \left( \sum_{n,m}   \  \chi \left( \frac{\lambda_n^2}{N^2} \right) \chi \left( \frac{\lambda_m^2}{N^2} \right) c_n \overline{c_m } g_n ( \omega ) \overline{ g_m (\omega ) }  \int_{\R^2} |\nabla|^s( h_n ) . |\nabla|^s( h_m) \ dx    \right)
\\ \geq \ &  E( |X|^2 ) \times \sum_{n} \chi^2 \left( \frac{\lambda_n^2}{N^2} \right) |c_n|^2 || \nabla^s(h_n) ||^2_{L^2( \R ^3 )}
\\ \geq \ & E( |X|^2 ) \times C_1^2 \sigma_N^2.
\end{align*}
De plus, grâce encore une fois à la proposition \ref{propre5}, on établit
\begin{align*}
& E \left(  ||  \chi \left( \frac{H}{N^2} \right) u _0  ||^4_{H^s(  \R^3 )}  \right)
\\  \leq \ & E \left( \sum_{n,m}   \  \chi \left( \frac{\lambda_n^2}{N^2} \right) \chi \left( \frac{\lambda_m^2}{N^2} \right) c_n \overline{c_m } \ g_n ( \omega ) \overline{ g_m (\omega ) }  \int_{\R^3}   \nabla^s( h_n ) \nabla^s( h_m) \ dx    \right) ^2
\\   & \hspace*{5cm} +  E \left( \sum_{n,m}   \  \chi \left( \frac{\lambda_n^2}{N^2} \right) \chi \left( \frac{\lambda_m^2}{N^2} \right) c_n \overline{c_m } \ g_n ( \omega ) \overline{ g_m (\omega ) }  \right) ^2
\\ \leq \ & E \left( \sum_{n_1,n_2,n_3,n_4}   \  \chi \left( \frac{\lambda_{n_1}^2}{N^2} \right) \chi \left( \frac{\lambda_{n_2}^2}{N^2} \right) \chi \left( \frac{\lambda_{n_3}^2}{N^2} \right) \chi \left( \frac{\lambda_{n_4}^2}{N^2} \right) c_{n_1} \overline{c_{n_2}} c_{n_3} \overline{c_{n_4}}  \right.
\\ & \hspace*{1cm} \times \left. g_{n_1} ( \omega ) \overline{g_{n_2} ( \omega ) } g_{n_3} ( \omega ) \overline{g_{n_4} ( \omega )}    \times  \int_{\R^3}   \nabla^s( h_{n_1} ) \nabla^s( h_{n_2}) \times  \int_{\R^3}   \nabla^s( h_{n_3} ) \nabla^s( h_{n_4}) \right)
\\  & +  E \bigg( \sum_{n_1,n_2,n_3,n_4}   \  \chi \left( \frac{\lambda_{n_1}^2}{N^2} \right) \chi \left( \frac{\lambda_{n_2}^2}{N^2} \right) \chi \left( \frac{\lambda_{n_3}^2}{N^2} \right) \chi \left( \frac{\lambda_{n_4}^2}{N^2} \right)
\\ & \hspace*{6cm} \times  c_{n_1} \overline{c_{n_2}} c_{n_3} \overline{c_{n_4}}   \times  g_{n_1} ( \omega ) \overline{g_{n_2} ( \omega ) } g_{n_3} ( \omega ) \overline{g_{n_4} ( \omega )}   \bigg)
\\ \leq \ & 4 E( |X|^4 ) \times \sum_{n,m} \chi \left( \frac{\lambda_n^2}{N^2} \right)^2 \chi \left( \frac{\lambda_m^2}{N^2}  \right) ^2 | c_n |^2 | c_m |^2 || \nabla^s( h_n ) ||^2_{L^2(\R^3)} || \nabla^s( h_m ) ||^2_{L^2(\R^3)}
\\ & \hspace*{5cm} + 4 E( |X|^4 ) \times \sum_{n,m} \chi \left( \frac{\lambda_n^2}{N^2} \right)^2 \chi \left( \frac{\lambda_m^2}{N^2}  \right) ^2 | c_n |^2 | c_m |^2
\\ \leq \ & 4 E( |X|^4 ) \times C_2^4 \sigma_N^4+ 4 E( |X|^4 ) \times \sigma_N^4 .
\end{align*}
Et, par conséquent,
\begin{equation*}
P \left( M ^2 \geq  \frac{1}{2} \ E(|X|^2) \times C_1 \ \sigma_{N }^2   \right)  \geq \frac{E(|X|^2) ^2}{E(|X|^4)}  \times   \left( \frac{C_1}{2} \right) ^4  \times \left( \frac{1}{C_2^4+1} \right) .
\end{equation*}
Puis en utilisant le théorème de convergence monotone, on trouve
\begin{equation*}
P( M = \infty ) \geq \frac{E(|X |^2) ^2}{E(|X|^4)}  \times \left( \frac{C_1}{2} \right) ^4  \times \left( \frac{1}{C_2^4+1} \right).
\end{equation*}
Et finalement d'après la proposition \ref{kolmogorov}, on a $ P ( M = \infty ) = 1. $ \hfill $  \boxtimes  $
\begin{thm} \label{continuite}
Pour tout $ s \geq 0 $, il existe une constante $ C > 0 $ telle que pour tout $ N \in \N ^* $ et $ u \in H^s(\R^3) $,
 \begin{equation*}
|| \chi \left(  \frac{H}{N^2}  \right) u ||_{H^s(  \R^3 )} \leq C ||u||_{H^s(\R^3)}.
\end{equation*}
\end{thm}
Le théorème \ref{continuite} et la proposition \ref{divergence} impliquent le théorème \ref{sobolev}.
\\\\En effet, si nous supposons que $ u_0^\omega \in H^s( \R^3 ) $ $ \omega $ presque surement alors par la proposition \ref{divergence}, on obtient 
\begin{equation*}
\sup_{N \in \N^*} ||  \chi \left(   \frac{H}{N^2}  \right) u_0^\omega ||_{H^s(  \R ^3 )} \leq C ||u_0^\omega ||_{  H^s(  \R^3 ) },
\end{equation*}
puis
\begin{equation*}
\sup_{N \in \N^*} ||  \chi \left(   \frac{H}{N^2}  \right) u_0^\omega ||_{H^s(  \R ^3 )} < \infty \ \omega \mbox{ presque surement.} 
\end{equation*}
Ce résultat contredit la proposition \ref{divergence} et finalement il suffit de prouver le théorème \ref{continuite}.
\\\\ \textit{Preuve du théorème \ref{continuite}.} En utilisant (\ref{semi-classique1}) et (\ref{semi-classique2}), il suffit de montrer que:
\\ $ \forall s \geq 0 , \ \exists \ C > 0 $ et $ h_0 $ tels que $ \forall \ 0 < h \leq h_0  , \ \forall u \in H^s(\R^3) $
 \begin{equation} \label{continuitebis}
|| \  \chi (  x^2 + (hD)^2   )  u ||_{H^s(  \R^3 )} \leq C ||u||_{H^s(\R^3)}.
\end{equation}
En effet, 
\begin{align*}
& ||  \chi \left(   \frac{H}{N^2}  \right) u ||_{H^s(  \R ^3 )} 
\\  \leq \ & || \ [ \chi (  x^2 + (hD)^2   )  ] \mathfrak{u}  ( \sqrt{h} .  ) ||_{H^s(  \R^3 )}
\\  \leq \ & || \ [ \chi (  x^2 + (hD)^2   )  ] \mathfrak{u}  ( \sqrt{h} .  ) ||_{L^2(  \R^3 )} + ||  \nabla^s \ [ \chi (  x^2 + (hD)^2   )  ] \mathfrak{u}  ( \sqrt{h} .  ) ||_{L^2(  \R^3 )}
\\ \leq \ & h^{-3/4} || \mathfrak{u} ||_{L^2( \R^3)} + h^{s/2 -3/4} || \ [ \chi (  x^2 + (hD)^2   )  ] \mathfrak{u}   ||_{H^s(  \R^3 )}
\\ \leq \ & ||u||_{L^2(\R^3)}+ h^{s/2 -3/4} || \mathfrak{u}  ||_{H^s(   \R^3  ) }
\\ \leq \ & ||u||_{  H^s( \R^3 ) }.
\end{align*}
Par interpolation, on peut limiter la preuve au cas où $ s $ est un entier. Grâce à la proposition \ref{approx} (avec $ N=s $), il existe une constante $  C > 0 $ telle que pour tout $ h \in ]0,1] $ et $ u \in L^2(\R^3) $, on a 
 \begin{equation*} 
|| \chi (  x^2 + (hD)^2    ) u  - \sum_{j=0}^N h^j Op_h (     \Psi_j(x,\xi)   ) u ||_{H^s(  \R^3 )} \leq C ||u||_{L^2(\R^3)},
\end{equation*}
avec $ Supp ( \Psi_j(x,\xi) ) \subset \left( (x,\xi) / x^2 + \xi^2 \in Supp ( \chi ) \right)  $.
\\\\Ainsi, pour obtenir (\ref{continuitebis}), il suffit d'obtenir que pour tout $ s \geq 0 $, il existe deux constantes $ C> 0 $ et $ h_0 \geq 1 $ telles que pour tout $ h \in ]0,h_0] $ et $ u \in H^s(\R^3) $,
 \begin{equation} \label{continuitetris}
|| Op_h (     \Psi_j(x,\xi)   ) u ||_{H^s(  \R^3 )} \leq C ||u||_{H^s(\R^3)}.
\end{equation}
Enfin, pour établir (\ref{continuitetris}), il suffit d'utiliser le théorème \ref{operateur} et de remarquer que
\begin{equation*}
( x , \xi  ) \longrightarrow \chi( x^2 + \xi ^2 ) \in S^0 \mbox{ and }  ( x , \xi  ) \longrightarrow \Psi_j ( x , \xi ) \in S^0. \ \ \ \boxtimes
\end{equation*}
Pour terminer cette partie, on évalue la norme Sobolev de la donnée initiale. Cela permettra d'établir que les théorèmes \ref{thm1} et \ref{thm3} seront valides pour pour des équations sur-critiques avec des données initiales "grandes". 
\begin{prop} Soit $ \sigma \geq 0, \  u_0 \in \overline{H}^\sigma(\R^3) $ et $ s \geq  \sigma $. Supposons que pour tout $ n \in \N $,
\begin{equation*}
  \lambda_n^{2s} \ |c_n|^2 \leq 1 
\end{equation*}
alors pour tout $ t \geq 0 $,
\begin{align*}
\mu \left( u_0 / || \chi \left(   \frac{H}{N^2}    \right)  u_0 ||_{  \overline{H}^s(\R^3) } \leq t   \right) & \leq 
\left\{
    \begin{array}{ll}
       e^{ t ^2 -   \frac{1}{2} || \chi \left(   \frac{H}{N^2}    \right)  u_0 ||^2_{\overline{H}^s(\R^3) } }  \mbox{ dans le cas Gaussien,} \\
      e^{ t ^2 -   \epsilon^2 || \chi \left(   \frac{H}{N^2}    \right)  u_0 ||^2_{\overline{H}^s(\R^3) } } \mbox{ dans le cas Bernoulli.}
    \end{array}
\right.
\end{align*}
\end{prop}
\textit{Preuve.} Effectuons la preuve dans le cas Gaussien. En utilisant que $ -\ln(1+u) \leq -  \frac{u}{2} $ pour tout $ u \in [0,1] $ et l'inégalité de Markov, on obtient
\begin{align*}
\mu \left(  u_0 \in \overline{H}^\sigma(\R^3) / || \chi \left(   \frac{H}{N^2}    \right)  u_0 ||_{   \overline{H}^s(\R^3) } \leq t   \right) & = P \left(   \omega \in \Omega /  e^{-|| \chi \left(   \frac{H}{N^2}    \right)  u^\omega_0 ||^2_{  \overline{H}^s(\R^3) } }  \geq e^{-t^2 }   \right)
\\ & \leq e^{t^2} E \left( e^{-|| \chi \left(   \frac{H}{N^2}    \right)  u_0 ||^2_{  \overline{H}^s(\R^3) } }  \right) 
\\ & \leq e^{t^2} \prod_{ n \in \N   } E \left( e^{- \chi^2 \left(   \frac{\lambda_n^2}{N^2}    \right)  \lambda_n^{2s} | c_n|^2 | X |^2 }   \right)
\\ & \leq e^{t^2} \prod_{ n \in \N   } \left( \frac{1}{1+\chi^2 (   \frac{\lambda_n^2}{N^2}    )  \lambda_n^{2s} | c_n|^2} \right)
\\ & \leq e^{t^2} \prod_{ n \in \N   }  \left( e^{  - \frac{1}{2}\chi^2 \left(   \frac{\lambda_n^2}{N^2}    \right)  \lambda_n^{2s} | c_n|^2 } \right)
\\ & \leq e^{ t^2 - \frac{1}{2}  || \chi \left(   \frac{H}{N^2}    \right)  u_0 ||^2_{\overline{H}^s(\R^3) }  }. \hspace*{2cm} \boxtimes 
\end{align*}
\textit{Remarques :} 1- Par exemple, si $ u_0 \notin \overline{H}^s(\R^3) $ et $ \lambda_n^{2s} \ |c_n|^2 \leq 1, \ \forall n \in \mathds{N} $, on obtient pour tout $ t \geq 0 $, 
\begin{equation*}
\lim_{N \rightarrow \infty} \mu \left( u_0 \in \overline{H}^\sigma(\R^3) / || \chi \left(   \frac{H}{N^2}    \right)  u_0 ||_{  \overline{H}^s(\R^d) } \leq t   \right) = 0 .
\end{equation*}
Encore une fois, cela signifie bien que la norme Sobolev de la donnée initiale n'est pas "petite".
\\2- Par exemple, pour N fixé, on peut choisir $ c_n = \frac{1}{N^{s-\epsilon}} $ pour n vérifiant $  \lambda_n \sim N $ et 0 sinon et obtenir pour $ t \geq 0$ ,
\begin{equation*}
\mu \left( u_0 \in \overline{H}^\sigma(\R^3) / || \chi \left(   \frac{H}{N^2}    \right)  u_0 ||_{  \overline{H}^s(\R^d) } \leq t    \right)  \leq \exp  \left( t^2 - \frac{N^\epsilon}{2}  \right),
 \end{equation*}
 alors que $ ||u_0||_{ \overline{H}^\sigma ( \R^3  )   } = N^{\sigma-s+\epsilon} \ll 1 $.
\section{L'argument de point fixe}
Dans cette partie, on établit les estimées qui vont nous servir à appliquer un théorème de point fixe. $ \psi $ désigne une fonction de $ C^\infty_0( \R ) $ égale à 1 sur $ [ - \frac{\pi}{4}, \frac{\pi}{4} ] $ et supportée dans $ [ - 2 \pi, 2 \pi  ] $.
\begin{prop} Il existe $ b' < \frac{1}{2} $ tel que pour tous $ b > \frac{1}{2} $ et $ s > \frac{1}{2} $, on ait l'existence de deux constantes $ C>0$ et $\kappa >0 $ telles que pour tout $ v \in \overline{X}^{s,b} $ et tous $ N_1 \geq N_2 \geq N_3 $,
\begin{equation} \label{estim1}
|| \Delta_{N_1}(v) \Delta _{N_2}(v) \Delta_{N_3} (v) ||_{\overline{X}^{s,-b'}} \leq C   N_1^{-\kappa} ||v||_{\overline{X}^{s,b}}^3.
\end{equation}
\end{prop}
\textit{Preuve.} Par dualité, il suffit de montrer qu'il existe une constante $\delta > 0 $ telle que  
\begin{equation*}
\int_{\R * \R^3} \Delta_{N_1}(v) \Delta _{N_2}(v) \Delta_{N_3} (v) \Delta_M(w) \leq C  N_1^{-\kappa} M ^{-\delta} ||v||_{\overline{X}^{s,b}}^3 || w ||_{ \overline{X}^{-s,b'} }.
\end{equation*}
Grâce au le lemme \ref{cas facile}, nous pouvons nous ramener au cas où $ M \leq N_1 ^ {1+ \delta } $.
\\\textbf{Cas $N_3 \leq M \leq N_1^{1+\delta}$:} En utilisant le théorème \ref{bilibourgain1}, on obtient
\begin{align*}
& \int_{\R * \R^3} \Delta_{N_1}(v) \Delta _{N_2}(v) \Delta_{N_3} (v) \Delta_M(w)
\\  \leq \ & || \Delta_{N_1}(v) \Delta _{N_2}(v) ||_ {L^2(\R, L^2(\R^3) ) } \times || \Delta_{M}(w) \Delta _{N_3}(v) ||_ {L^2(\R , L^2(\R^3)) }
\\  \leq \ &  (N_2 N_3)^{1/2+\delta}  (\frac{N_2}{N_1})^{1/2-\delta}  (\frac{N_3}{M})^{1/2-\delta}  
\\ & \hspace*{3cm} \times || \Delta_{N_1}(v) ||_{ \overline{X}^{0,b} } || \Delta_{N_2}(v) ||_{ \overline{X}^{0,b} } || \Delta_{N_3}(v) ||_{ \overline{X}^{0,b} } || \Delta_{M}(w) ||_{ \overline{X}^{0,b'} }
\\  \leq \ &  (\frac{N_2}{N_1})^{1/2-\delta} (\frac{N_3}{M})^{1/2-\delta} (N_2 N_3)^{1/2+ \delta-s} (   \frac{M}{N_1} )^{s} \times || v ||_{ \overline{X}^{s,b} }^3 || w ||_{ \overline{X}^{-s,b'} }
\\  \leq \ &  (N_2 N_3)^{1-s}  M ^{-1/2+\delta} N_1^{-1/2+(1+s)\delta} \times || v ||_{ \overline{X}^{s,b} }^3 || w ||_{ \overline{X}^{-s,b'} }
\\  \leq \ & M^{-\delta} M^{1/2-s+(1+s)\delta}   N_1^{1/2-s+2\delta} \times || v ||_{ \overline{X}^{s,b} }^3 || w ||_{ \overline{X}^{-s,b'} }
\\  \leq \ & M^{-\delta} N_1^{1-2s+ (3+s)\delta} \times || v ||_{ \overline{X}^{s,b} }^3 || w ||_{ \overline{X}^{-s,b'} }.
\end{align*}
\textbf{Cas $ M \leq N_3 $:} En utilisant le théorème \ref{bilibourgain1}, on établit
\begin{align*}
& \int_{\R * \R^3} \Delta_{N_1}(v) \Delta _{N_2}(v) \Delta_{N_3} (v) \Delta_M(w)
\\ \leq \ & || \Delta_{N_1}(v) \Delta _{N_2}(v) ||_ {L^2(\R , L^2(\R^3) ) } \times  || \Delta_{M}(w) \Delta _{N_3}(v) ||_ {L^2(\R , L^2(\R^3) ) }
\\ \leq \ & (N_2 M )^{1/2+\delta}  (\frac{N_2}{N_1})^{1/2-\delta} (\frac{M}{N_3})^{1/2-\delta} 
\\ & \hspace*{3cm} \times || \Delta_{N_1}(v) ||_{ \overline{X}^{0,b} } || \Delta_{N_2}(v) ||_{ \overline{X}^{0,b} } || \Delta_{N_3}(v) ||_{ \overline{X}^{0,b} } || \Delta_{M}(w) ||_{ \overline{X}^{0,b'} }
\\ \leq \ & N_2 M (\frac{1}{N_3N_1})^{1/2-\delta}  (   \frac{M}{N_1N_2 N_3} )^{s} \times || v ||_{ \overline{X}^{s,b} }^3 || w ||_{ \overline{X}^{-s,b'} }
\\ \leq \ & M^{-\delta}  N_2^{1-s} N_3^{1/2+2 \delta} \times (   \frac{1}{N_1} )^{1/2+s-\delta}   \times || v ||_{ \overline{X}^{s,b} }^3 || w ||_{ \overline{X}^{-s,b'} }
\\  \leq \ &  M^{-\delta}  N_1^{1-2s+3 \delta}   \times  || v ||_{ \overline{X}^{s,b} }^3 || w ||_{ \overline{X}^{-s,b'} }. \hspace*{7.5cm} \boxtimes 
\end{align*}
\begin{prop} Il existe $ b' < \frac{1}{2} $ tel que pour tous $ b > \frac{1}{2} $ et $ s > \frac{1}{2} $, on ait l'existence de deux constantes $ C,\kappa >0 $ telles que si pour un certain $ \lambda > 0 $, on a pour tout $ N $,
\begin{equation*}
|| \Delta_N( e^{itH} u_0 ) ||_{ L^2 (  [-2 \pi, 2 \pi], L^\infty ( \R^3 )) } \leq \lambda N^{-1/6}
\end{equation*}
alors pour tout $ v \in \overline{X}^{s,b} $ et tous $ N_1 \geq N_2 \geq N_3 $,
\begin{equation} \label{estim2}
|| \Delta_{N_1}(v) \Delta _{N_2}(v) \Delta_{N_3} (\psi(t) e^{-itH}u_0) ||_{\overline{X}^{s,-b'}} \leq C N_1^{-\kappa} \times ( ||v||_{\overline{X}^{s,b}}^3 + \lambda ^3 ),
\end{equation}
\begin{equation} \label{estim3}
|| \Delta_{N_1}(v) \Delta _{N_2}(\psi(t) e^{-itH}u_0) \Delta_{N_3} (v) ||_{\overline{X}^{s,-b'}} \leq C N_1^{-\kappa} \times ( ||v||_{\overline{X}^{s,b}}^3 + \lambda ^3 ).
\end{equation}
\end{prop}
\textit{Preuve.} On montre (\ref{estim2}), la preuve de (\ref{estim3}) étant similaire.
\\Par dualité, il suffit de montrer qu'il existe une constante $ \delta > 0 $ telle que 
\begin{equation*}
\int_{\R * \R^3}\Delta_{N_1}(v) \Delta _{N_2}(v) \Delta_{N_3} (\psi(t) e^{-itH}u_0) \Delta_M(w) \leq C N_1^{-\kappa} M^{-\delta} \times ( ||v||_{\overline{X}^{s,b}}^3 + \lambda ^3 ) \times ||w ||_{  \overline{X}^{-s,b'} }.
\end{equation*}
Grâce au lemme \ref{cas facile}, nous pouvons nous ramener au cas où $ M \leq N_1^{1+\delta} $.
\\\textbf{Cas $N_2 \leq M \leq N_1^{1+\delta}$: } En utilisant le théorème \ref{bilibourgain1} et la proposition \ref{bourgain1}, on obtient
\begin{align*}
& \int_{\R * \R^3}\Delta_{N_1}(v) \Delta _{N_2}(v) \Delta_{N_3} (\psi(t) e^{-itH}u_0) \Delta_M(w)
 \\  \leq \ & ||\Delta_{N_2}(v) \Delta _{M}(w)  ||_{L^2( \R, L^2(\R^3))} \times   ||\Delta_{N_3} (\psi(t) e^{-itH}u_0) ||_{L^2( [-2 \pi, 2 \pi ] ,  L^\infty ( \R^3 )) } 
 \\ & \hspace*{8cm} \times  ||\Delta_{N_1}(v)  ||_{L^\infty( [-2 \pi, 2 \pi ], L^2(\R^3))}
\\  \leq \ &  N_2^{1/2+\delta} (\frac{N_2}{M})^{1/2-\delta}   
\\ & \hspace*{1.3cm} \times   ||  \Delta_{N_1}(v) ||_{  \overline{X}^{0,b} } || \Delta _{N_2}(v) ||_{  \overline{X}^{0,b} }  ||\Delta_{N_3} ( e^{itH}u_0) ||_{L^2( [-2 \pi, 2 \pi] , L^\infty(\R^3))} ||\Delta_{M}(w)  ||_{  \overline{X}^{0,b'} }
\\  \leq  \ &  N_2^{1/2+\delta} (\frac{N_2}{M})^{1/2-\delta} (\frac{M}{N_1N_2})^{s}   N_3^{-1/6} \times  \lambda ||  v ||^2_{  \overline{X}^{s,b} }    ||w  ||_{ \overline{X}^{-s,b'} }
\\  \leq  \ & N_2^{1-s}  M^{s-1/2+\delta}  N_1^{-s}  N_3^{-1/6} \times \lambda ||  v ||^2_{  \overline{X}^{s,b} }   ||w  ||_{ \overline{X}^{-s,b'} }
\\  \leq \ &   M^{-\delta}   N_2^{1-s}  N_1^{-1/2+(3+s) \delta}  N_3^{-1/6} \times \lambda ||  v ||^2_{  \overline{X}^{s,b} }   ||w  ||_{ \overline{X}^{-s,b'} }
\\  \leq  \ &   N_1^{1/2-s+(3+s)\delta} \times \lambda ||  v ||^2_{  \overline{X}^{s,b} }   ||w  ||_{ \overline{X}^{-s,b'} }.
\end{align*}
\textbf{Cas $ M \leq N_2 $:} En utilisant le théorème \ref{bilibourgain1} et la proposition \ref{bourgain1}, on établit
\begin{align*}
& \int_{\R * \R^3}\Delta_{N_1}(v) \Delta _{N_2}(v) \Delta_{N_3} (\psi(t) e^{-itH}u_0) \Delta_M(w)
 \\ \leq \ & ||\Delta_{N_2}(v) \Delta _{M}(w)  ||_{L^2( \R, L^2(\R^3))} \times  ||\Delta_{N_3} (\psi(t) e^{-itH}u_0) ||_{L^2( [-2 \pi, 2 \pi ] ,  L^\infty ( \R^3 )) } 
 \\ & \hspace*{8cm} \times  ||\Delta_{N_1}(v)  ||_{L^\infty( [-2 \pi, 2 \pi ], L^2(\R^3))}
\\  \leq  \ & M^{1/2+\delta} (\frac{M}{N_2})^{1/2-\delta}   
\\ &\hspace*{1.3cm}  \times    ||  \Delta_{N_1}(v) ||_{  \overline{X}^{0,b} } || \Delta _{N_2}(v) ||_{  \overline{X}^{0,b} }  ||\Delta_{N_3} ( e^{itH}u_0) ||_{L^2( [-2 \pi, 2 \pi] , L^\infty(\R^3))} ||\Delta_{M}(w)  ||_{  \overline{X}^{0,b'} }
\\  \leq  \ & M^{1/2+\delta}  (\frac{M}{N_2})^{1/2-\delta} (\frac{M}{N_1N_2})^{s}   N_3^{-1/6} \times \lambda ||  v ||^2_{  \overline{X}^{s,b} }    ||w  ||_{ \overline{X}^{-s,b'} }
\\  \leq \ &  N_2^{-1/2-s+\delta} M^{1+s} N_1^{-s} N_3^{-1/6} \times \lambda ||  v ||^2_{  \overline{X}^{s,b} }   ||w  ||_{ \overline{X}^{-s,b'} }
\\  \leq  \ &  N_2^{1/2+\delta}  N_1^{-s} N_3^{-1/6} \times \lambda ||  v ||^2_{  \overline{X}^{s,b} }   ||w  ||_{ \overline{X}^{-s,b'} }
\\  \leq \ &   M^{-\delta }   N_1^{1/2-s+2\delta} \times \lambda ||  v ||^2_{  \overline{X}^{s,b} }   ||w  ||_{ \overline{X}^{-s,b'} }.  \hspace*{7cm} \boxtimes 
\end{align*}
\begin{prop} Il existe $ b' < \frac{1}{2} $ tel que pour tous $ b > \frac{1}{2} $ et $ s > \frac{1}{2} $, on ait l'existence d'une constante $ C >0 $ telle que si pour un certain $ \lambda > 0 $, on a pour tout $ N $,
\begin{align*}
|| \Delta_N( e^{itH} u_0 )  ||_{ L^4 (  [-2 \pi, 2 \pi], L^\infty ( \R^3 )) } \leq \lambda N^{-1/6}\ \mbox{et} \ ||  \ [ e^{itH} u_0 ]^2 \  ||_{ L^4 (  [-2 \pi, 2 \pi], H^s ( \R^3 )) } \leq \lambda ^2 
\end{align*}
alors pour tout $ v \in \overline{X}^{s,b} $,
\begin{equation} \label{estim4}
|| v* \psi(t) e^{-itH}u_0  * \psi(t) e^{-itH}u_0 ||_{\overline{X}^{s,-b'}} \leq C  ( ||v||_{\overline{X}^{s,b}}^3 + \lambda ^3 ).
\end{equation}
\end{prop}
\textit{Preuve.} En utilisant les propositions \ref{bourgain3} et  \ref{bourgain1}, on trouve
\begin{align*}
  &  || v * \psi(t) e^{-itH}u_0 * \psi(t) e^{-itH}u_0 ||_{\overline{X}^{s,-b'}}  
  \\  \leq \ & || v * \psi(t) e^{-itH}u_0 * \psi(t) e^{-itH}u_0 ||_{L^{1+\delta}( \R , \overline{H}^{s}(\R^3))} 
\\ \leq \ & || v * \psi(t) e^{-itH}u_0 * \psi(t) e^{-itH}u_0 ||_{L^{1+\delta}( [-2\pi,2\pi] , \overline{H}^{s}(\R^3))} 
\\  \leq \ &  ||v||_{L^\infty ( [-2 \pi, 2 \pi ], \overline{H}^s(\R^3))} \times ||  e^{itH} u_0 ||^2_{L^4 ( [-2 \pi, 2 \pi], L^\infty(\R^3))} 
\\  & \hspace*{3cm} +  ||v||_{L^2 ( [-2 \pi, 2 \pi ], L^\infty(\R^3))} \times ||  \ [ e^{itH} u_0 ] ^2 \  ||_{L^4 ( [-2 \pi, 2 \pi],  \overline{H}^s(\R^3))} 
\\  \leq  \ & \lambda ^2 ||v||_{ \overline{ X } ^{s,b}} +  \lambda^2 ||v||_{L^2 ( \R,   \overline{W}^{s,6} (\R^3))}  
\\  \leq \ &  C  ( ||v||_{\overline{X}^{s,b}}^3 + \lambda ^3 ). & \boxtimes
\end{align*}
\begin{prop} 
Il existe $ b' < \frac{1}{2} $ tel que pour tous $ b > \frac{1}{2} $ et $ s > \frac{1}{2} $, on ait l'existence d'une constante $ C >0 $ telle que si pour un certain $ \lambda > 0 $, on a 
\begin{align*}
|| \ [ e^{itH} u_0 ]^3 ||_{L^4 ([-2 \pi , 2 \pi], \overline{ H } ^{s} (\R^3)) } \leq \lambda^3 
\end{align*}
alors pour tout $ v \in \overline{ X} ^{s,b}$,
\begin{equation} \label{estim5}
|| \psi(t) e^{-itH}u_0 * \psi(t) e^{-itH}u_0 * \psi(t) e^{-itH}u_0  ||_{\overline{X}^{s,-b'}} \leq C \lambda ^3.
\end{equation}
\end{prop} 
\textit{Preuve.} En utilisant la proposition \ref{bourgain3}, on établit
\begin{align*}
& || \psi(t) e^{-itH}u_0 * \psi(t) e^{-itH}u_0  * \psi(t) e^{-itH}u_0  ||_{\overline{X}^{s,-b'}} 
\\ \leq \ & || \psi(t) e^{-itH}u_0 * \psi(t) e^{-itH}u_0  * \psi(t) e^{-itH}u_0  ||_{L^{1+\delta} ( \R , \overline{ H}^{s}(\R^3))}
\\  \leq \ & C || \ [ e^{itH}u_0 ]^3  \  ||_{L^4 ( [-2 \pi , 2 \pi ] , \overline{ H } ^{s}(\R^3))}
\\  \leq \ & C \lambda ^3. & \boxtimes 
\end{align*}
\begin{prop} \label{Rfixé}
Il existe $ b' < \frac{1}{2} $ tel que pour tous $ b > \frac{1}{2} $ et $ s \in ] \frac{1}{2} , 1 [  $, on ait l'existence de deux constantes $ C>0, \ \kappa >0 $ et d'un réel $ R \in [ 2,\infty[ $ tels que si pour un certain $ \lambda > 0 $, on a pour tout N,
\begin{equation*}
\left\{
    \begin{array}{ll}
        ||  u_0 ||_{L^2(\R^3)) } \leq \lambda, \\
       || \Delta_N(e^{itH} u_0 ) ||_{L^4( [-2\pi,2\pi] , L^\infty (\R^3)) } \leq \lambda  N^{-1/6}, \\
       || \Delta_N( e^{itH} u_0 ) ||_{L^R( [-2\pi,2\pi] , \overline{W}^{s,4} (\R^3)) } \leq \lambda  N^{s-1/4},
    \end{array}
\right.
\end{equation*}
alors pour tout $ v \in \overline{ X}^{s,b}$ et tous $ N_1 \geq N_2 \geq N_3 $,
\begin{equation} \label{estim6} 
|| \Delta_{N_1}(\psi(t) e^{-itH}u_0) \Delta _{N_2}(v) \Delta_{N_3} (v) ||_{\overline{X}^{s,-b'}} \leq C N_1^{-\kappa} \times (  \lambda ^3  + || v ||_{  \overline{ X } ^{s,b}}^3  ).
\end{equation}
\end{prop} 
\textit{Preuve.} Soit $ \delta > 0  $ assez petit, fixé par la suite.
\\\textbf{Cas $  N_1 \geq (N_2 N_3 )^{  \frac{1-s}{1-s-4 \delta} } $:} Par dualité, il suffit d'établir
\begin{equation*}
\int_{\R * \R^3} \Delta_{N_1}(\psi(t) e^{itH}u_0) \Delta _{N_2}(v) \Delta_{N_3} (v) \Delta_M(w) \leq C  N_1^{-\kappa} M ^{-\delta} \times  ||w||_{\overline{X}^{-s,b'}} \times ( \lambda ^3 +  || v ||^3_{ \overline{X}^{s,b} } ).
\end{equation*}
Grâce au lemme \ref{cas facile}, on peut se ramener au cas où $ M \leq N_1^{1+\delta} $.
\\Si $ N_3 \leq M $ alors en utilisant les théorèmes \ref{bilibourgain1} et \ref{bilibourgain2}, on obtient
\begin{align*}
& \int_{\R * \R^3} \Delta_{N_1}(\psi(t) e^{-itH}u_0) \Delta _{N_2}(v) \Delta_{N_3} (v) \Delta_M(w) 
\\  \leq & \ || \Delta_{N_1}(\psi(t) e^{-itH}u_0) \Delta _{N_2}(v) ||_{L^2( \R , L ^2(\R^3)) } \times ||  \Delta_{N_3} (v) \Delta_M(w) ||_{L^2( \R , L ^2(\R^3)) }
\\  \leq  & \ (N_2 N_3)^{1/2+\delta} (\frac{N_2}{N_1})^{1/2-\delta} (\frac{N_3}{M})^{1/2-\delta}
\\ & \hspace*{3cm} \times || \Delta_{N_1}(u_0) ||_{L^2} || \Delta_{N_2}(v)||_{ \overline{ X}^{0,b}} || \Delta_{N_3}(v)||_{ \overline{ X}^{0,b}} || \Delta_{M}(w)||_{ \overline{ X}^{0,b'}}
\\  \leq & \ (N_2 N_3)^{1/2+\delta} (\frac{N_2}{N_1})^{1/2-\delta} (\frac{N_3}{M})^{1/2-\delta} (\frac{M}{N_2N_3})^s 
\\ & \hspace*{3cm} \times \lambda || \Delta_{N_2}(v)||_{ \overline{ X}^{s,b}} || \Delta_{N_3}(v)||_{ \overline{ X}^{s,b}} || \Delta_{M}(w)||_{ \overline{ X}^{-s,b'}}
\\  \leq & \ M ^{-\delta} N_1^{-\delta} N_1^{-1+s+4\delta} (N_2 N_3)^{1-s} \times  \lambda   || v||^2_{ \overline{ X}^{s,b}} || w ||_{ \overline{ X}^{-s,b'}}
\\  \leq & \ M ^{-\delta} N_1^{-\delta} \times  \lambda   || v||^2_{ \overline{ X}^{s,b}} || w ||_{ \overline{ X}^{-s,b'}}.
\end{align*}
Puis, si $ N_3 \geq M $, en utilisant les théorèmes \ref{bilibourgain1} et \ref{bilibourgain2}, on trouve
\begin{align*}
& \int_{\R * \R^3} \Delta_{N_1}(\psi(t) e^{-itH}u_0) \Delta _{N_2}(v) \Delta_{N_3} (v) \Delta_M(w) 
\\  \leq & \ || \Delta_{N_1}(\psi(t) e^{-itH}u_0) \Delta _{N_2}(v) ||_{L^2( \R , L ^2(\R^3)) } \times ||  \Delta_{N_3} (v) \Delta_M(w) ||_{L^2( \R , L ^2(\R^3)) }
\\  \leq  & \ (N_2 M)^{1/2+\delta} (\frac{N_2}{N_1})^{1/2-\delta} (\frac{M}{N_3})^{1/2-\delta} 
\\ & \hspace*{3cm} \times || \Delta_{N_1}(u_0) ||_{L^2} || \Delta_{N_2}(v)||_{ \overline{ X}^{0,b}} || \Delta_{N_3}(v)||_{ \overline{ X}^{0,b}} || \Delta_{M}(w)||_{ \overline{ X}^{0,b'}}
\\  \leq & \ (N_2 M)^{1/2+\delta} (\frac{N_2}{N_1})^{1/2-\delta} (\frac{M}{N_3})^{1/2-\delta} (\frac{M}{N_2N_3})^s 
\\ & \hspace*{3cm} \times \lambda || \Delta_{N_2}(v)||_{ \overline{ X}^{s,b}} || \Delta_{N_3}(v)||_{ \overline{ X}^{s,b}} || \Delta_{M}(w)||_{ \overline{ X}^{-s,b'}}
\\  \leq & \ M^{-\delta} N_1^{-\delta} N_1^{-1+s+4\delta} (N_2 N_3)^{1-s} \times  \lambda   || v||^2_{ \overline{ X}^{s,b}} || w ||_{ \overline{ X}^{-s,b'}}
\\  \leq & \ M^{-\delta} N_1^{-\delta} \times  \lambda   || v||^2_{ \overline{ X}^{s,b}} || w ||_{ \overline{ X}^{-s,b'}}.
\end{align*}
\textbf{Cas $  N_1 \leq (N_2 N_3 )^{  \frac{1-s}{1-s-4 \delta} } $:} En utilisant les propositions \ref{bourgain3} et  \ref{bourgain1}, on établit
\begin{align*}
& \ || \Delta_{N_1}(\psi(t) e^{-itH}u_0) \Delta _{N_2}(v) \Delta_{N_3} (v) ||_{  \overline{X}^{s,-b'} }
\\  \leq & \  || \Delta_{N_1}(\psi(t) e^{-itH}u_0) \Delta _{N_2}(v) \Delta_{N_3} (v) ||_{  L^{1+\delta}( \R , \overline{H}^s(\R ^3 ) )  }
\\  \leq & \  || \Delta_{N_1}(e^{-itH}u_0) \Delta _{N_2}(v) \Delta_{N_3} (v) ||_{  L^{1+\delta}( [-2\pi,2\pi] , \overline{H}^s(\R ^3 ) )  }
\\  \leq & \ || \Delta_{N_1}(e^{itH}u_0) ||_{  L^\frac{(1+\delta)(1+2\delta)}{\delta} ( [-2 \pi, 2 \pi ]  ,  \overline{W}^{s,4}(\R^3)) }  \times  \prod_{i=1}^2  || \Delta _{N_i}(v) ||_{  L^{2(1+2\delta)}( \R , L^8(\R ^3 ) )  } 
\\  + & || \Delta_{N_1}(e^{itH}u_0) ||_{  L^4 ( [-2 \pi, 2 \pi ]  , L^\infty (\R^3)) }   || \Delta _{N_2}(v) ||_{  L^\infty( [-2 \pi, 2 \pi ] ,  \overline{H}^s (\R ^3 ) )  } || \Delta_{N_3} (v) ||_{  L^2( [-2 \pi, 2 \pi ] ,  \overline{W}^{s,6} (\R ^3 ) )  }
\\  + &  || \Delta_{N_1}(e^{itH}u_0) ||_{  L^4 ( [-2 \pi, 2 \pi ]  , L^\infty (\R^3)) }   || \Delta_{N_2} (v) ||_{  L^2( [-2 \pi, 2 \pi ] ,  \overline{W}^{s,6} (\R ^3 ) )  }|| \Delta _{N_3}(v) ||_{  L^\infty( [-2 \pi, 2 \pi ] ,  \overline{H}^s (\R ^3 ) )  }
\\  \leq & \ N_1^{s-1/4}  (N_2N_3)^{ \frac{9}{8} - \frac{1}{1+2\delta} -s  } \lambda ||v||^2_{  \overline{X}^{s,b} } +  N_1^{-1/6} \times  \lambda ||v||^2_{  \overline{X}^{s,b} }
\\ \leq & \ (N_2 N_3 )^{  \frac{(1-s)(s-1/4+\delta)}{1-s-4 \delta} } (N_2N_3)^{ \frac{9}{8} - \frac{1}{1+\delta} -s  } N_1^{-\delta} \times \lambda ||v||^2_{  \overline{X}^{s,b} } +  N_1^{-1/6} \lambda ||v||^2_{  \overline{X}^{s,b} },
\end{align*}
avec
\begin{align*}
 \frac{(1-s)(s-1/4+\delta)}{1-s-4 \delta} + \frac{9}{8} - \frac{1}{1+2\delta} -s & = s- \frac{1}{4}+  \delta + \frac{4 \delta (s- \frac{1}{4} + \delta )}{1-s-4 \delta} + \frac{9}{8} - \frac{1}{1+2\delta} -s 
\\ & = \frac{7}{8}+ \frac{4 \delta (s- \frac{1}{4} + \delta  )}{1-s-4 \delta}  - \frac{1}{1+2\delta}
\\ & = - \frac{1}{8} + \circ (\delta) < 0.
 \end{align*}
Et finalement la proposition est démontrée avec $ R = \frac{(1+\delta)(1+2\delta)}{\delta} $. \hfill $ \boxtimes $
\\\begin{dfn}
Soit $ \lambda > 0 $ et définissons $ E_0( \lambda ) $ comme l'ensemble des fonctions $ u_0 \in L^2(\R^3)$ qui vérifient
\begin{align*}\left\{
    \begin{array}{ll}
      ||u_0||_{ L^2(\R^3) } \leq \lambda   
      \\ || \ [e^{itH} u_0]^2 \ ||_{  L^4( [-2\pi,2\pi] , \overline{H}^s (\R^3)) } \leq  \lambda ^2 
       \\ || \ [e^{itH} u_0]^3 \ ||_{ L^4( [-2\pi,2\pi] , \overline{H}^s (\R^3))  } \leq \lambda ^3
       \\  || \Delta_N(e^{itH} u_0 ) ||_{L^4( [-2\pi,2\pi] , L^\infty (\R^3)) } \leq \lambda  N^{-1/6} \ , \  \forall N
       \\      || \Delta_N( e^{itH} u_0 ) ||_{L^R( [-2\pi,2\pi] , \overline{W}^{s,4} (\R^3)) } \leq \lambda  N^{s-1/4}   \ , \ \forall N  
    \end{array} 
\right.
\end{align*}
où R est fixé par la proposition \ref{Rfixé}.
\end{dfn} 
\begin{thm} Soient $ \frac{1}{2} < s < 1 $ et $ K \in \lbrace -1, 1 \rbrace $ alors il existe une constante $ C > 0 $ et un réel $ b > 1/2 $ tels que si $ u_0 \in E_0(\lambda ) $ avec $ \lambda > 0 $ alors pour tout $ v \in \overline{X}^{s,b} $,
 \begin{align*}
\bigg| \bigg|  \psi(t)  \int_0^t   e^{-i(t-s)H} \psi(s) K \cos ( 2s ) |  \psi(s) e^{-isH} u_0 + v |^2 \times (   \psi(s) e^{-isH} u_0 + v ) \ ds   \bigg| \bigg|_{  \overline{X}^{s,b} }  
\\ \leq C \times ( \lambda ^3 +  ||v||^3_{  \overline{X}^{s,b} } ).
\end{align*} 
\end{thm} 
\textit{Preuve.} Pour tout $ b > \frac{1}{2} $, en utilisant les propositions \ref{bourgain4} et \ref{bourgain5}, on trouve
\begin{align*}
& ||  \psi(t)  \int_0^t   e^{-i(t-s)H} K \cos ( 2s ) \psi(s) |  \psi(s) e^{-isH} u_0 + v |^2 \times (   \psi(s) e^{-isH} u_0 + v ) \ ds   ||_{  \overline{X}^{s,b} }  
\\  \leq \  C & || \ K \cos ( 2s ) \psi(s)  |  \psi(s) e^{-isH} u_0 + v |^2 \times (   \psi(s) e^{-isH} u_0 + v )   ||_{  \overline{X}^{s,b-1} }  
\\ \leq \  C & || \ |  \psi(s) e^{-isH} u_0 + v |^2 \times (   \psi(s) e^{-isH} u_0 + v )   ||_{  \overline{X}^{s,b-1} }.
\end{align*}
Puis en utilisant (\ref{estim1}), (\ref{estim2}), (\ref{estim3}), (\ref{estim4}), (\ref{estim5}), (\ref{estim6}), on établit l'existence d'un entier $ b' < \frac{1}{2} $ tel que pour tout $ u_0 \in E_0(\lambda) $,
\begin{align*}
|| \  |  \psi(s) e^{-isH} u_0 + v |^2 \times (   \psi(s) e^{-isH} u_0 + v )   ||_{  \overline{X}^{s,-b'} }  \leq C ( \lambda^3 +   ||v||^3_{  \overline{X}^{s,b} }   ).
\end{align*}
Il suffit alors de choisir $ b = 1 - b' > \frac{1}{2} $ et la proposition est démontrée. \hfill $ \boxtimes $
\begin{thm} \label{point fixe etape 1}
Soient $ \frac{1}{2} < s < 1 $ et $ K \in \lbrace -1, 1 \rbrace $ alors il existe une constante $ C > 0 $ et un réel $ b > 1/2 $ tels que si $ u_0 \in E_0(\lambda ) $ avec $ \lambda > 0 $, alors pour tout $ v \in \overline{X}_T^{s,b} $,
\begin{align*}
\bigg| \bigg|  \psi(t)   \int_0^t   e^{-i(t-s)H}  K \cos (2s) \psi(s) \  |  \psi(s) e^{-isH} u_0 + v |^2 \times (   \psi(s) e^{-isH} u_0 + v ) \ ds   \bigg| \bigg|_{  \overline{X}_T^{s,b} }  
\\ \leq C ( \lambda ^3 +  ||v||^3_{  \overline{X}_T^{s,b} } ).
\end{align*} 
\end{thm} 
\textit{Preuve.} Soit $ w \in \overline{X}^{s,b} $ telle que $ w |_{[-T,T]}=v $ alors
\begin{align*}
& ||  \psi(t)   \int_0^t   e^{-i(t-s)H}  K\cos (2s) \psi(s) |  \psi(s) e^{-isH} u_0 + v |^2 \times (   \psi(s) e^{-isH} u_0 + v ) \ ds   ||_{  \overline{X}_T^{s,b} } 
\\  \leq \ & ||  \psi(t)  \int_0^t   e^{-i(t-s)H} K\cos (2s) \psi(s) |  \psi(s) e^{-isH} u_0 + w |^2 \times (   \psi(s) e^{-isH} u_0 + w ) \ ds   ||_{  \overline{X}^{s,b} } 
\\  \leq \ & C ( \lambda ^3 +  ||w||^3_{  \overline{X}^{s,b} } ) \ \mbox{pour tout} \ w. & \boxtimes 
\end{align*}
De manière similaire, on pourrait démontrer le théorème suivant :
\begin{thm} \label{point fixe etape 2}
Soient $ \frac{1}{2} < s < 1 $ et $ K \in \lbrace -1, 1 \rbrace $ alors il existe une constante $ C > 0 $ et un réel $ b > 1/2 $ tels que si $ u_0 \in E_0(\lambda ) $ avec $ \lambda > 0 $, alors pour tous $ v_1,v_2 \in \overline{X}_T^{s,b} $,
 \begin{align*}
& \bigg| \bigg|  \psi(t)   \int_0^t   e^{-i(t-s)H}  K \cos (2s) \psi(s) |  \psi(s) e^{-isH} u_0 + v_1 |^2 \times (   \psi(s) e^{-isH} u_0 + v_1 ) \ ds  
\\ & \hspace*{0.1cm} - \psi(t)   \int_0^t   e^{-i(t-s)H} K \cos (2s) \psi(s) |  \psi(s) e^{-isH} u_0 + v_2 |^2 \times (   \psi(s) e^{-isH} u_0 + v_2 ) \ ds  \bigg| \bigg|_{  \overline{X}_T^{s,b} }  
\\  \leq \ C  & ||v_1-v_2||_{ \overline{X}^{s,b}_T } \times ( \lambda ^2 +  ||v_1||^2_{  \overline{X}_T^{s,b} } + ||v_2||^2_{  \overline{X}_T^{s,b} } ).
\end{align*} 
\end{thm}
\section{Solutions globales pour l'équation (\ref{schrodinger})}
Dans cette partie, on applique un théorème de point fixe pour établir l'existence de solutions globales pour l'équation (\ref{schrodinger}). On démontre également l'unicité des solutions, ainsi que quelques propriétés qu'elles vérifient comme le scattering.
\\Commençons par établir l'existence, pour cela considérons l'équation suivante: 
\begin{equation} \label{schrodingerH} \tag{NLSH}
  \left\{
      \begin{aligned}
        i \frac{ \partial u }{ \partial t } - H u &= K \cos (2 t) |u|^2 u,
       \\ u(0,x)&=u_0(x).
      \end{aligned}
    \right.
\end{equation}
\begin{thm}  \label{existencebis}
Soit $ \frac{1}{2} < s < 1 $ alors il existe une constante $ C > 0 $ et un réel $ b > 1/2 $ tels que si $ u_0 \in E_0(\lambda ) $ avec $ \lambda <   \frac{1}{ 2  \sqrt{C} } $ alors il existe unique solution à l'équation (\ref{schrodingerH}) avec donnée initiale $ u_0 $ dans l'espace $ e^{-itH}u_0 + B_{\overline{X}^{s,b}_T } (0, \frac{1}{2} \sqrt{  \frac{1}{C} } ) $.
\end{thm}
\textit{Preuve.} On définit 
\begin{equation*}
L : v  \rightarrow  -i \psi(t)  \int_0^t   e^{-i(t-s)H} K \cos (2s) \psi(s)  |  \psi(s) e^{-isH} u_0 + v(s) |^2 ( \psi(s) e^{-isH} u_0 + v(s) ) \ ds,
\end{equation*}
$ u = e^{-itH} u_0 + v $ est l'unique solution de (\ref{schrodingerH}) dans l'espace $ e^{-itH}u_0 + B_{\overline{X}^{s,b}_T}(0 , R ) $ si et seulement v est l'unique point fixe de L dans l'espace $ B_{\overline{X}^{s,b}_T}(0 , R ) $.
\\\\ Selon les théorèmes \ref{point fixe etape 1} et \ref{point fixe etape 2}, il existe une constante $ C > 0 $ telle que 
\begin{align*}
& || L (v)  ||_{  \overline{ X } ^{s,b}_T  } \leq C ( \lambda ^3 + || v||^3 _{ \overline{ X } ^{s,b}_T }), 
\\ & || L (v_1)-L(v_2)  ||_{  \overline{ X } ^{s,b}_T  } \leq C || v_1-v_2|| _{ \overline{ X } ^{s,b}_T } ( \lambda ^2 + || v_1||^2 _{ \overline{ X } ^{s,b}_T }+ || v_2||^2 _{ \overline{ X } ^{s,b}_T } ).
\end{align*}
Ainsi, si $ \lambda <   \frac{1}{2 \sqrt{ C} } $ alors L est une application contractante de l'espace complet \\$ B_{\overline{X}^{s,b}_T } (0, \frac{1}{2} \sqrt{  \frac{1}{C} } ) $ et admet donc un unique point fixe. \hfill $ \boxtimes $
\begin{thm}  \label{existence}
Soit $ \frac{1}{2} < s < 1 $ alors il existe deux constantes $ C,c > 0 $ telles que si $ u_0 \in E_0(\lambda ) $ avec $ \lambda <   \frac{1}{2 \sqrt{C} } $ alors il existe une unique solution globale à l'équation (\ref{schrodinger}) avec donnée initiale $ u_0 $ dans l'espace $ e^{it\Delta}u_0 + B_{X^s } (0, \frac{c}{2} \sqrt{  \frac{1}{C} } ) $.
\\La constante $ C $ est donnée par les théorèmes \ref{point fixe etape 1} et \ref{point fixe etape 2} et la constante c par la proposition \ref{toutvabien}.
\end{thm}
\textit{Preuve.} Soit u donnée par le théorème \ref{existencebis} et définissons
\begin{equation*}
 \tilde{u} (t,x) = \left( \frac{1}{\sqrt{1+4t^2}} \right)  ^{3/2} \times u \left( \frac{1}{2} \arctan(2t) , \frac{x}{\sqrt{1+4t^2} } \right)   \times e^{ \frac{ix^2t}{1+4t^2}  }.
\end{equation*}
D'après le théorème \ref{lechangementdevariable}, comme u est solution de (\ref{schrodingerH}) sur $ ]-T,T[ $ alors $ \tilde{u} $ est solution globale de (\ref{schrodinger}).
\\\\Ainsi, pour obtenir le théorème, il suffit de remarquer que
\begin{equation*}
(e^{it\Delta} u_0)(t,x) = \left( \frac{1}{\sqrt{1+4t^2}} \right) ^{3/2} \times ( e^{-itH}  u_0 ) \left( \frac{1}{2} \arctan(2t) , \frac{x}{\sqrt{1+4t^2} } \right)  \times e^{ \frac{ix^2t}{1+4t^2}  },
\end{equation*}
et d'utiliser la proposition \ref{toutvabien}. \hfill $ \boxtimes $
\\\\L'existence de solutions étant prouvée, on démontre ensuite qu'elles sont uniques.
\begin{thm} \label{unicite}
Soient $ \frac{1}{2} < s < 1 $ et $ u_0 \in E_0(\lambda) $ avec $\lambda > 0 $ alors si $ \tilde{u}_1 $ et $ \tilde{u}_2 $ sont deux solutions de (\ref{schrodinger}) de l'espace $  e^{it\Delta} u_0 + X^s $, alors
\begin{equation*}
 \tilde{u}_1  = \tilde{u}_2 \ \mbox{dans} \ C^0( \R , L^2(\R^3) )
\end{equation*} 
\end{thm}
\textit{Preuve.} Il suffit de montrer que $ \tilde{u}_1(t)  = \tilde{u}_2(t) $ pour $ t \geq 0 $ (on remarque que $ x_i(t) := \tilde{u}_i(-t) $ vérifie $ i \frac{\partial x_i }{ \partial t } - \Delta x_i = -K|x_i|^2 x_i  $ et on pourra faire la même preuve pour obtenir que $ x_1(t)=x_2(t)$ pour $ t\geq 0$, c'est à dire $ \tilde{u}_1(t)  = \tilde{u}_2(t) $ pour $ t \leq 0 $).
\\\\Pour tout $ t \in \R$,
\begin{align*}
\partial_t || \tilde{u}_1(t) - \tilde{u}_2(t) ||^2_{ L^2(\R^3) } & = 2 \Re \left( <  \partial_t(\tilde{u}_1(t)-\tilde{u}_2(t)) , \tilde{u}_1(t)-\tilde{u}_2(t)  >_{L^2(\R^3) \times L^2(\R^3) } \right)
\\ & \leq 2 |  <  | \tilde{u}_1(t)| ^2 \tilde{u}_1(t) - |\tilde{u}_2(t)|^2  \tilde{u}_2(t)  , \tilde{u}_1(t)-\tilde{u}_2(t)  >_{L^2(\R^3) \times L^2(\R^3) } |
\\ & \leq 2 || \tilde{u}_1(t) - \tilde{u}_2(t) ||_{L^2(\R^3)} \times || \ |\tilde{u}_1(t)| ^2 \tilde{u}_1(t) - |\tilde{u}_2(t)|^2 \tilde{u}_2(t) ||_{L^2(\R^3)} 
\\ & \leq 4 ||\tilde{u}_1(t) - \tilde{u}_2(t) ||^2_{L^2(\R^3)} \times (  ||  \tilde{u}_1(t) ||^2_{L^\infty(\R^3)} +   ||  \tilde{u}_2(t) ||^2_{L^\infty(\R^3)}  ).
\end{align*}
Par le lemme de Gronwall, le théorème est démontré si $ || \tilde{u}_1(t) ||^2_{L^\infty(\R^3)} +   ||  \tilde{u}_2(t) ||^2_{L^\infty(\R^3)} \in L^1_{loc}(\R^+) $ car  $ || \tilde{u}_1(0) - \tilde{u}_2(0) ||^2_{ L^2(\R^3) } = 0 $. 
\\\\Le théorème est donc clair puisque grâce à la proposition \ref{toutvabienbis},
\begin{align*}
|| \tilde{u}_i ||_{ L^2(\R, L^\infty(\R^3)  ) } & \leq || e^{it\Delta} u_0 ||_{ L^2(\R, L^\infty(\R^3) )  } +|| \tilde{v}_i ||_{ L^2(\R, L^\infty(\R^3) ) }
\\ & \leq C (     || e^{itH} u_0 ||_{ L^4([-2 \pi,2 \pi ], L^\infty(\R^3) )  } + || \tilde{v}_i ||_{ X^s } )
\\ & \leq C (    \lambda + ||  \tilde{v}_i ||_{ X^s } ). & \boxtimes
\end{align*}
On prouve ensuite que les solutions construites diffusent en $ \infty $ et $ - \infty $. Pour prouver ce résultat, commençons par établir un lemme préliminaire.
\begin{lem} \label{changevariables} 
Pour tout $ t \in \R $ et $ x \in \R ^3 $, on a 
\begin{align*}
& \left( e^{it\Delta} F(\frac{1}{2} \arctan 2 t,.) \right) (  t  , x )
\\  = & \left( \frac{1}{\sqrt{1+4t^2}} \right) ^{3/2} \times ( e^{-itH} F(t,.)) \left( \frac{1}{2} \arctan(2t) , \frac{x}{\sqrt{1+4t^2} } \right)  \times e^{ \frac{ix^2t}{1+4t^2}  }.
\end{align*}
\end{lem}
\textit{Preuve.} Rappelons la formule de Mehler, pour $ t \in ]-T,T[$ et $ x \in \R^3$, on a  
\begin{equation*}
( e^{-itH} f ) (t,x ) =  \left( \frac{1}{2 \pi i  \sin 2 t } \right) ^{  \frac{3}{2} } \int_{\R ^3} e^{ \ \frac{i}{2} \ \left(  \frac{\cos 2t}{\sin 2t }  x^2 - \frac{xy}{\sin 2 t}  +  \frac{\cos 2t}{\sin 2t }  y^2 \right)     } f(y) \ d y.
\end{equation*}
Par conséquent,
\begin{align*}
& \left( \frac{1}{\sqrt{1+4t^2}} \right) ^{3/2} \times ( e^{-itH} F(t,.)) \left( \frac{1}{2} \arctan(2t) , \frac{x}{\sqrt{1+4t^2} } \right)   \times e^{ \frac{ix^2t}{1+4t^2}  }  
\\&  =  \left( \frac{1}{4 \pi i t } \right) ^{  \frac{3}{2} } \times  \int_{\R ^3} e^{ \ \frac{i(x-y)^2}{4t}    } F \left( \frac{1}{2} \arctan(2t) , y \right) \  d y 
\\ & =  \left( e^{it\Delta} F(\frac{1}{2} \arctan 2 t,.)  \right) (   t , x ). & \boxtimes 
\end{align*}
\begin{thm} \label{scattering}
Soit $ \tilde{u} $ l'unique solution de (\ref{schrodinger}) construite dans le théorème \ref{existence}, alors il existe $ L^+ \in \overline{H}^s(\R^3) $ et $ L_- \in \overline{H}^s(\R^3)$ deux fonctions telles que
\begin{align*}
& \lim_{t \rightarrow \infty } || \tilde{u}(t) -e^{it\Delta} u_0 - e^{it\Delta} L^+ ||_{H^s(\R^3)} = 0, 
\\ & \lim_{t \rightarrow - \infty } || \tilde{u}(t) -e^{it\Delta} u_0 - e^{it\Delta} L_- ||_{H^s(\R^3)} = 0.
\end{align*}
\end{thm}
\textit{Preuve.} On a montré que
\begin{align*}
-i \psi(t)   \int_0^t   e^{-i(t-s)H} K \cos (2s) \psi(s)  & |  \psi(s) e^{-isH} u_0 + v(s) |^2 ( \psi(s) e^{-isH} u_0 + v(s) ) \ ds 
\\ & \in \overline{X}^{s,b}_T.
\end{align*}
Ainsi, par le lemme \ref{bourgain7},
\begin{align*}
- i e^{-itH} \int_0^t   e^{isH} K \cos (2s) \psi(s) & |  \psi(s) e^{-isH} u_0 + v(s) |^2 ( \psi(s) e^{-isH} u_0 + v(s) ) \ ds 
\\ & \in C ^0 ( [-T,T]  , \overline{H}^s(\R^3)).
\end{align*}
Et donc, il existe une fonction $ L \in \overline{H}^s(\R^3) $ telle que 
\begin{align*}
& \lim_{ t \rightarrow T } \bigg|\bigg| L- i e^{-itH} \int_0^t   e^{isH} K \cos (2s) \psi(s)  |  \psi(s) e^{-isH} u_0 + v(s) |^2 ( \psi(s) e^{-isH} u_0 + v(s) ) ds  \bigg|\bigg|_{ \overline{H}^s } 
\\  = & \lim_{ t \rightarrow T } \bigg|\bigg|  e^{itH} L - i \int_0^t   e^{isH} K \cos (2s) \psi(s)  |  \psi(s) e^{-isH} u_0 + v(s) |^2 ( \psi(s) e^{-isH} u_0 + v(s) ) ds  \bigg|\bigg|_{ \overline{H}^s } 
\\  = &\lim_{ t \rightarrow T } \bigg|\bigg|   e^{iTH} L  - i \int_0^t   e^{isH} K \cos (2s) \psi(s)  |  \psi(s) e^{-isH} u_0 + v(s) |^2 ( \psi(s) e^{-isH} u_0 + v(s) ) ds  \bigg|\bigg|_{ \overline{H}^s} 
\\ & \hspace*{4cm} = 0 .
\end{align*}
Or, pour $ t \in [-T,T]$,
\begin{align*}
 u(t) = & e^{-itH} u_0 
\\ & -i  e^{-itH} \int_0^t   e^{isH} K \cos (2s) \psi(s)  |  \psi(s) e^{-isH} u_0 + v(s) |^2 ( \psi(s) e^{-isH} u_0 + v(s) ) \ ds.
\end{align*}
Donc, par le lemme \ref{changevariables}, on obtient
\begin{align*}
\tilde{u}(t) & =   e^{it\Delta} u_0 
\\ \ +  e^{it\Delta} & \left[-i\int_0^{ \frac{1}{2} \arctan 2t  }   e^{isH} K \cos (2s) \psi(s)  |  \psi(s) e^{-isH} u_0 + v(s) |^2 ( \psi(s) e^{-isH} u_0 + v(s) ) \ ds \right]
\\ & =   e^{it\Delta} u_0 +  e^{it\Delta}  F( t ),
\end{align*}
avec 
\begin{equation*}
\lim_{t \rightarrow \infty} || F(t)- e^{iTH} L  ||_{ \overline{H}^s(\R^3)  } = 0. 
\end{equation*}
Et le théorème est démontré avec
\begin{equation*}
L^+ =  e^{iTH} L \in \overline{H}^s(\R^3),
\end{equation*}
car
\begin{align*}
\lim_{ t \rightarrow \infty } || e^{it\Delta}  F( t ) - e^{it \Delta} L ^+ ||_{  H^s(\R^3) } & = \lim_{ t \rightarrow \infty } ||  F( t ) - L ^+ ||_{  H^s(\R^3) }
\\ & \leq C \lim_{ t \rightarrow \infty } ||  F( t ) - L ^+ ||_{  \overline{H}^s(\R^3) } = 0. & \boxtimes 
 \end{align*}
Enfin, pour conclure cette partie, on démontre que le flot de l'équation est lipschitzien en un certain sens.
\begin{thm} \label{continuiteinitial}
Soient $ u_0^1 , u_0^2 \in E_0(\lambda) $ avec $ \lambda $ donné par le théorème \ref{existence} et soient $ \tilde{u_1},  \tilde{u_2} $ les solutions de (\ref{schrodinger}), alors il existe une constante $ C_\lambda > 0 $ telle que
\begin{equation*}
|| \tilde{u_1} -\tilde{u_2} ||_{ X^0 } \leq C_\lambda ||u_0^1-u_0^2 ||_{ L^2 (\R^3)}.
\end{equation*}
\end{thm}
\textit{Preuve.} Grâce à la proposition \ref{toutvabientris}, il suffit d'établir que
\begin{equation*}
|| u_1 -u_2 ||_{ \overline{X}^0_T } \leq C_\lambda ||u_0^1-u_0^2 ||_{ L^2 (\R^3)}.
\end{equation*}
Rappelons que
\begin{align*}
& u_1(t) -u_2(t)= e^{-itH}(u_0^1-u_0^2) - i \psi(t) \int_0^t e^{-i(t-s)H} \psi(s) K \cos 2s  \times  (    | \psi(s) e^{-isH}u_0^1 + v_1(s)  |^2
\\ & \hspace*{1cm} * (\psi(s) e^{-isH}u_0^1 + v_1(s) )  - |\psi(s) e^{-isH}u_0^2 + v_2(s) |^2 (\psi(s) e^{-isH}u_0^2 + v_2(s) ) )     ds.  
\end{align*}
Ainsi, en utilisant les estimées de Strichartz et par la proposition \ref{bourgain1}, le fait que
\begin{equation*}
|| v_i ||_{  L^2( [T,T] , L^{\infty} (\R^3)) } \leq || v_i ||_{  L^2( [-T,T] , \overline{W}^{s,6} (\R^3)) }  \leq C || v_i ||_{  \overline{X}^{s,b}_T } \leq  C \lambda,
\end{equation*}
on obtient
\begin{align*}
|| u_1 -u_2 ||_{ \overline{X}^0_T } & \leq C  ( ||u_0^1-u_0^2 ||_{ L ^2 (\R^3)} + ||v_1-v_2 ||_{ L^\infty ( [-T,T ] , L^2 (\R^3) ) } ) 
\\ & \hspace*{1cm} \times \left(  1 +  \sum_{j=1,2}   ||e^{itH} u_0^j ||^2_{  L^2( [-2\pi,2\pi] , L^{\infty} (\R^3))}  + || v_j ||^2_{  L^2( [-T,T] , L^{\infty} (\R^3))}  \right)
\\ &  \leq C  \lambda ^2 ( ||u_0^1-u_0^2 ||_{ L^2 (\R^3)} + ||v_1-v_2 ||_{ L^\infty ( [-T,T ] , L^2 (\R^3) ) } ) .
\end{align*} 
Ainsi, pour prouver le résultat, il suffit de montrer que
\begin{equation*}
||v_1-v_2 ||_{ L^\infty ( [-T,T ] , L^2 (\R^3) ) }  \leq C ||u_0^1-u_0^2 ||_{ L^2 (\R^3)}.
\end{equation*}
Mais comme 
\begin{equation*}
i \partial_t v_i - H v_i = K \cos 2t | e^{-itH} u_0^i + v_i  |^2 * (  e^{-itH} u_0^i + v_i   ),
\end{equation*}
on en déduit
\begin{align*}
  & \hspace*{0.7cm} \partial_t || v_1(t)-v_2(t) ||^2_{ L^2 (\R^3)   }  
  \\ & = 2 \Re \left(  \int_{\R^3}   \partial_t(v_1(t)-v_2(t)) \ . \ \overline{ v_1(t)-v_2(t) }  \right)  
\\ & =  2 \Re ( -i K \cos 2 t *  \int_{\R^3}  (  | e^{-itH} u_0^1 + v_1  |^2  (  e^{-itH} u_0^1 + v_1   ) 
\\ & \hspace*{5cm} - | e^{-itH} u_0^2 + v_2  |^2  (  e^{-itH} u_0^2 + v_2   )     ) \ . \  \overline{ v_1(t)-v_2(t) } \ dx   ) 
\\ & \leq || v_1(t)-v_2(t) ||_{ L^2 (\R^3)   } 
\\ & \hspace*{0.1cm} \times || \  | e^{-itH} u_0^1 + v_1(t)  |^2 * (  e^{-itH} u_0^1 + v_1(t)   ) - | e^{-itH} u_0^2 + v_2(t)  |^2 * (  e^{-itH} u_0^2 + v_2(t)  ) \       ||_{L^2(\R^3)}
\\ & \leq || v_1(t)-v_2(t) ||_{ L^2 (\R^3)   }  \times (     || v_1(t)-v_2(t) ||_{ L^2 (\R^3)   } +  ||u_0^1-u_0^2 ||_{ L^2 (\R^3)}  ) 
\\ & \hspace*{6.5cm} \times  \left( \sum_{j=1,2}   ||e^{itH} u_0^j ||^2_{  L^{\infty} (\R^3)}  + || v_j ||^2_{  L^{\infty} (\R^3)}  \right)
\end{align*}
Ainsi, on a établi
 \begin{align*}
& \partial_t || v_1(t)-v_2(t) ||_{ L^2 (\R^3)   }  \leq  (     || v_1(t)-v_2(t) ||_{ L^2 (\R^3)   } +  ||u_0^1-u_0^2 ||_{ L^2 (\R^3)}  ) 
\\ & \hspace*{7cm} \times  \left( \sum_{j=1,2}   ||e^{itH} u_0^j ||^2_{  L^{\infty} (\R^3)}  + || v_j ||^2_{  L^{\infty} (\R^3)}   \right) .
\end{align*}
Finalement, en utilisant le lemme de Gronwall, on obtient
\begin{align*}
& || v_1-v_2 ||_{ L^\infty( [0,T,]L^2 (\R^3))   }  \leq C ||u_0^1-u_0^2 ||_{ L^2 (\R^3)}    
\\&  \hspace*{5cm} \times e^{  \underset{j=1,2}{\sum} ||e^{itH} u_0^j ||^2_{ L^2( [-2\pi,2\pi] ,  L^{\infty} (\R^3))}  + || v_j ||^2_{ L^2( [-T,T] , L^{\infty} (\R^3) ) }}
\\ & \hspace*{4cm} \leq  C e^{ C \lambda^2} \times  ||u_0^1-u_0^2 ||_{ L^2 (\R^3)}. & \boxtimes 
\end{align*}
\section{Estimation de la régularité de la donnée initiale aléatoire}
Dans cette partie, on estime la régularité de la donnée aléatoire en démontrant des estimée de types grandes déviations. En particulier, on établit que $ u^\omega_0 \in \underset{ n \in \N}{ \bigcup } E_0( n) $ $ \omega $ presque surement. On supposera que $ g_n \sim \mathcal{N}_\C(0,1) $ ou $ g_n \sim \mathcal{B}( \frac{1}{2}  ) $ car dans ce dernier cas, il suffira de remplacer $ u_0 $ par $ \epsilon * u_0 $ pour revenir au cas où $ \frac{ g_n}{\epsilon} \sim \mathcal{B}( \frac{1}{2}  ) $.
\\\\On définit pour $ t > 0 $,
\begin{align*}
\Omega_t =  & \left( \omega \in \Omega  / ||u^\omega_0||_{ \overline{H}^\sigma(\R^3) } \leq t , || \ [e^{itH} u^\omega_0]^2 \ ||_{  L^4( [-2\pi,2\pi] , \overline{H}^s (\R^3)) } \leq  t ^2 \right. ,
\\ & \hspace*{0.2cm} || \ [e^{itH} u^\omega_0]^3 \ ||_{ L^4( [-2\pi,2\pi] , \overline{H}^s (\R^3))  } \leq t^3, \underset{ N }{ \cap } \ || \Delta_N(e^{itH} u^\omega_0 ) ||_{L^4( [-2\pi,2\pi] , L^\infty (\R^3)) } \leq t  N^{-1/6},  
\\ & \hspace*{5.2cm} \left. \underset{ N }{ \cap } \ || \Delta_N( e^{itH} u^\omega_0 ) ||_{L^R( [-2\pi,2\pi] , \overline{W}^{s,4} (\R^3)) } \leq t  N^{s-1/4} \right)
\end{align*}
et l'objectif de cette partie est de démontrer le théorème suivant:
\begin{thm} \label{deviation}Il existe deux constantes $ C, c > 0$ telles que pour tout $ t \geq  0 $ et $ u_0 \in \overline{H}^{\sigma} (\R^3) $, 
\begin{equation} \label{inegalite gauss 4} 
P( \Omega_t ^ c ) \leq C e ^{ -  \frac{c t^2}{ ||u_0||^2_{ \overline{H}^{\sigma} (\R^3)  }   }  }.
\end{equation}
\end{thm}
On commence par établir des inégalités de type chaos de Wiener dans le cas où $ g_n \sim \mathcal{N}_\C(0,1) $ ou $ g_n \sim \mathcal{B}( \frac{1}{2}  ) $. Pour cela, on introduit deux définitions et on démontre 3 lemmes préliminaires.
\begin{dfn}
Pour $ p \in \N ^* $, on définit 
\begin{equation*}
\mathcal{A}_{2p} = \left\{  \sigma \in \mathcal{S}_{2p} \ / \ \sigma ^2 = Id \mbox{ et } \sigma (i) \neq i , \ \forall i \lbrace 1, .. , 2p \rbrace  \right\},
\end{equation*}
et pour $ p \in \N ^* $ et $ \sigma \in \mathcal{A}_{2p} $, on pose
\begin{equation*}
 \mathcal{I}(\sigma,p)= \left\{ i  \in \lbrace 1, ... , p \rbrace \ / \ \sigma(2i)=2i-1   \right\}   .
\end{equation*}
\end{dfn}
\begin{lem} \label{wienerbenoulli1} Soit $ X_n, \ n \in \N $ une suite de variables indépendantes telles que $ E(X_i^{2k+1}) = 0 $ pour tout $ k \in \N $, alors pour tout 2p-upplet $ (n_1, ... n_{2p} )  \in \N^{2p} $,
\\si $ E( X_{n_1} \times ... \times X_{n_{2p}}  ) \neq 0 $ alors il existe une permutation $ \sigma \in \mathcal{A}_{2p} $ telle que 
\\$ n_{\sigma(i) } = n_i , \ \forall i \in \lbrace 1, .. , 2p\rbrace $.
\end{lem}
\textit{Preuve.} Le lemme se démontre facilement par récurrence sur p. \hfill $ \boxtimes $
\begin{lem} \label{wienerbenoulli2}
Il existe une constante $ C> 0 $ telle que pour tout $ p \in \N ^* $,
\begin{equation*}
Card (  \mathcal{A}_{2p}  ) \leq (Cp)^p. 
\end{equation*}
\end{lem}
\textit{Preuve.} Il suffit d'utiliser la formule de Sterling. En effet, on a 
\vspace*{0.3cm}\\ \hspace*{1cm} \hfill $
Card (  \mathcal{A}_{2p}  ) = (2p-1)  \times  .... \times 3  \times 1 = \frac{(2p)!}{  2^p .  p!} \leq (Cp)^p.$ \hfill $ \boxtimes $
\begin{dfn}On définit
\begin{equation*}
l^2 = \left\{ c = (c_{n,m})_{n,m} \in \mathcal{F}( \N \times \N ) \ / \   || c||_{l^2} : =  \sqrt{\sum_{m,n} |c_{m,n}|^2 } < \infty \right\},
\end{equation*}
et
\begin{equation*}
\tilde{l}^1 = \left\{ c = (c_{n,m})_{n,m} \in \mathcal{F}( \N \times \N ) \ / \   || c||_{\tilde{l}^1} : = \sum_{n} |c_{n,n}| < \infty \right\}.
\end{equation*}
\end{dfn}
\begin{lem} \label{wienerbenoulli3}
Pour tout $ p \in \N^* $, $ \sigma \in \mathcal{A}_{2p} $ et $ c^1 , c^2 , ... c^p $ dans $ \tilde{l}^1\cap l^2  $,
\begin{align*}
\sum_{  \underset{n_{\sigma(i)} = n_i}{n_1,...,n_{2p},}  }   |c^1_{n_1,n_2}| \ ... \ |c^p_{n_{2p-1},n_{2p}}| & \leq  \prod^p_{  \underset{ i \in \mathcal{I}(\sigma,p ) }{i=1,}  } ||c^i||_{ \tilde{l}^1 } \ \times \ \prod^p_{  \underset{ i \notin \mathcal{I}(\sigma,p ) }{i=1,}  } ||c^i||_{ l^2 }
\\ & \leq \prod_{i=1}^p \bigg( \ ||c^i||_{ \tilde{l}^1 }  + || c^i ||_{l^2} \ \bigg).
\end{align*}
\end{lem}
\textit{Preuve.} La seconde inégalité est triviale à partir de la première. On démontre cette dernière par récurrence sur p. Les cas p=1, p=2 et p=3 sont clairs. Soit $ p \geq 4 $ et supposons le résultat établi pour $ q \in \lbrace  1, ... ,  p-1 \rbrace $.
\\\\ \textbf{-Cas $ \mathcal{I}(\sigma,p) \neq \varnothing  $.}
\\Si $ p \in \mathcal{I}(\sigma,p) $ alors il existe une permutation $ \sigma ' = \sigma|_{ \lbrace 1, ...,2(p-1) \rbrace } \in \mathcal{A}_{2(p-1)} $ telle que
\begin{align*}
\sum_{  \underset{n_{\sigma(i)} = n_i}{n_1,...,n_{2p},}  }   |c^1_{n_1,n_2}| \ ... \ |c^p_{n_{2p-1},n_{2p}}| & = \left( \sum_{  \underset{n_{\sigma'(i)} = n_i}{n_1,...,n_{2(p-1)},}  }   |c^1_{n_1,n_2}| \ ... \ |c^{p-1}_{n_{2p-3},n_{2p-2}}| \right) \times ||c^p||_{ \tilde{l}^1 }.
\end{align*}
et le résultat est prouvé par récurrence en remarquant que $ \mathcal{I}( \sigma',p-1) = \mathcal{I}( \sigma,p) \setminus \lbrace p \rbrace  $.
\\Sinon $ p \notin \mathcal{I}(\sigma,p) $ et il existe un entier $ i \in \lbrace 1, ..., p-1 \rbrace $ tel que $ i \in \mathcal{I}(\sigma,p)$.
\\Dans cette situation, posons
\begin{align*}
\begin{array}{rclcl}
\gamma : & \lbrace 1, ..., 2p  \rbrace \setminus  \lbrace 2i-1,2i \rbrace   & \rightarrow  & \lbrace 1, ..., 2p-2 \rbrace
\\ & k  & \mapsto & k & \mbox{ si } k \notin \lbrace 2p-1,2p \rbrace,
\\  & 2p-1 & \mapsto & 2i-1,
\\  & 2p & \mapsto & 2i,
\end{array}
\end{align*}
$ \tau = \sigma |_{ \lbrace 1, ..., 2p  \rbrace \setminus  \lbrace 2i-1,2i \rbrace }$ et $ \sigma' = \gamma  \circ \tau \circ \gamma^{-1}  $ pour obtenir que 
\begin{align*}
& \hspace*{1cm} \sum_{  \underset{n_{\sigma(i)} = n_i}{n_1,...,n_{2p},}  }   |c^1_{n_1,n_2}| \ ... \ |c^p_{n_{2p-1},n_{2p}}| 
\\ = \ & \bigg( \sum_{  \underset{n_{\sigma'(i)} = n_i}{n_1,...,n_{2(p-1)},}  }   |c^1_{n_1,n_2}| \ ...  \cancel{|c^{i}_{n,n}|} ... |c^{p-1}_{n_{2p-3},n_{2p-2}}| \ |c^{p}_{n_{2i-1},n_{2i}}| \bigg) \times ||c^i||_{ \tilde{l}^1 }.
\end{align*}
Remarquons que $ \sigma' \in \mathcal{A}_{2(p-1)}$ et que $ \mathcal{I}( \sigma',p-1) = \mathcal{I}( \sigma,p) \setminus \lbrace i \rbrace  $. Ainsi, nous pouvons appliquer l'hypothèse de récurrence à $ c^1 = c^1, ... , c^i = c^p, ..., c^{p-1} = c^{p-1} $ et le résultat suit.
\\\\\textbf{Cas $ \mathcal{I}(\sigma,p) = \varnothing  $ :}
\\Nous devons prouver que
\begin{equation*}
\sum_{  \underset{n_{\sigma(i)} = n_i}{n_1,...,n_{2p},}  }   |c^1_{n_1,n_2}| \ ... \ |c^p_{n_{2p-1},n_{2p}} | \leq \prod_{i=1}^p ||c^i||_{l^2}.
\end{equation*}
Ainsi, on remarque que le rôle des $ c^i $ est symétriques et qu'il est donc possible d'intervertir leurs positions. Par conséquent, il est possible de supposer que $ \sigma( 2p ) = 2p-2  $ et $ \sigma( 2p-1  ) =  2p-3 $ ou $ 2p-4$.  
\\\\\textbf{Sous-cas $ \sigma( 2p ) = 2p-2 $ et $ \sigma( 2p-1  ) =  2p-3 $ :}
\\Il existe une permuation $ \sigma ' = \sigma|_{ \lbrace 1, ...,2(p-2) \rbrace } \in \mathcal{A}_{2(p-2)} $ telle que
\begin{align*}
\sum_{  \underset{n_{\sigma(i)} = n_i}{n_1,...,n_{2p},}  }   |c^1_{n_1,n_2}| \ ... \ |c^p_{n_{2p-1},n_{2p}}| & = \left( \sum_{  \underset{n_{\sigma'(i)} = n_i}{n_1,...,n_{2(p-2)},}  }   |c^1_{n_1,n_2}| \ ... \ |c^{p-2}_{n_{2p-3},n_{2p-2}}| \right) \times \sum_{n,k}  c^{p-1}_{n,k} c^p_{n,k} .
\end{align*}
Puis, il suffit d'utiliser l'inégalité de Cauchy-Schwarz et l'hypothèse de récurrence pour prouver le résultat.
\\\\\textbf{Sous-cas $ \sigma( 2p ) = 2p-2 $ et $ \sigma( 2p-1  ) =  2p-4 $ :}
\\Posons $ \tau = \sigma|_{ \lbrace 1, .., 2p-5 \rbrace \cup \lbrace 2p-3  \rbrace } $,
\begin{align*}
\begin{array}{rclc}
\gamma : & \lbrace 1, ..., 2(p-2)  \rbrace  & \rightarrow  & \lbrace 1, .., 2p-5 \rbrace \cup \lbrace 2p-3  \rbrace
\\ & k  & \mapsto & k  \mbox{ si } k \in \lbrace 1, ..., 2p-5 \rbrace,
\\  & 2p-4 & \mapsto & 2p-3,
\end{array}
\end{align*}
et $ \sigma'= \gamma^{-1} \circ \tau \circ \gamma \in \mathcal{A}_{2(p-2)} $.
\\Comme $ 2p-4 = \sigma'(2p-4) \ \Leftrightarrow \ 2p-3 = \sigma (2p-3) $ et $ \sigma' = \sigma $ sur $ \lbrace 1, ... ,2p-5 \rbrace $, alors
\begin{align*}
\sum_{  \underset{n_{\sigma(i)} = n_i}{n_1,...,n_{2p},}  }   |c^1_{n_1,n_2}| \ ... \ |c^p_{n_{2p-1},n_{2p}}| = \sum_{  \underset{n_{\sigma'(i)} = n_i}{n_1,...,n_{2(p-2)},}  } \sum_{n,m}  |c^1_{n_1,n_2}| \ ... \ |c^{p-2}_{n_{2p-5},m}| \ |c^{p-1}_{n_{2p-4},n}| \ |c^p_{m,n}|.
\end{align*}
Puis, grâce à l'inégalité de Cauchy-Schwarz, on obtient  
\begin{equation*}
\sum_{n,m}   |c^{p-2}_{n_{2p-5},m}| \ |c^{p-1}_{n_{2p-4},n}| \ |c^p_{m,n}|  \leq \sqrt{ \sum_m |c^{p-2}_{n_{2p-5},m}|^2 . \sum_n |c^{p-1}_{n_{2p-4},n}|^2 }  \times  || c^p||_{l^2}  ,
\end{equation*}
pour ensuite appliquer l'hypothèse de récurrence à $ p-2 $, $ \sigma ' \in \mathcal{A}_{2(p-2)} $ et $ c^1 = c^1 , ..., c^{p-3}=c^{p-3}$ et
\begin{equation*}
  \tilde{c}^{p-2}_{k,l} = \displaystyle{ \sqrt{ \sum_m |c^{p-2}_{k,m}|^2 \times \sum_n |c^{p-1}_{l,n}|^2 } }.
\end{equation*}
Il est important de remarquer qu'il n'est pas clair que $ \mathcal{I}( \sigma', p-2 ) = \varnothing $ et qu'il est possible que $  p-2 \in \mathcal{I}( \sigma', p-2 ) $. On obtient alors
\begin{align*}
\sum_{  \underset{n_{\sigma(i)} = n_i}{n_1,...,n_{2p},}  }   |c^1_{n_1,n_2}| \ ... \ |c^p_{n_{2p-1},n_{2p}}|  \leq \prod_{i=1}^{p-3} || c^i ||_{l^2} \times  || \tilde{c}^{p-2} ||_{l_i^{?}} \times  || c^p ||_{l^2} ,
\end{align*}
où $ l^{?}_i $ désigne la norme $ l^2 $ ou $ \tilde{l}^1 $.
\\\\Finalement, pour conclure, on note que 
\begin{align*}
|| \tilde{c}^{p-2}||_{l^2}  = || \sqrt{ \sum_m |c^{p-2}_{k,m}|^2 \times  \sum_n |c^{p-1}_{l,n}|^2 } ||_{l^2} = || c^{p-2} ||_{l^2} \times || c^{p-1} ||_{l^2} ,
\end{align*}
et par l'inégalité de Cauchy Schwarz que
\begin{align*}
|| \tilde{c}^{p-2}||_{  \tilde{l}^1 }  = || \sqrt{ \sum_m |c^{p-2}_{k,m}|^2 . \sum_n |c^{p-1}_{l,n}|^2 } ||_{\tilde{l} ^1} \leq || c^{p-2} ||_{   l^2} \times || c^{p-1} ||_{l^2}.
\end{align*}
Ce qui achève la récurrence. \hfill $ \boxtimes $
\\\\Maintenant, grâce aux lemmes \ref{wienerbenoulli1}, \ref{wienerbenoulli2} et \ref{wienerbenoulli3}, on peut démontrer les estimées de chaos de Wiener qui nous serviront à estimer $P( \Omega_t ^ c ) $.
\begin{prop} Supposons que $ g_n \sim \mathcal{N}_\C(0,1) $ ou $ g_n \sim \mathcal{B}( \frac{1}{2} ) $ alors il existe une constante $ C> 0 $ telle que pour tout $ q \geq 2 $ et $ (c_n)_{n} \in l^2(\N), \ (c_{n,m})_{n,m} \in l^2(\N \times \N) $ et $ (c_{n,m,k})_{n,m,k} \in l^2(\N \times \N \times \N ) $, 
\begin{align}  
& \label{zygmound} \bigg| \bigg|  \sum_{ n \in \N   }  c_n \times g_n(\omega) \bigg| \bigg|_{L^q(\Omega)} \leq C \times  q^{ \frac{1}{2} } \times  \sqrt{ \sum_{ n \in \N } |c_n |^2  },
\\ & \label{zygmound'}  \bigg| \bigg|  \ \sum_{ n , m \in \N  }  c_{n,m} \times  g_n(\omega) \times  g_m(\omega)  \ \bigg| \bigg| _{L^q(\Omega)} \leq C  \times  q  \times \left( \sqrt{ \sum_{n,m \in \N} |c_{n,m} |^2 } + L \sum_{n  \in \N } |c_{n, n} | \right) ,
\\ & \label{zygmound''} \bigg| \bigg| \sum_{ n,m,k \in \N  }  c_{n,m,k } \times g_n(\omega)  \times  g_m(\omega)  \times g_k(\omega)  \ \bigg| \bigg| _{L^q(\Omega)} \leq C  \times  q^{ \frac{3}{2}  }  \times \left( \sqrt{ \sum_{n,m,k \in \N} |c_{n,m,k} |^2 } +\right.
\end{align}
\begin{equation*}
 \hspace*{1.5cm} + \left. L \times \left( \sum_{n \in \N} \sqrt{\sum_{m \in \N} |c_{n,n,m}|^2  } + \sum_{n \in \N} \sqrt{\sum_{m \in \N} |c_{n,m,n}|^2  } + \sum_{n \in \N} \sqrt{\sum_{m \in \N} |c_{m,n,n}|^2  } \right) \right) ,
\end{equation*}
où \begin{equation*} 
L = \left\{ \begin{array}{ll}
0 \mbox{ dans le cas Gaussien complexe,}
\\ 1 \mbox{ dans le cas Bernoulli (ou Gaussien réel).}
\end{array}
\right.
\end{equation*}
\end{prop}
\textit{Preuve.} \\ 1- \ Dans \cite{burq4}, il est montré que si
\begin{equation} \label{zygmound1} 
\exists \delta > 0  / \forall \alpha \in \R \mbox{ et } n \in \N , \ E(e^{ \alpha g_n    }) \leq e^{ \delta \alpha ^2 },
\end{equation}
alors (\ref{zygmound}) est satisfait. Ainsi, (\ref{zygmound}) est vérifiée puisque (\ref{zygmound1}) est satisfait si $ g_n \sim \mathcal{N}_\R(0,1) $ ou $ g_n \sim \mathcal{B}( \frac{1}{2} ) $.
\\\\ 2- Pour le cas Gaussien, d'après \cite{thomann1}, proposition 2.4 (Wiener chaos estimates), il existe une constante $ C > 0  $ telle que pour tout $ q \geq 2 $,
\begin{equation} \label{sauvetagedevie}
\bigg| \bigg|  \ \sum_{ n,m \in \N   }  c_{n,m} \times  g_n(\omega) \times  g_m(\omega)  \ \bigg| \bigg|_{L^q(\Omega)} \leq C  \times  q  \times \bigg| \bigg|  \ \sum_{ n,m \in \N   }  c_{n,m} \times  g_n(\omega) \times  g_m(\omega)  \ \bigg| \bigg|_{L^2(\Omega)}.
\end{equation}
Or
\begin{align*}
& \bigg| \bigg|  \  \sum_{ n,m \in \N   }  c_{n,m} \times  g_n(\omega) \times  g_m(\omega)  \ \bigg| \bigg| ^2 _{L^2(\Omega)}
\\  = \hspace*{1cm}   & \sum_{ n, n' , m , m' \in \N   }  c_{n,m} \times \overline{c_{n',m'} } \times  E \left( g_n(\omega) \times  g_m(\omega) \times \overline{ g_{n'}(\omega) \times  g_{m'}(\omega) } \right)
\\  \leq  \hspace*{1cm} & \sum_{ n = n' =  m = m' \in \N   } | \ | + \sum_{ n = n' , m = m' \in \N   } | \ | + \sum_{ n = m  ,   n' = m' \in \N   } | \ | + \sum_{ n = m'  ,   n' = m \in \N   }| \ | .
\end{align*}
Mais
\begin{align*}
 & \sum_{ n=n'=m=m'  \in \N   } \bigg| c_{n,m} \times \overline{c_{n',m'} } \times  E \left( g_n(\omega) \times  g_m(\omega) \times \overline{ g_{n'}(\omega) \times  g_{m'}(\omega) } \right) \bigg| 
\\ = & \hspace*{0.68cm} \sum_{ n \in \N   } |c_{n,n}|^2 E ( |g_n(\omega)|^4 )
\\ \leq & \hspace*{0.6cm}  \sum_{n,m \in \N} |c_{n,m} |^2 ,
\end{align*}
et en utilisant que $ E( g_n(\omega)^2  ) = 0 $, on trouve
\begin{align*}
& \sum_{ n = m , n'= m' \in \N   } \bigg| c_{n,m} \times \overline{c_{n',m'} } \times  E \left( g_n(\omega) \times  g_m(\omega) \times \overline{ g_{n'}(\omega) \times g_{m'}(\omega) } \right) \bigg| 
\\ = & \hspace*{0.5cm} \sum_{ n,m \in \N   } | c_{n,n} | \times | c_{m,m}| \times \bigg|  E \left( g_n(\omega)^2 \times  \overline{ g_m(\omega) ^2 } \right)\bigg|
\\ = & \hspace*{0.5cm} \sum_{ n \in \N   } |c_{n,n}|^2 \times  E \left( |g_n(\omega)|^4 \right)
\\  \leq & \hspace*{0.45cm} \sum_{n,m \in \N} |c_{n,m} |^2 .
\end{align*}
Ainsi (\ref{zygmound'}) est démontrée dans le cas $ g_n \sim \mathcal{N}_{\C}(0,1) $. Nous pouvons procéder de la même manière pour obtenir (\ref{zygmound''}) puisque l'inégalité (\ref{sauvetagedevie}) est vraie pour un produit quelconque de variables aléatoires. 
\\\\ 3- Dans le cas Bernoulli, démontrons (\ref{zygmound'}). Nous pouvons limiter la preuve au cas où $ q = 2 p $ avec $ p \in  \N^* $. Il suffit alors de montrer que
\begin{equation*}
\bigg| \bigg|  \ \sum_{ n,m  \in \N   }  c_{n,m} \times  g_n(\omega) \times  g_m(\omega)  \ \bigg| \bigg| ^{2p}_{L^{2p}(\Omega)} \leq (Cp)^{2p}   \times \left( \sqrt{ \sum_{n,m \in \N} |c_{n,m} |^2 } +   \sum_{n \in \N } |c_{n,n} |  \right)^{2p}.
\end{equation*}
Nous pouvons utiliser successivement les lemmes \ref{wienerbenoulli1}, \ref{wienerbenoulli3} et \ref{wienerbenoulli2} pour obtenir 
\begin{align*}
& \bigg| \bigg|  \ \sum_{ n,m \in \N   }  c_{n,m} \times  g_n(\omega) \times  g_m(\omega)  \ \bigg| \bigg| ^{2p}_{L^{2p}(\Omega)} 
\\ = \ &   \sum    \limits_{  { n_1,...,n_{4p}  } } c_{n_1,n_2} ... c_{n_{2p-1},n_{2p}} . \overline{ c_{n_{2p+1},n_{2p+2}} }...\overline{c_{n_{4p-1},n_{4p}}} \times   E \left( \prod_{i=1}^{4p}  g_{n_i}   \right)
\\  \leq  \ &  \sum    \limits_{  { n_1,... ,n_{4p}  } }   | c_{n_1,n_2}| ...  |c_{n_{4p-1},n_{4p}}| \times \bigg| E \left( \prod_{i=1}^{4p}  g_{n_i}   \right) \bigg| 
\\ \leq   \  & \sum_{\sigma \in \mathcal{A}_{4p} }  \sum    \limits_{ \underset{n_{\sigma(i)} = n_i }{ n_1,...,n_{4p},  } }   | c_{n_1,n_2}| ...  |c_{n_{4p-1},n_{4p}}|
\\  \leq \ & Card (  \mathcal{A}_{4p} ) \times   \left( \sqrt{ \sum_{n,m \in \N} |c_{n,m} |^2 } +   \sum_{n \in \N } |c_{n,n} |  \right)^{2p}
\\ \leq \ &   (  Cp )^{2p} \times   \left( \sqrt{ \sum_{n,m \in \N} |c_{n,m} |^2 } +   \sum_{n \in \N } |c_{n,n} |  \right)^{2p}.
\end{align*}
Ce qui démontre (\ref{zygmound'}). Pour obtenir (\ref{zygmound''}), on peut remarquer que le lemme \ref{wienerbenoulli3} reste vrai pour un produit de 3 variables aléatoires (il suffit de faire la preuve avec des suites de 3 variables). \hfill $ \boxtimes $
\\\\\textit{Preuve du théorème \ref{deviation}.}
\begin{align*}
P( \Omega_t^c  ) & \leq P \left( \omega \in \Omega /  ||u^\omega_0||_{ \overline{H}^\sigma(\R^3) } \geq t \right)
\\ & + P \left( \omega \in \Omega / || \ [e^{itH} u^\omega_0]^2 \ ||_{  L^4( [-2\pi,2\pi] , \overline{H}^s (\R^3)) } \geq  t ^2 \right)
\\ &  +  P \left( \omega \in \Omega / || \ [e^{itH} u^\omega_0]^3 \ ||_{ L^4( [-2\pi,2\pi] , \overline{H}^s (\R^3))  } \geq t ^3 \right)
\\ & +  P \left( \underset{ N }{ \bigcup } \left\{ \omega \in \Omega / || \Delta_N(e^{itH} u^\omega_0 ) ||_{L^4( [-2\pi,2\pi] , L^\infty (\R^3)) } \geq t  N^{-1/6} \right\} \right)
\\ & +  P \left( \underset{ N }{ \bigcup } \left\{ \omega \in \Omega /  || \Delta_N( e^{itH} u^\omega_0 ) ||_{L^R( [-2\pi,2\pi] , \overline{W}^{s,4} (\R^3)) } \geq t  N^{s-1/4}  \right\} \right)
\end{align*}
ainsi il suffit d'établir la majoration (\ref{inegalite gauss 4}) pour chacun des termes. Effectuons la démonstration pour le second et le quatrième terme (pour les autres termes, la démarche est identique).
\\\\ \textbf{I/ Cas $ || \ [e^{itH} u^\omega_0]^2 \ ||_{  L^4( [-2\pi,2\pi] , 
 \overline{H}^s (\R^3)) } \geq  t ^2 $ : }
\\\\ Comme $ s \in ]  \frac{1}{2} , \frac{1}{2} + \sigma [$, il suffit d'obtenir le lemme suivant:
\begin{lem}  \label{lemme d'estimation du terme bilinéaire} Pour tout $ \epsilon > 0 $, il existe 2 constantes $ C,c > 0 $ telles que pour tout $ t > 0 $ et $ u_0 \in  \overline{H}^{\sigma} (\R^3) $,
\begin{equation*}
P \left( \omega \in \Omega /    || \ [e^{itH} u^\omega_0]^2 \ ||_{  L^4( [-2\pi,2\pi] , 
 \overline{H}^{\sigma + 1 /2 -2\epsilon} (\R^3)) } \geq  t ^2  \right)  \leq C e ^{ -  \frac{c t^2}{ ||u_0||^2_{ \overline{H}^{\sigma} (\R^3)  }   }  } .
\end{equation*}
\end{lem}
\textit{Preuve du lemme \ref{lemme d'estimation du terme bilinéaire}.} Remarquons qu'il suffit d'obtenir l'estimation pour $ t \geq C ||u_0||_{ \overline{H}^{\sigma}(\R^3)  } $. Grâce aux inégalités de Markov et Minkowski, on obtient pour $ q \geq 4 $,
\begin{align*}
 & P \left( \omega \in \Omega /    || \ [e^{itH} u^\omega_0]^2 \ ||_{  L^4( [-2\pi,2\pi] , 
 \overline{H}^{\sigma + 1 /2 -2\epsilon} (\R^3)) } \geq  t ^2  \right) 
\\  \leq \ & P \left( \omega \in \Omega / || \  \ H^{  \sigma/2+1/4-\epsilon }  \   [e^{itH} u^\omega_0]^2 \ ||^q_{  L^4( [-2\pi,2\pi] , L^2 (\R^3)) } \geq  t ^{2q} \right)
\\  \leq \ & t^{-2q} \times E_\omega  \left( || \  \ H^{  \sigma/2+1/4-\epsilon }  \  [e^{itH} u^\omega_0]^2 \ || ^q _{  L^4( [-2\pi,2\pi] , L^2 (\R^3)) } \right)
\\ \leq \ & t^{-2q} \times  || \  \ H^{  \sigma/2+1/4-\epsilon }  \  [e^{itH} u^\omega_0]^2 \ || ^q _{ L ^q( \Omega ,   L^4( [-2\pi,2\pi] , L^2 (\R^3))) } 
\\  \leq \ & t^{-2q} \times  || \  \ H^{  \sigma/2+1/4-\epsilon }  \  [e^{itH} u^\omega_0]^2 \ || ^q _{   L^4( [-2\pi,2\pi] , L^2 (\R^3, L ^q( \Omega) )) } .
\end{align*}
Puis, par (\ref{zygmound'}), on établit
\begin{align*}
 & || \  H^{  \sigma/2+1/4-\epsilon }  \  [e^{itH} u^\omega_0]^2 \ || _{   L ^q( \Omega) } 
 \\  \leq \  &  \bigg| \bigg|   \ H^{  \sigma/2+1/4-\epsilon }  \  \left[  \sum_{n,m \in \N} e^{it(\lambda_n^2+\lambda_m^2)} c_n c_m  h_n(x) h_m(x)   g_n(\omega )   g_m(\omega) \ \right] \ \bigg| \bigg| _{   L ^q( \Omega) } 
\\ \leq \ &  \bigg| \bigg|  \  \sum_{n,m \in \N} e^{it(\lambda_n^2+\lambda_m^2)} c_n c_m \times  H^{  \sigma/2+1/4-\epsilon }  \  [  h_n(x) h_m(x)  ]  \times g_n(\omega )   g_m(\omega)  \ \bigg| \bigg|  _{   L ^q( \Omega) } 
\\ \leq \ & C \times q \times \bigg(  \sqrt{   \sum_{n,m \in \N } |c_n|^2   |c_m|^2 |  \ H^{  \sigma/2+1/4-\epsilon }  \  [  h_n(x) h_m(x)  ] \   |^2     } 
\\ & \hspace*{5cm} + L \sum_{ n \in \N  } |c_n|^2  |  H^{  \sigma/2+1/4-\epsilon }  \  [  h_n^2(x)   ]   |  \bigg) .
\end{align*}
Prenons $ L=0$ pour simplifier les calculs (le terme avec L=1 s'estime de la même façon). Par l'inégalité triangulaire, on trouve
\begin{align*}
& P \left( \omega \in \Omega / || \ [e^{itH} u^\omega_0]^2 \ ||_{  L^4( [-2\pi,2\pi] , \overline{H}^{\sigma+1/2-2\epsilon} (\R^3)) } \geq  t ^2 \right)
\\  \leq \ & \left( \frac{C q }{ t^2 }  \right) ^q \times  \bigg| \bigg| \sum_{ n,m \in \N } |c_n|^2 \times |c_m|^2 \times |  \ H^{  \sigma/2+1/4-\epsilon }  \  [  h_n(x) h_m(x)  ] \   |^2  \bigg| \bigg| _{L^2( [-2\pi,2\pi] , L^1 (\R^3 ) )     }  ^{q/2}
\\  \leq \ & \left(  \frac{C q }{ t^2 } \right) ^q \times \left( \sum_{ n , m \in  \N} |c_n|^2 \times |c_m|^2 \times || \  H^{  \sigma/2+1/4-\epsilon }  \  [  h_n(x) h_m(x)  ] \     ||^2_{L^4( [-2\pi,2\pi] , L^2 (\R^3 ) )     } \right) ^{q/2}.
\end{align*}
Puis grâce à la proposition \ref{propre6} avec $ \delta = \epsilon  $, on arrive à
\begin{align*} 
|| \  H^{  \sigma/2+1/4-\epsilon }  \  [  h_n(x) h_m(x)  ] \     ||^2_{L^4( [-2\pi,2\pi] , L^2 (\R^3 ) )     }  & \leq C_\epsilon \times \max( \lambda_n , \lambda_m )^{2 ( \sigma - \epsilon ) } 
\\ & \leq C \times \max( \lambda_n , \lambda_m )^{2  \sigma  }.
\end{align*}
Et finalement, on a 
\begin{align*}
& P \left( \omega \in \Omega  / || \ [e^{itH} u^\omega_0]^2 \ ||_{  L^4( [-2\pi,2\pi] , \overline{H}^{\sigma+1/2-2\epsilon} (\R^3)) } \geq  t ^2 \right)
\\  \leq \ & \left(  \frac{C q }{ t^2 } \right) ^q \times \left( \sum_{n,m \in \N} |c_n|^2 \times |c_m|^2 \times \max( \lambda_n , \lambda_m )^{2  \sigma  }     \right)  ^{q/2}
\\ \leq \ & \left(  \frac{C q   ||u_0||^2_{ \overline{H}^{\sigma}  (\R^3)} }{ t^2 }  \right) ^q.
\end{align*}
Il suffit ensuite de choisir $ q = \frac{t^2}{2C||u_0||^2_{ \overline{H}^{\sigma } (\R^3)  }} \geq 4 $ puisque $ t \gg 1 $ pour conclure . \hfill $ \boxtimes $
\\\\\textbf{II/ Cas $ || \Delta_N[ e^{itH} u^\omega_0 ] ||_{  L^4 ( [-2\pi,2\pi] , 
L^\infty (\R^3)) } \geq  N^{-1/6 } t  $ : }Commençons par établir le lemme suivant:
\begin{lem} \label{lemme d'estimation frequence tronque} Pour tout $ p_1,p_2 \in [ 2,  \infty [  $ et $ \epsilon > 0 $, il existe deux constantes $ C,c > 0 $ telles que pour tout $ t > 0 $, $ N \geq 1 $ et $ u_0 \in \overline{H}^{\sigma-\epsilon}(\R^3)$,
\begin{equation*}
 P \left( \omega \in \Omega / || \Delta_N(e^{itH} u^\omega_0 ) ||_{L^{p_1}( [-2\pi,2\pi] , W^{\epsilon,p_2} (\R^3)) } \geq t  N^{-1/6-\sigma+2 \epsilon} \right) \leq C e^{-  \frac{c t^2}{ ||  \Delta_N( u_0 ) ||^2_{ \overline{H}^{\sigma-\epsilon}(\R^3)  }   }    }.
\end{equation*}
\end{lem}
\textit{Preuve du lemme \ref{lemme d'estimation frequence tronque}} Quitte à remplacer $ u_0 $ par $ \Delta_N (u _ 0 ) $, on se ramène à démontrer que 
\begin{equation*}
 P \left( \omega \in \Omega / || \Delta'_N(e^{itH} u^\omega_0 ) ||_{L^{p_1}( [-2\pi,2\pi] , W^{\epsilon,p_2} (\R^3)) } \geq t  N^{-1/6-\sigma+2 \epsilon} \right) \leq C e^{-  \frac{c t^2}{ || u_0  ||^2_{ \overline{H}^{\sigma-\epsilon}(\R^3)  }   }    }.
\end{equation*}
Comme pour le lemme \ref{lemme d'estimation du terme bilinéaire}, il suffit de prouver l'estimation pour $ t \geq C ||u_0||_{ \overline{H}^{\sigma-\epsilon}(\R^3)  }    $.
\\\\Grâce aux inégalités de Markov et Minkowski, on obtient pour $ q \geq p_1,p_2 $,
\begin{align*}
&  P \left( \omega \in \Omega  / || \Delta'_N(e^{itH} u^\omega_0 ) ||_{L^{p_1}( [-2\pi,2\pi] , W^{\epsilon,p_2} (\R^3)) } \geq  t N^{-1/6-\sigma+2\epsilon}  \right)
\\ \leq \ &  \left( \frac{N^{1/6+\sigma-2\epsilon}  \times E_\omega (    || \Delta'_N(e^{itH} u^\omega_0 ) ||_{L^{p_1}( [-2\pi,2\pi] , W^{\epsilon,p_2} (\R^3)) }  ) }{t} \right) ^ q
\\  \leq  \ & \left( \frac{ N^{1/6+\sigma-2\epsilon} }{t} \right) ^q  \times    || H^{  \epsilon/2 }  \Delta'_N(e^{itH} u^\omega_0 ) ||_{ L^q(\Omega)  ,  L^{p_1}( [-2\pi,2\pi] ) , L^{p_2} (\R^3) }   ^q
\\ \leq \ & \left( \frac{ N^{1/6+\sigma-2\epsilon} }{t} \right) ^q  \times      || H^{  \epsilon/2 }  \Delta'_N(e^{itH} u^\omega_0 ) ||_{  L^{p_1}( [-2\pi,2\pi] ) , L^{p_2} (\R^3) , L^q(\Omega) }   ^q.
\end{align*}
Or, d'après (\ref{zygmound}),
\begin{align*}
|| H^{  \epsilon/2 }  \Delta'_N(e^{itH} u^\omega_0 ) ||_{   L^q(\Omega) }  & \leq \bigg| \bigg|  \sum_{ n \in \N } \phi \left( \frac{\lambda_n^2}{N^2} \right) \lambda_n^{\epsilon}  e^{it\lambda_n^2} c_n  h_n(x) g_n( \omega ) \bigg| \bigg|_{   L^q(\Omega) }  
\\ & \leq \bigg| \bigg|  \sum_{ \lambda_n \sim  N } \phi  \left( \frac{\lambda_n^2}{N^2}  \right) \lambda_n^{\epsilon}  e^{it\lambda_n^2} c_n h_n(x) g_n ( \omega ) \bigg| \bigg| _{   L^q(\Omega) }  
\\ & \leq C \times \sqrt{q} \times  \sqrt{ \sum_{ \lambda_n \sim  N } \phi^2 \left( \frac{\lambda_n^2}{N^2} \right) \lambda_n^{2 \epsilon}  \times | c_n |^2 \times | h_n(x) |^2 }.
\end{align*}
Alors, par (\ref{propre1}) et (\ref{propre2}), on trouve
\begin{align*}
& || H^{  \epsilon/2 }   \Delta'_N(e^{itH} u^\omega_0 ) ||_{  L^{p_1}( [-2\pi,2\pi] ) , L^{p_2} (\R^3) , L^q(\Omega) }  
 \\  \leq \ & C \times \sqrt{q} \times \bigg| \bigg| \sum_{ \lambda_n \sim  N } \phi^2 \left( \frac{\lambda_n^2}{N^2} \right) \times \lambda_n^{2 \epsilon}  \times | c_n |^2 \times | h_n(x) |^2 \bigg| \bigg| ^{1/2}_{ L^{p_1/2}( [-2\pi,2\pi] ) , L^{p_2/2} (\R^3) }
\\  \leq \ & C \times \sqrt{q} \times \sqrt{  \sum_{ \lambda_n \sim  N } \phi^2 \left( \frac{\lambda_n^2}{N^2} \right)  \times \lambda_n^{2 \epsilon} \times | c_n |^2 \times    || h_n(x)  ||^2_{ L^{p_2} (\R^3) } }
\\ \leq \ & C \times \sqrt{q} \times \sqrt{  \sum_{ \lambda_n \sim  N } \phi^2 \left( \frac{\lambda_n^2}{N^2} \right) \times \lambda_n^{2 \epsilon-1/3}  \times | c_n |^2   }
\\  \leq \ & C \times \sqrt{q} \times N^{-\sigma - 1/6 + 2\epsilon}  \times \sqrt{  \sum_{ \lambda_n \sim  N } \phi^2 \left( \frac{\lambda_n^2}{N^2} \right) \times \lambda_n^{2 (\sigma - \epsilon )  }  \times | c_n |^2   }
\\ \leq \ & C \times \sqrt{q} \times N^{-\sigma - 1/6 + 2\epsilon}  \times ||  \Delta'_N (  u_0 ) ||_{ \overline{H}^{\sigma-\epsilon}(\R^3) }
\\ \leq \ & C \times \sqrt{q} \times N^{-\sigma - 1/6 + 2\epsilon}  \times ||   u_0  ||_{ \overline{H}^{\sigma-\epsilon}(\R^3) }.
\end{align*}
Finalement, on a pour tout $ q \geq p_1,p_2 $,
\begin{align*}
 P & \left( \omega \in \Omega / || \Delta'_N(e^{itH} u^\omega_0 ) ||_{L^{p_1}( [-2\pi,2\pi] , W^{\epsilon,p_2} (\R^3)) } \geq  t N^{-1/6-\sigma+2\epsilon}  \right) 
\\ & \hspace*{5cm} \leq  \left( \frac{C \times \sqrt{q} \times ||  u_0 ||_{ \overline{H}^{\sigma-\epsilon}(\R^3) } }{t}  \right) ^q.
\end{align*}
Il suffit alors de choisir $ q =  \left( \frac{t }{4C ||u_0||_{ \overline{H}^\sigma(\R^3) }   } \right) ^2 \geq p_1,p_2 $ pour prouver le lemme \ref{lemme d'estimation frequence tronque}. \hfill $ \boxtimes $
\\\\Ensuite, pour $ p_2 = \frac{3}{\epsilon} + \epsilon  $ , on a 
\begin{equation*}
W^{ \epsilon, p_2 } (\R^3) \hookrightarrow L^{\infty } ( \R^3).
\end{equation*}
Ainsi pour tout $ p_1 \in [2,\infty[$ et $ \epsilon > 0  $, il existe deux constantes $ C,c $ telles que pour tout $ t > 0 $, $ N\geq 1 $ et $ u_0 \in \overline{H}^{\sigma-\epsilon}(\R^3)   $
\begin{equation*}
 P \left( \omega \in \Omega / || \Delta_N(e^{itH} u^\omega_0 ) ||_{L^{p_1}( [-2\pi,2\pi] , L^{\infty} (\R^3)) } \geq t  N^{-1/6-\sigma+2\epsilon} \right) \leq C e^{-  \frac{c t^2}{ ||  \Delta_N ( u_0 ) ||^2_{ \overline{H}^{\sigma-\epsilon}(\R^3)  }   }    }
\end{equation*}
Ensuite, nous pouvons choisir $ p_1 = 4 $ et utiliser que $ || \Delta_N( u_0 ) ||^2_{  \overline{H}^{\sigma-\epsilon} ( \R^3 ) }  \leq N^{-2\epsilon}  ||u_0||^2_{  \overline{H}^{\sigma} ( \R^3 ) } $ pour obtenir pour tout $ \epsilon > 0 $ l'existence de deux constantes $ C,c > 0 $ telles que pour tout $ t > 0 $, $ N\geq 1 $ et $ u_0 \in \overline{H}^{\sigma}(\R^3)   $,
\begin{align} \label{inegalite gauss 1}
 P \left( \omega \in \Omega / || \Delta_N(e^{itH} u^\omega_0 ) ||_{L^4( [-2\pi,2\pi] , L^{\infty} (\R^3)) } \geq t  N^{-1/6} \right) \leq C e^{-  \frac{c N^{2(\sigma-\epsilon)} t^2}{ || u_0  ||^2_{ \overline{H}^{\sigma}(\R^3)  }   }    }.
\end{align}
Nous devons prouver qu'il existe deux constantes $ C,c > 0  $ telles que pour tout $ t > 0 $, $ N\geq 1 $ et $ u_0 \in \overline{H}^{\sigma}(\R^3)   $,
\begin{equation} \label{inegalite gauss 2}
P  \left( \underset{ N }{ \bigcup }  \left\{ \omega \in \Omega / || \Delta_N(e^{itH} u^\omega_0 ) ||_{L^4( [-2\pi,2\pi] , L^\infty (\R^3)) } \geq t  N^{-1/6} \right\} \right) \leq C e^{-  \frac{c t^2}{ ||  u_0  ||^2_{ \overline{H}^{\sigma}(\R^3)  }   }    }.
\end{equation}
Depuis
\begin{equation*}
P  \left( \underset{ N }{ \bigcup } \left\{ \omega \in \Omega / || \Delta_N(e^{itH} u^\omega_0 ) ||_{L^4( [-2\pi,2\pi] , L^\infty (\R^3)) } \geq t  N^{-1/6} \right\} \right) \leq 1,
\end{equation*}
il suffit de montrer (\ref{inegalite gauss 2}) pour $ t \geq C ||u_0||_{  \overline{H}^{\sigma}(\R^3) } $.
\\\\On pose $ \alpha =  \left( \frac{t}{C||u_0||_{  \overline{H}^{\sigma}(\R^3) } }  \right) ^2 $ et en choisissant $ \epsilon < \sigma $ dans (\ref{inegalite gauss 1}), il suffit de montrer que 
\begin{equation} \label{inegalite gauss 3}
\forall \delta >0, \ \exists \ C,c >0 / \ \forall \alpha \geq 1 , \ \sum_N e^{-\alpha N^\delta } \leq C e^{-c\alpha }.
\end{equation}
Or 
\begin{align*}
\sum_N e^{-\alpha N^\delta } & = \sum_{ k \geq 1}e^{-\alpha (2^{k-1})^\delta }  \leq  \int_0^\infty e^{-\alpha (2^{x-1})^\delta } dx  \leq \int _{1/2} ^\infty   \frac{ e^{-\alpha y^\delta } }{y} dy    \\ & \leq   \int _{(1/2)^\delta} ^\infty   \frac{ e^{-\alpha z } }{z} \frac{dz}{\delta}   \leq \int _{\alpha (1/2)^\delta} ^\infty   \frac{ e^{-t } }{t} \frac{dt}{\delta}  \leq  \frac{2^\delta}{\delta\alpha} e^{  - \frac{\alpha}{2^\delta} } \leq  C e^{ - c \alpha  },
\end{align*}
et (\ref{inegalite gauss 2}) est prouvé ainsi que l'estimation du terme dans le cas \textbf{II}. 
\section{Preuves des théorèmes}
\subsection{Preuve du théorème \ref{thm1}}
Pour montrer le théorème \ref{thm1}, il suffit de montrer que pour tout $ t> 0 , \  P( \Omega_t ) > 0 $. Pour cela, on commence par établir qu'il suffit de montrer le résultat pour un nombre fini de termes dans la donnée initiale. On commence par introduire la définition suivante:
\begin{dfn}
Pour $ u_0 =   \displaystyle{  \sum_{n  \in \N} c_n h_n(x) }  $ une fonction de $ L^2 ( \R^3 ) $, on définit pour $ K \in \N^*$,
\begin{align*}
& [u_0]_K = \sum_{\lambda_n < K } c_n h_n(x),
\\ & [u_0]^K = \sum_{\lambda_n  \geq K} c_n h_n(x).
\end{align*}
\end{dfn}
Pour ensuite énoncer le théorème.
\begin{thm} \label{deviation2}
Pour tout $ t > 0 $ et $  \alpha \in ]0,1]  $, il existe un entier $ K \in \N^* $ tel que    
\begin{align*}
 \hspace*{1cm} P( \Omega_t ) \geq  (1 - \alpha ) \times 
\\  \mu \bigg(  u_0 \in \overline{H}^\sigma(\R^3 ) /   &  ||\ [u_0]_K||_{ \overline{H}^\sigma(\R^3) } \leq \frac{t}{2} , \  || \ (e^{itH} [u_0]_K )^2 \ ||_{  L^4( [-2\pi,2\pi] , \overline{H}^s (\R^3)) } \leq  \frac{t^2}{2} ,
\\  & || \ (e^{itH} [u_0]_K)^3 \ ||_{ L^4( [-2\pi,2\pi] , \overline{H}^s (\R^3))  } \leq \frac{t^3}{2} , 
\\  & \hspace*{0.5cm} \underset{ N }{ \bigcap } \left\{   || \Delta_N(e^{itH} [u_0]_K ) ||_{L^4( [-2\pi,2\pi] , L^\infty (\R^3)) } \leq   \frac{t N^{-1/6}}{2} \right\} ,
 \\ & \hspace*{1cm} \left. \underset{ N }{ \bigcap } \left\{   || \Delta_N( e^{itH} [u_0]_K ) ||_{L^R( [-2\pi,2\pi] , \overline{W}^{s,4} (\R^3)) } \leq  \frac{t  N^{s-1/4}}{2} \right\} \right) .
\end{align*}
\end{thm}
\textit{Preuve.} Par indépendance, on trouve
\begin{align*}
 P( \Omega_t ) \geq  & \mu \bigg(  || \  [ u_0]^K ||_{ \overline{H}^\sigma(\R^3) } \leq \frac{t}{2}  , \   || \ (e^{itH} [u_0]^K ) ^2  ||_{  L^4( [-2\pi,2\pi] , \overline{H}^s (\R^3)) } \leq  \frac{t^2}{2}  , 
\\ & \hspace*{1cm}  || \ (e^{itH} [u_0]^K)^3 \ ||_{ L^4( [-2\pi,2\pi] , \overline{H}^s (\R^3))  }  \leq \frac{t^3}{2}, \
\\  & \hspace*{2cm} \underset{ N }{ \bigcap }  \left\{  || \Delta_N(e^{itH} [u_0]^K ) ||_{L^4( [-2\pi,2\pi] , L^\infty (\R^3)) } \leq \frac{t  N^{-1/6} }{ 2 } \right\} , 
\\ & \hspace*{3cm} \underset{ N }{ \bigcap } \left\{  || \Delta_N( e^{itH} [u_0]^K ) ||_{L^R( [-2\pi,2\pi] , \overline{W}^{s,4} (\R^3)) } \leq \frac{t  N^{s-1/4} }{2 } \right\} \bigg) 
\\  \times & \mu \bigg(  || \ [ u_0]_K ||_{ \overline{H}^\sigma(\R^3) } \leq \frac{t}{2}, \  || \ (e^{itH} [u_0]_K )^2 \ ||_{  L^4( [-2\pi,2\pi] , \overline{H}^s (\R^3)) } \leq  \frac{t^2}{2},
\\ & \hspace*{1cm}  || \ (e^{itH} [u_0]_K)^3 \ ||_{ L^4( [-2\pi,2\pi] , \overline{H}^s (\R^3))  } \leq \frac{t^3}{2} , 
\\ & \hspace*{2cm} \underset{ N }{ \bigcap }  \left\{  || \Delta_N(e^{itH} [u_0]_K ) ||_{L^4( [-2\pi,2\pi] , L^\infty (\R^3)) } \leq \frac{t  N^{-1/6} }{2} \right\},
\\ &  \hspace*{3cm} \underset{ N }{ \bigcap } \left\{  || \Delta_N( e^{itH} [u_0]_K ) ||_{L^R( [-2\pi,2\pi] , \overline{W}^{s,4} (\R^3)) } \leq \frac{t  N^{s-1/4} }{2}  \right\}  \bigg) .
\end{align*}
Notons par $ P_{t,K} $ le premier terme probabiliste de cette inégalité. Alors par le théorème \ref{deviation}, pour tout $ t\geq 0 $ et $ K \in \N^* $,  
\begin{equation*}
P_{t,K} \geq 1 - C e ^{  - \frac{c t ^2 } {|| \ [u_0]^K \ ||^2_{ \overline{H}^\sigma ( \R^3 ) }} },
\end{equation*}
avec
\begin{equation*}
 \underset{K \rightarrow \infty }{\lim} || \ [u_0]^K \ ||^2_{ \overline{H}^\sigma ( \R^3 ) } = 0.
\end{equation*}
Par conséquent, il existe bien un entier non nul $ K $ tel que $ C e ^{  - \frac{c t ^2 } {|| \ [u_0]^K \ ||^2_{ \overline{H}^\sigma ( \R^3 ) }} } \leq \alpha. $ \hfill $\boxtimes$ 
\\\\\textit{Remarque :} Dans ce théorème, s'il existe $ \sigma ' > \sigma $ avec $ u_0 \in \overline{H}^{\sigma'}(\R^3) $ alors nous pouvons choisir 
\\$ M \geq ||u_0||_{ \overline{H}^{\sigma'}(\R^3)} \times \sqrt{ \frac{ \log  \left( \frac{C}{\alpha} \right) }{c} }   $ pour avoir $ P_{MK^{-(\sigma'-\sigma)},K} \geq  1 - \alpha. $
\\\\Puis, nous pouvons choisir $ K  \geq \left( \frac{M}{t}  \right)^{ \frac{1}{\sigma'-\sigma}  }  $ pour obtenir
\begin{equation*}
 P_{t,K}  \geq P_{MK^{-(\sigma'-\sigma)},K} \geq  1 - \alpha.
\end{equation*}
Et finalement $ K = \left( \frac{  ||u_0||_{ \overline{H}^{\sigma'}(\R^3)} } {t} \times \sqrt{  \frac{ \log  \left( \frac{C}{\alpha} \right) }{c} }   \right)^{ \frac{1}{\sigma'-\sigma}  } $ satisfait le théorème \ref{deviation2}.
\\\\Enfin, on démontre le résultat pour un nombre fini de termes dans la donnée initiale.
\begin{prop} \label{esti6} Il existe une constante $ C > 0 $ telle que pour tout $ K\in \N $,
\begin{align} \label{esti1}
\left\{  \omega \in \Omega / \sum_{\lambda_n \leq  K} \lambda_n^{2\sigma } | c_n |^2 |g_\lambda(\omega) |  ^2 \leq \frac{t ^2}{4}  \right\}   \subset   \left\{ \omega \in \Omega / || \ [ u^\omega_0]_K||_{ \overline{H}^\sigma(\R^3) } \leq \frac{t}{2} \right\} ,
\end{align}
\begin{center} et  \end{center}
\begin{equation*}
\left\{  \omega \in \Omega / \sum_{\lambda_n \leq  K} \lambda_n^{2\sigma } | c_n |^2 |g_n(\omega) |  ^2 \leq \frac{t ^2}{C K^6 }  \right\}  \subset A 
\end{equation*}
 \begin{center}
où
\end{center}
\begin{align}\label{esti2}
A = \bigcap & \ \ \left\{ \omega \in \Omega / || \ (e^{itH} [ u^\omega_0]_K)^2 \ ||_{  L^4( [-2\pi,2\pi] , \overline{H}^s (\R^3)) } \leq \frac{t^2}{2} \right\} ,
\\ \label{esti3}
 \bigcap  & \ \ \left\{ \omega \in \Omega / || \ (e^{itH} [ u^\omega_0]_K )^3 \ ||_{  L^4( [-2\pi,2\pi] , \overline{H}^s (\R^3)) } \leq \frac{t^3}{2} \right\} ,
\\ \label{esti4}
  \bigcap & \ \ \underset{N} \bigcap   \left\{   \omega \in \Omega /  || \Delta_N(e^{itH} [ u^\omega_0]_K ) ||_{L^4( [-2\pi,2\pi] , L^\infty (\R^3)) } \leq \frac{t}{2}  N^{-1/6} \right\} ,
\\ \label{esti5}
  \bigcap  & \ \ \underset{N} \bigcap  \left\{  \omega \in \Omega /   || \Delta_N(e^{itH} [ u^\omega_0]_K ) ||_{L^R( [-2\pi,2\pi] , \overline{W}^{s,4} (\R^3)) } \leq \frac{t}{2}  N^{s-1/4}  \right\} .
\end{align}
\end{prop}
\textit{Preuve.} On remarque qu'il existe deux constantes $ C_1,C_2 > 0 $ telles que $ C_1 n^{  \frac{1}{3} } \leq  \lambda^2_n \leq C_2 n^{  \frac{1}{3} } $ puis que $  | \lbrace \lambda_n \leq K  \rbrace | \leq C K^{ 6 }$. Le résultat suit alors grâce à la proposition \ref{dispersive} et à l'inégalité de Cauchy-Schwarz. \hfill $ \boxtimes $
\\\\ \textit{Preuve du théorème \ref{thm1}.} Par le théorème \ref{deviation2} et la proposition \ref{esti6}, il suffit de montrer que pour tout entier  $ K \geq 1 $ et réel $ t > 0 $,
\begin{equation*}
P \left(  \omega \in \Omega / \sum_{\lambda_n \leq  K} \lambda_n^{2\sigma } | c_n |^2 |g_n(\omega) |  ^2 \leq \frac{ t^2}{C K^6 }  \right) > 0.
\end{equation*}
Mais, par indépendance,
\begin{align*}
& P \left(  \omega \in \Omega / \sum_{\lambda_n \leq  K} \lambda_n^{2\sigma } | c_n |^2 |g_n(\omega) |  ^2 \leq \frac{ t^2}{C K^6 }  \right) 
\\ \geq \ &  P \left(  \underset{\lambda_n \leq K }{\bigcap}  \left\{ \omega \in \Omega /  |g_n(\omega) |  ^2 \leq \frac{ t^2}{C K ^{12} ||u_0||^2_{  \overline{H}^\sigma(\R^3)  }  } \right\} \right)
\\ \geq  \ &  \underset{\lambda_n \leq K }{\prod} \ P \left( \omega \in \Omega /   |g_n(\omega) |  ^2 \leq \frac{ t^2}{C K ^{12} ||u_0||^2_{  \overline{H}^\sigma(\R^3)  } } \right) \hspace*{1cm} > 0,
\end{align*}
car pour tout $ R > 0 $ et $ n  \in \N \ , \  P( \omega  \in \Omega / |g_n(\omega)| \leq R  ) > 0$. \hfill $ \boxtimes $
\\\\\textit{Remarque :} Dans le cas Gaussien, comme 
\begin{equation*}
P \left( \omega \in \Omega  / |g_n (\omega ) |^2 \leq R \right) = 1- e^{-R},
\end{equation*}
on en déduit
\begin{equation*}
P \left(  \omega \in \Omega / \sum_{\lambda_n \leq  K} \lambda_n^{2\sigma } | c_n |^2 |g_n(\omega) |  ^2 \leq \frac{ t^2}{C K^2 }  \right) \geq \left( 1- e^{ - \frac{t^2}{C K ^{12} ||u_0||^2_{ \overline{H}^\sigma(\R^3)    }  }   }    \right)^{K^6}.
\end{equation*}
Par conséquent, pour $ \sigma' > \sigma $, en utilisant la remarque du théorème \ref{deviation2}, on obtient qu'il existe deux constantes $  R > 0  $ et $ C > 0 $ telles que si $ ||u_0||_{ \overline{H}^{\sigma'} (\R^3) } \geq R $ alors 
\begin{equation*}
P( \Omega') \geq (1-\alpha)  \times  e^{ -C  ||u_0||^{ \frac{6}{\sigma'-\sigma}   }_{ \overline{H}^{\sigma'} (\R^3) } \log ||u_0||_{ \overline{H}^{\sigma'} (\R^3) }   }.
\end{equation*}
\subsection{Preuve du théorème \ref{thm2}}
En utilisant les théorèmes \ref{existence}, \ref{unicite} et \ref{scattering}, pour prouver (\ref{proba1}), il suffit d'établir que pour tout $ \lambda > 0 $,
\begin{equation}
\lim _{  \eta \rightarrow 0 } \ P \left( \omega \in \Omega_\lambda    | \ ||u^\omega_0||_{\overline{H}^\sigma {(\R^3)}   }   \leq \eta   \right) = 1.
\end{equation}
En adaptant la preuve de l'appendice A.2 de \cite{burq7}, on peut obtenir que 
\begin{equation*}
P \left( \omega \in \Omega^c_\lambda   \ | \ ||u_0^\omega||_{ \overline{H}^\sigma(\R^3) } \leq \eta \right) \leq C e ^{  -c \frac{\lambda^2}{\eta^2} },
\end{equation*}
et (\ref{proba1}) est prouvé.
\\\\Ensuite, pour obtenir (\ref{proba2}), il suffit de montrer que pour tout $ \lambda > 0 $ et $ \epsilon > 0 $,
\begin{align*}
  \lim_{ \eta \rightarrow 0 }  \ \mu_1 \otimes \mu_2 \left(   u_0^1 \in E_0(\lambda) , u_0^2 \in E_0(\lambda) , ||  \tilde{u}_1 - \tilde{u}_2 ||_{X^0} \leq  \epsilon | \right. & \ ||u_0^1 ||_{ \overline{H}^\sigma(\R^3)}  \leq \eta 
  \\ & \hspace*{0.2cm} \left. , ||u_0^2 ||_{ \overline{H}^\sigma(\R^3) } \leq \eta  \right) = 1.
\end{align*}
En utilisant le théorème \ref{continuiteinitial}, il suffit de démontrer que pour tout $ \lambda > 0 $ et $ \epsilon > 0 $,
\begin{align*}
   \lim_{ \eta \rightarrow 0 } \ \mu_1 \otimes \mu_2  \left(   u_0^1 \in E_0(\lambda) , u_0^2 \in E_0(\lambda) 
 | \ ||u_0^1 ||_{ \overline{H}^\sigma(\R^3)}  \leq \eta , ||u_0^2 ||_{ \overline{H}^\sigma(\R^3) } \leq \eta  \right)   = 1.
\end{align*}
Ainsi, il suffit d'établir que pour tout $ \lambda > 0 $ et $ \epsilon > 0 $,
\begin{equation*}
 \lim_{ \eta \rightarrow 0 } \ \mu_1 \otimes \mu_2 \left(   u_0^1 \in E_0(\lambda)^c  | \ ||u_0^1 ||_{ \overline{H}^\sigma(\R^3)}  \leq \eta , ||u_0^2 ||_{ \overline{H}^\sigma(\R^3) } \leq \eta  \right) = 0.
\end{equation*}
Mais
\begin{align*}
 &  \lim_{ \eta \rightarrow 0 } \ \mu_1 \otimes \mu_2 \left(   u_0^1 \in E_0(\lambda)^c  | \ ||u_0^1 ||_{ \overline{H}^\sigma(\R^3)}  \leq \eta , ||u_0^2 ||_{ \overline{H}^\sigma(\R^3) } \leq \eta   \right) 
 \\  = & \lim_{ \eta \rightarrow 0 } \ \mu_1   \left(   u_0^1 \in E_0(\lambda)^c  | \ ||u_0^1 ||_{ \overline{H}^\sigma(\R^3)}  \leq \eta  \right)  = 0,
\end{align*}
et (\ref{proba2}) est démontré.
\subsection{Preuve du théorème \ref{thm3}}
Pour prouver le théorème \ref{thm3}, il suffit de montrer que pour tout $ t> 0 $ et $ \alpha \in ]0,1] $, il existe $ \epsilon > 0 $ tel que $ P( \Omega^\epsilon_t ) \geq 1 - \alpha  $. Ce résultat est clair en vu du théorème \ref{deviation} puisque
\begin{equation*}
P( \Omega^\epsilon_t  ) \geq 1 -  C_1 e ^{ -  \frac{C_2}{ \epsilon^2 ||u_0||^2_{ \overline{H}^{\sigma} (\R^3)  }   }  } \geq 1 - \alpha \mbox{ si } \epsilon \ll 1.
\end{equation*}
\section{Généralisation du résultat}
Dans cette partie, on suppose la dimension d'espace $  d \geq 2 $ et on donne une généralisation du théorème \ref{thm1}. Si l'on suppose que la suite de variables aléatoires $ (g_n)_{n \in \N } $ est identiquement distribuée de mêmes lois gaussiennes complexes standards et que p est un entier impair dans l'équation (\ref{schrodinger}), on peut alors établir le théorème suivant :
\begin{thm} 
Soient $ \sigma \in ]\frac{d-1}{2} - \frac{2}{p-1} , \frac{d}{2} - \frac{2}{p-1} [ $, $ u_0 \in \overline{ H }  ^{\sigma}  ( \R^d ) $ et $ s \in ] \frac{d}{2} - \frac{2}{p-1} ,  \sigma  + \frac{1}{2} [ $ alors il existe un ensemble $ \Omega ' \subset \Omega $ vérifiant les conditions suivantes: 
\\ i) $ P( \Omega ' ) > 0 $,
\\ii) Pour tout élément $  \omega \in \Omega'$, il existe une unique solution globale $ \tilde{u} $ à l'équation (\ref{schrodinger}) dans l'espace $ e^{it\Delta} u_0( \omega,.) + X^s $ avec donnée initiale $ u_0(\omega,.)$.
\\ iii) Pour tout élément $  \omega \in \Omega'$, il existe $ L ^+ $ et $ L_ - \in \overline{H}^s(\R^d) $ telles que
\begin{align*}
 \lim_{t \rightarrow \infty } || \tilde{u}(t) - e^{it\Delta} u_0 ( \omega,.) - e^{it\Delta} L^+ ||_{  H^s(\R^d) } = 0,
\\  \lim_{t \rightarrow -\infty } || \tilde{u}(t) - e^{it\Delta} u_0 ( \omega,.) - e^{it\Delta} L_- ||_{  H^s(\R^d) } = 0.
\end{align*}
\end{thm}
Les points clefs de la démonstration du théorème \ref{thm1} sont l'existence de l'estimée bilinéaire de type Bourgain pour l'oscillateur harmonique et la transformation de lentille, propriétés vraies en dimension quelconque plus grande que 2. Rappelons que les inégalités du chaos de Wiener pour les variables aléatoires gaussiennes sont établies pour une n-linéarité quelconque dans \cite{thomann1}. Ainsi, dans la preuve du théorème \ref{thm1}, le fait que $ p = 3 $ intervient surtout dans l'application du théorème de point fixe de Picard. Vu que les arguments de bases restent vrais en dimension $ d \geq 2$, ce dernier reste applicable pour p quelconque impair à condition de vérifier que $ u_0( \omega ,. ) \in L^\infty( \R, L^\infty(\R^d)) $ $ \omega $ presque surement (si on utilise deux fois l'estimée bilinéaire, les termes restants sont à évaluer dans $ L^\infty_t $, mais pour p=3, il n'y a pas de termes restants), ce que nous expliquons ici. On peut montrer que
\begin{equation*}
 e^{itH} u_0( \omega, .) \in   L^\infty ( [-2 \pi, 2\pi] , \overline{W}^{\frac{1}{6} + \sigma -  , \infty } ( \R^d )) \ \omega \mbox{ presque surement} . 
 \end{equation*}
Grâce à l'inégalité de Minkowsky et les injections de Sobolev, on a pour tout $ \epsilon > 0 $, 
\begin{align*}
|| e^{itH} u_0 ||_{ L^\infty ( [-2\pi,2\pi] , L^p(\R^d) )  }  & \leq  || e^{itH} u_0 ||_{   L^p(\R^d  , L^\infty ( [-2\pi,2\pi] )) } 
\\ & \leq C || e^{itH} u_0 ||_{L^p( \R^d , W^{1/p+\epsilon,p}( [-2\pi,2\pi]  ))}
\\ & \leq C ||H^{1/p+\epsilon} e^{itH} u_0 ||_{L^p ( \R^d , L^p( [-2\pi,2\pi]  ))}.
\end{align*}
Puis nous pouvons remplacer $ u_0 $ par $ H^{  \frac{s}{2} } u _ 0 $ pour obtenir que
\begin{align*}
|| e^{itH} u_0 ||_{ L^\infty ( [-2\pi,2\pi] , \overline{W} ^{s,p}(\R^d) )  } &  \leq C || H^{  \frac{s}{2} + 1/p+ \epsilon } e^{itH} u_0 ||_{ L^p ( \R^d , L^p( [-2\pi,2\pi]  ))}
\\ &  \leq C || H^{  \frac{s}{2} +1/p + \epsilon } e^{itH} u_0 ||_{ L^p( [-2\pi,2\pi] ,L^p(\R^d ))}
\\ & \leq C ||  e^{itH} u_0 ||_{ L^p( [-2\pi,2\pi] , \overline{W} ^{s+2/p+2\epsilon , p}(\R^d ))}.
\end{align*}
Ainsi, si $ \epsilon >  \frac{d}{p} $ alors
\begin{align*}
|| e^{itH} u_0 ||_{ L^\infty ( [-2\pi,2\pi] , \overline{W} ^{\sigma+1/6-4\epsilon,\infty}(\R^d) )  }  & \leq C  || e^{itH} u_0 ||_{ L^\infty ( [-2\pi,2\pi] , \overline{W} ^{\sigma+1/6-3\epsilon,p}(\R^d) )  } 
\\ & \leq C ||  e^{itH} u_0 ||_{ L^p( [-2\pi,2\pi] , \overline{W} ^{\sigma+1/6 , p}(\R^d ))}.
\end{align*}
\nocite{*}
\bibliographystyle{short}
\bibliography{biblioarticle1}
\end{document}